\documentclass[12pt]{article}

\usepackage[latin1]{inputenc}

\usepackage{pgf,tikz}
\usetikzlibrary{arrows}
\usetikzlibrary{patterns}

\usepackage{multirow} 
\usepackage{subfigure}

\usepackage{amssymb}

\usepackage{amsmath,epsfig}
\def\X{\textrm{X}}
\def\Y{\textrm{Y}}

\def\sgn{\textrm{sgn}}

\newtheorem{Theorem}{Theorem}[section]
\newtheorem{Corollary}[Theorem]{Corollary}

\newtheorem{Lemma}[Theorem]{Lemma}
\newtheorem{Example}[Theorem]{Example}

\numberwithin{equation}{section}

\newcommand{\Z}{{\mathbb Z}}

\def\pf{\textrm{Pf}}

\newcommand{\beq}{\begin{equation}}
\newcommand{\eeq}{\end{equation}}

\def\qed{\hfill$\Box$\\ \medskip}

\newcommand{\merge}{{\ast}}

\usepackage{multirow}

\def\mybox{\hbox to 12.0pt}

\def\mybigbox{\hbox to 35.0pt}
\def\myverybigbox{\hbox to 60.0pt}

\def\ov{\overline}
\def\sc{\scriptstyle}
\def\ssc{\scriptscriptstyle}

\def\wgt{{\rm wgt}}

\def\={\!=\!}
\def\+{\!+\!}
\def\-{\!-\!}

\def\x{{\bf x}}
\def\y{{\bf y}}

\def\a{{\bf a}}

\def\0{{\bf 0}}

\def\1{{\bf 1}}

\def\red{\textcolor{red}}
\def\blue{\textcolor{blue}}
\def\magenta{\textcolor{magenta}}
\def\cyan{\textcolor{cyan}}
\def\brown{\textcolor{brown}}
\def\green{\textcolor{green}}

\title{Determinantal and Pfaffian identities for ninth variation skew Schur functions and Q-functions}

\author{
Ang\`ele M. Foley\thanks{ 
Department of Physics and Computer Science,
Wilfrid Laurier University,
Waterloo, Ontario, N2L 3C5, Canada ({\tt ahamel@wlu.ca})}
\and 
Ronald C. King\thanks{
Mathematical Sciences, University of Southampton, 
Southampton SO17 1BJ, England ({\tt r.c.king@soton.ac.uk})}}

\begin{document}

\maketitle

\begin{abstract}
Recently Okada defined algebraically ninth variation skew Q-functions,
in parallel to Macdonald's ninth variation skew Schur functions. Here we introduce
a skew shifted tableaux definition of these ninth variation skew Q-functions, and 
prove by means of a non-intersecting lattice path model a Pfaffian outside 
decomposition result in the form of a ninth variation version of Hamel's Pfaffian 
outside decomposition identity. As corollaries to this we derive Pfaffian identities 
generalizing those of Josefiak-Pragacz, Nimmo, and most recently Okada.
As a preamble to this we present a parallel development based on (unshifted) 
semistandard tableaux that leads to a ninth variation version of the outside 
decomposition determinantal identity of Hamel and Goulden. In this case the 
corollaries we offer include determinantal identities generalizing the Schur and skew 
Schur function identities of Jacobi-Trudi, Giambelli, Lascoux-Pragacz, Stembridge, and 
Okada.
\end{abstract}

\section{Introduction}\label{Sec-introduction}

Macdonald's ninth variation Schur function is the last type of Schur function introduced in his seminal ``Schur Functions: Theme and Variations'' paper \cite{Mac92}.  In some sense it is the ultimate generalization as it includes as special cases five of the other generalizations listed in that 1992 paper.  Close to 30 years later, Okada \cite{Oka2} has introduced a ninth variation $Q$-function, providing an algebraic definition and a number of Pfaffian results in analogy to the classical cases, such as the row-based Pfaffian of Schur \cite{Sch} and the skew row-based Pfaffian of Jozefiak-Pragacz \cite{JP} and Nimmo \cite{Nim}.    We continue to exploit this new type of $Q$-function by creating a combinatorial shifted tableau definition and deriving a general Pfaffian result along the lines of Hamel's outside decomposition result for $Q$-functions \cite{Ham}.

Although the ninth variation $Q$-functions are in their infancy, study of the ninth variation Schur functions has been an active area of research for several years.
Several  authors  (e.g.\ Kim and Yoo \cite{KY}, Nakagawa  {\em et al.} \cite{NNSY}, Furuwara and Moriyama \cite{FM}) have used Macdonald's algebraic definition to derive a host of determinantal results in the style of classical Schur function identities. At the same time a number of other authors (e.g.\ Nakagawa  {\em et al.} \cite{NNSY}, Nakasuji {\em et al.} \cite{NPY}, Bachmann and Charleton \cite{BC}) have employed various tableaux weightings to provide combinatorial results. We are able to extend our tableaux for ninth variation $Q$-functions to the ninth variation Schur function case as well and to do so we introduce a new combinatorial definition that 
involves a 
weighting of the tableau that combines in a very general way its dependence on  the {\em entry} in a tableau box and 
the {\em location} of the box. This weighting unifies a number of the previous weightings which can be shown to be special cases of ours.  The power of our definition lies in the fact that it leads very naturally to a Lindst\"om-Gessel-Viennot lattice path interpretation and from there may be used to prove a wide range of determinantal results. A key feature of our approach is the systematic exploitation of the fact the nonintersecting lattice paths correspond to strips in the tableau diagram that are subdiagrams of the cutting strip defining the relevant outside decomposition, with the box of each strip inheriting its content from its position in the original tableau diagram.

Although our motivation is the ninth variation $Q$-functions, we reverse the exposition to deal with the more basic ninth variation Schur functions first.
Our tableau version of Macdonald's ninth variation of Schur functions~\cite{Mac92},
proceeds by
introducing weight parameters $\X=(x_{kc})$ for $k\in \{1,2,\ldots,n\}$ and $c\in\Z$ so as to define
\begin{equation}\label{eqn-sfnXa}
s_{\lambda/\mu}(\X)= \sum_{T\in{\cal T}^{\lambda/\mu}} \prod_{(i,j)\in F^{\lambda/\mu}} \wgt(t_{ij}),
\end{equation}
where the weight of the entry $t_{ij}$ at position $(i,j)$ is given by $x_{kc}$ with entries $k=t_{ij}$ and content $c=j-i$. (A full explanation of terminology is in Section \ref{Sec-math-background}).

Similarly, in the case of  $Q$-functions, 
 our generalisation in the spirit of Macdonald's ninth variation for $Q$-functions with
weights parametrised by $\X=(x_{kc})$ and $\Y=(y_{kc})$ for entries $k$ or $k'$ with $k\in \{1,2,\ldots,n\}$ and content $c\in\Z$, takes the form
\begin{equation}\label{eqn-QfnXYa}
   Q_{\lambda/\mu}(\X,\!\Y)= \sum_{T\in{\cal ST}^{\lambda/\mu}} \prod_{(i,j)\in SF^{\lambda/\mu}} \wgt(t_{ij}),
\end{equation}
where now the sum is over shifted tableaux and $\wgt(t_{ij})=x_{kc}$ if $t_{ij}=k$ and $\wgt(t_{ij})=y_{kc}$ if $t_{ij}=k'$ with  
$k\in\{1,2,\ldots,n\}$ and   $c=j-i$ in both cases.

We remark that this weighting is highly general, and that a large number of known tableau weightings are simply special cases of this.
Apart from the obvious specialisation of the above tableaux based formulae whereby weights $\wgt(t_{ij})=x_{kc}$ and $y_{kc}$ are replaced by
just $x_k$,
perhaps the most important specialisations are those involving factorial parameters. The definitions of these and their origin may be 
tabulated as follows where $t_{ij}$ is as usual the entry in a box at position $(i,j)$ having content $c=j-i$ and the parameters
$\x=(x_1,x_2,\ldots,x_n)$ and $\y=(y_1,y_2,\ldots,y_n)$ 
are supplemented by factorial parameters $\a=(\ldots,a_{-2},a_{-1},a_0,a_1,a_2,\ldots)$. These functions, along with several of the other special cases, appear in the table below.

\begin{equation}
\begin{array}{|c|c|l|l|}
\hline
\mbox{Entries $t_{ij}$}&\mbox{Weights $\wgt(t_{ij})$}&\mbox{Function}& \mbox{Reference}\cr
\hline
\hline
k&x_k&\mbox{Schur functions}&\cite{Sch,Mac92}\cr
\hline
k&x_k-k-c+1&\mbox{Original factorial Schur functions}&\cite{BL,CL,GH}\cr
\hline
k&x_k+a_{k+c}&\mbox{Factorial Schur functions}&\cite{GG,Mac92}\cr
\hline
k&k^{-a_c}&\mbox{Schur multiple zeta functions}&\cite{NPY,BC}\cr
\hline
k&u_k^{(c)}&\mbox{9th Variation Schur functions}&\cite{NNSY}\cr
\hline
k&f(k,a_c)&\mbox{Diagonal 10th Variation Schur fns}&\cite{BC}\cr
\hline
\hline
\begin{array}{c}k\cr k'\cr\end{array}&\begin{array}{c}x_k\cr x_k\cr\end{array}&\mbox{$Q$-functions}&\cite{Ste89}\cr
\hline
\begin{array}{c}k\cr k'\cr\end{array}&\begin{array}{c}x_k-c\cr x_k+c\cr\end{array}&\mbox{Original factorial $Q$-functions}&\cite{Iva01}\cr
\hline
\begin{array}{c}k\cr k'\cr\end{array}&\begin{array}{c}x_k-a_{c+1}\cr x_k+a_{c+1}\cr\end{array}&\mbox{Factorial $Q$-functions}&\cite{Iva05}\cr
\hline
\begin{array}{c}k\cr k'\cr\end{array}&\begin{array}{c}x_k\cr y_k\cr\end{array}&\mbox{Generalised $Q$-functions}&\cite{HK07}\cr
\hline
\begin{array}{c}k\cr k'\cr\end{array}&\begin{array}{c}x_k+a_{k+c}\cr y_k-a_{k+c}\cr\end{array}&\mbox{Factorial generalised $Q$-functions}&\cite{HK15,FK18b}\cr
\hline
\hline
\end{array}
\label{weights}
\end{equation}

It might be noted that all of these specialisations lead to symmetric functions in the sense that they are invariant under independent 
permutations of the components of $\x$ and $\y$, and those involving both $\x$ and $\y$ are supersymmetric by virtue of the required 
additional property~\cite{Ste85} that if one sets $x_1=t$ and $y_1=-t$ then they are independent of $t$.

We want to further explain that Bachmann and Charlton \cite{BC} define what we term a {\em tenth variation}, going beyond, as it does, the ninth variation. It is equivalent in our context to setting $X=(x_{k,(i,j)})$ 
However, it is only the restricted ``diagonal'' case that leads to determinantal identities, where by {\em diagonal} is meant a weighting dependent diagonally on the content $c=j-i$, rather than independently on both $i$ and $j$. The more general non-diagonal tenth generation weighting, with a weight function $\X=(x_{k,(i,j)})$, does not allow any lattice path interpretation of outside decompositions that could lead in a natural manner to determinantal identities.
This will be made clear at the end of Section \ref{The-HG}.

The paper is organized as follows:  Section \ref{Sec-previous-work} surveys the relevant previous work, including a summary of previous determinantal results for both classical Schur and $Q$--functions.  Section \ref{Sec-math-background} gives the mathematical definitions and results that will be generalized in the following section.
Section \ref{Sec-main-theorems} contains the main theorems for the tableaux version ninth variation Schur functions and ninth variation $Q$-functions.  These theorems are proved using a Lindstr\"om-Gessel-Viennot approach.
Section \ref{Sec-cor} shows the connections to previously known results by demonstrating that they are corollaries of the main theorems in Section \ref{Sec-main-theorems}.

\section{Previous Work}\label{Sec-previous-work}

In working with a tableau definition of ninth variation Schur functions and  $Q$-functions, rather than confining ourselves to the classical Jacobi-Trudi and Giambelli identities, it seems an appropriate time to draw together and to generalise as far as possible the determinantal and Pfaffian 
expressions for all manner of Schur and $Q$-functions.  

To this end, we first recall a general form of determinantal and Pfaffian results for   the classical Schur functions, $s_\lambda(\x)$, and  $Q$-functions, $Q_\lambda(\x)$, and their skew counterparts 
$s_{\lambda/\mu}(\x)$ and $Q_{\lambda/\mu}(\x)$, where $\x=(x_1,x_2,\ldots,x_n)$ is a sequence of indeterminates and $\lambda$ and $\mu$
are partitions. To prove these it is convenient to exploit well known combinatorial definitions of these functions expressed in terms of  
certain tableaux constructed through filling the boxes of skew Young diagrams $F^{\lambda/\mu}$ and their shifted versions $SF^{\lambda/\mu}$
with entries  whose weights are determined by $\x$. Complete definitions of these terms can be found in Section \ref{Sec-math-background}.

In the hands of Hamel and Goulden~\cite{HG} and of Hamel~\cite{Ham} these combinatorial definitions
have given rise, through the use of non-intersecting path models of the relevant tableaux, to the following identities, 
each involving the concept of an outside decomposition and certain strip shapes to be defined later:

\begin{Theorem} [Hamel and Goulden~\cite{HG}] \label{The-HG} 
Let $\lambda$ and $\mu$ be partitions such that $\mu\subseteq\lambda$, and let 
$\theta=(\theta_1,\theta_2,\ldots,\theta_s)$ be an outside decomposition of $F^{\lambda/\mu}$.
Then for all $\x=(x_1,x_2,\ldots,x_n)$
\begin{equation}\label{eqn-HGsfn}
    s_{\lambda/\mu}(\x) = \left|\,  s_{(\theta_p\#\theta_q)}(\x) \,\right|_{1\leq p,q \leq s}\,,
\end{equation}
where $(\theta_p\#\theta_q)$ is either undefined, or a single edge, or a strip
formed by adjoining the strips $\theta_p$ and $\theta_q$, with some overlap if necessary,
while preserving their shape and their content, that is to say 
the diagonals on which their boxes lie.
\end{Theorem}

and 

\begin{Theorem} [Hamel~\cite{Ham}] \label{The-Ham}
Let $\lambda$ and $\mu$ be strict partitions with $\mu\subseteq\lambda$, and let $\theta=(\theta_1,\theta_2,\ldots,\theta_r,\theta_{r+1},\ldots,\theta_s)$ 
be an outside decomposition of $SF^{\lambda/\mu}$ in which $\theta_p$ is strip with a box on the main diagonal if and only if $i\leq r$. 
Then for all $\x=(x_1,x_2,\ldots,x_n)$
\begin{equation}\label{eqn-HpfQ}
    Q_{\lambda/\mu}(\x) = \pf\left( \begin{array}{cc} Q_{(\ov{\theta}_p,\ov{\theta}_q)}(\x)&Q_{(\theta_k\#\theta_p)}(\x)\cr
		                                      -\,^t(Q_{(\theta_k\#\theta_p)}(\x))&0\cr \end{array} \right) 
\end{equation}
with $1\leq p,q\leq s$ and $r+1\leq k\leq s$, where $\ov{\theta}_p$ is a strip produced by extending $\theta_p$ to the 
main diagonal in a manner determined by $\theta$, while $(\theta_k\#\theta_p)$ is as before and 
$(\ov{\theta}_p,\ov{\theta}_q)$ is the skew shape formed by abutting $\ov{\theta}_p$ and $\ov{\theta}_q$ with the 
box on the main diagonal of $\ov{\theta}_p$ lying immediately above and to the left of that of $\ov{\theta}_q$. 
\end{Theorem}

These two theorems contain within them a remarkable number of classical identities, differing
only in the choice of appropriate outside decompositions as made explicit in~\cite{HG} and~\cite{Ham}. 
Determinantal Schur function and skew Schur function identities based on strips are an old idea: 
the Jacobi-Trudi identity~\cite{Jac,Tru} invokes a determinant based on rows; the dual 
Jacob-Trudi identity~\cite{Nag,Kos,Ait}, based on columns; the Giambelli identity~\cite{Gia}, based on hooks; 
its skew version~\cite{LP84}, based on hooks and broken hooks; 
the Lascoux-Pragacz identity~\cite{LP88}, based on rim ribbons. 
All of these have been shown to be corollaries of Theorem~\ref{The-HG} in~\cite{HG}.
We remark that among our corollaries we have what we believe is the first statement of Lascoux and Pragacz's outer rim identity in an explicit form in terms of parts of the partition (see Corollary \ref{Cor-LPX}).

In addition, we show here that a more recent determinantal identity due to Okada~\cite{Oka1} is a special case of
of Theorem~\ref{The-HG} based on a different set of rim ribbons (see Corollary~\ref{Cor-OkaX}),
while the ``shifted" Giambelli identities of Matsuno and Moriyama~\cite{MM}, although not discussed in detail here,
is based on hooks, shifted with respect to those of the Giambelli identity.
The following table shows various known determinantal results for Schur functions:

\begin{equation}
\begin{array}{|l|l|l|}
\hline
\mbox{Deteterminant results}& \mbox{Definition} & \mbox{Reference}\cr
\hline
\hline
\mbox{Jacobi-Trudi}      &  \mbox{algebraic} & \mbox{\cite{Jac}, \cite{Tru}} \cr \hline
\mbox{dual Jacobi-Trudi}   &  \mbox{algebraic} & \mbox{\cite{Nag}}  \cr \hline
\mbox{Giambelli} &  \mbox{algebraic} & \mbox{\cite{Gia}} \cr \hline
 \mbox{skew Giambelli} &  \mbox{algebraic} & \mbox{\cite{LP84}} \cr \hline
 \mbox{outer rim ribbon} &  \mbox{algebraic} & \mbox{\cite{LP88}} \cr \hline
\mbox{Jacobi-Trudi} &  \mbox{tableaux} &  \mbox{\cite{GV}}\cr \hline
 \mbox{dual Jacobi-Trudi} & \mbox{tableaux} & \mbox{\cite{Uen}} \cr \hline
\mbox{skew Giambelli}&  \mbox{tableaux} & \mbox{\cite{Ste90}} \cr \hline
\mbox{outside decomposition} & \mbox{tableaux} & \mbox{\cite{HG}} \cr 
\mbox{({\em which includes all others in this table}) }  & &  \cr \hline
  \mbox{``shifted'' hook }   &   \mbox{tableaux  }  & \mbox{  \cite{MM} } \cr \hline
  \mbox{inner rim ribbon} & \mbox{tableaux} & \mbox {\cite{Oka1}}  \cr \hline 
  \hline
\end{array}
\end{equation}

Pfaffians play the role of determinants in a parallel theory of $Q$-functions and skew $Q$-functions.
Indeed $Q$-functions were originally defined as Pfaffians~\cite{Sch} 
and the generalisation to the skew case is a Pfaffian identity due to J\'ozefiak and Pragacz~\cite{JP} and Nimmo \cite{Nim}  that may be 
viewed as a special case of Theorem~\ref{The-Ham} based on strips corresponding to rows~\cite{Ham}.
It will be shown that the recent Pfaffian skew $Q$-function identity due to Okada~\cite{Oka1} is
another special case of Theorem~\ref{The-Ham} based this time on rim ribbon strips, and as we shall see there
is yet another new identity of the same type based on a different set of rim ribbons.
 The following table shows  various known Pfaffian results for $Q$-functions:

  \begin{equation}
\begin{array}{|l|l|l|}
\hline
\mbox{Pfaffian results}&\ \mbox{Definition} & \mbox{Reference}\cr
\hline
\hline
 \mbox{row-based Pfaffian} &  \mbox{algebraic} & \mbox{\cite{Sch}}  \cr \hline
\mbox{skew row-based Pfaffian}& \mbox{algebraic} &  \mbox{\cite{JP}, \cite{Nim}} \cr \hline
 \mbox{skew row-based Pfaffian}  &  \mbox{shifted tableaux} &  \mbox{\cite{Ste89}} \cr \hline
\mbox{outside decomposition}& \mbox{shifted tableaux} & \mbox{\cite{Ham}} \cr 
\mbox{({\em which includes all others in this table}) } &  &  \cr \hline
\hline
\end{array}
\end{equation}

All of these determinantal and Pfaffian identities will be generalised in Theorems \ref{The-HGX} and \ref{The-HamXY} of Section \ref{Sec-main-theorems} in a quite different way. 
As discussed in Section \ref{Sec-introduction}, these generalisations come about, not by varying the choice of outside decomposition, but by varying the weight assigned to each entry 
in the semistandard and the primed shifted tableaux that are used in the combinatorial definitions of skew Schur functions 
and skew $Q$-functions, respectively. For this purpose it should be noted that a box in the $i$th row and $j$th
column of either type of tableau is said to have content $c=j-i$, which is not to be confused with its entry $k$ or
$k'$ for some $k\in\{1,2,\ldots,n\}$.
In exactly the same way the Pfaffian identities of Theorem~\ref{The-Ham} which involve weights $x_k$ for entries $k$ 
or $k'$ in the primed shifted tableaux that define $Q_{\lambda/\mu}(\x)$ will be extended here to
Pfaffian identities for $Q_{\lambda/\mu}(\X,\!\Y)$ with the weight of an entry $k$ or $k'$ in a box of content $c$ being given by 
or $x_{kc}$ or $y_{kc}$, respectively. This is true of all possible outside decompositions and for all 
$\X=(x_{kc})$ and $\Y=(y_{kc})$ with $k\in\{1,2,\ldots,n\}$ and $c\in\Z$. 
The table (\ref{weights}) in Section \ref{Sec-introduction} shows the various known weightings.

There are several other papers that concern determinantal identities for ninth variation Schur or $Q$-functions and it is  important to distinguish our contribution from those. The following table illustrates the differences between the various results.

\begin{equation}
\begin{array}{|l|c|l|l|}
\hline
\mbox{Determinant and Pfaffian results}&\mbox{Function}& \mbox{Definition} & \mbox{Reference}\cr
\hline
\hline
\mbox{skew JT, skew Dual JT,}   &\mbox{ 9th variation }&\mbox{tableaux} & \mbox{\cite{NNSY} }\cr  
\mbox{Giambelli---all proved algebraically} & \mbox{Schur functions}& & \cr \hline
 \mbox{skew JT, skew Dual JT, Giambelli   } &   \mbox{multiple zeta  }   & \mbox{tableaux   }  &  \mbox{\cite{NPY}  }\cr
   & \mbox{ functions}& & \cr \hline
  \mbox{Hamel-Goulden determinants} &   \mbox{ diagonal 10th}   & \mbox{tableaux   }  &  \mbox{\cite{BC} }\cr
  \mbox{(i.e.\ arbitrary outside decompositions), } &\mbox{variation} & & \cr
  \mbox{but only when $k_{ij}$ depends on $c=j-i$.}&  \mbox{Schur functions}& & \cr \hline
   \mbox{Hamel-Goulden determinants }&   \mbox{  9th variation  }   & \mbox{algebraic  }  &   \mbox{\cite{KY}   } \cr 
    \mbox{(arbitrary outside decompositions )  } & \mbox{Schur functions}& & \cr \hline
 \mbox{JT, Giambelli and some  }    &   \mbox{9th variation  }   & \mbox{algebraic  }  &    \mbox{\cite{FM} } \cr 
 \mbox{variations thereof}   & \mbox{Schur functions} & &    \cr \hline 
    \mbox{Hamel-Goulden determinants}   &   \mbox{9th variation  }   & \mbox{tableaux   }  & \mbox{This paper } \cr 
      \mbox{(arbitrary outside decompositions),} &  \mbox{Schur functions} & & \cr
    \mbox{also includes Okada inner rim ribbon   }  &  & & \cr \hline \hline
     \mbox{Jozefiak-Pragacz-Nimmo   }  &   \mbox{9th variation  }   & \mbox{algebraic }  & \mbox{\cite{Oka2}   }  \cr 
         & \mbox{$Q$ functions}& & \cr \hline
   \mbox{Hamel Pfaffians}     &   \mbox{9th variation }    & \mbox{tableaux   }  &  \mbox{This paper } \cr  
  \mbox{(i.e.\ arbitrary outside decompositions),} & \mbox{ generalised} & & \cr
    \mbox{also includes JPN}  & \mbox{$Q$-functions }& & \cr \hline
\hline
\end{array}
\end{equation}

\section{Mathematical Background}~\label{Sec-math-background}
\subsection{Schur function and $Q$-function definitions}~\label{Sec-sfn-Qfun}

Since this work takes place inside the world of partitions, tableaux, outside decompositions, strips and cutting strips,
some definitions are in order. First we give those pertaining to the definitions of Schur functions and $Q$-functions. 

Let $\lambda$ be a partition of length $\ell(\lambda)=l$ so that $\lambda=(\lambda_1,\lambda_2,\ldots,\lambda_l)$ 
where each part $\lambda_i$ is a positive integer and $\lambda_1\geq\lambda_2\geq\cdots\geq\lambda_l>0$. 
Each such partition $\lambda$ specifies a Young diagram $F^\lambda$ of shape $\lambda$
consisting from top to bottom of consecutive rows of boxes of lengths $\lambda_i$ for $i=1,2,\ldots,l$ left adjusted to a vertical line. 
The conjugate partition $\lambda'=(\lambda'_1,\lambda'_2,\ldots,\lambda'_{l'})$ of length $\ell(\lambda')=l'=\lambda_1$
has parts $\lambda'_j$ for $j=1,2,\ldots,l'$ equal to the lengths from left to right of consecutive columns of $F^\lambda$.
 
The box of $F^\lambda$ in the $i$th row and $j$th column is denoted by $(i,j)$ and is said
to have content $j-i$.  The {\em diagonal} with content $c$ is the set of all boxes in $F^\lambda$ having content $c$.
Boxes on the main diagonal have content $0$, while those above and below have positive and negative content, respectively.
We shall adopt a typographically convenient notation whereby any negative content $c=-a$ is written as $\ov{a}$.

If the number of boxes $(i,i)$ on the main diagonal of $F^\lambda$ is $p$ then 
in Frobenius notation, we write $\lambda=(\alpha|\beta)$, with $\alpha=(\alpha_1,\alpha_2,\ldots,\alpha_p)$ 
and $\beta=(\beta_1,\beta_2,\ldots,\beta_p)$, where $\alpha_i=\lambda_i-i$ and $\beta_i=\lambda'_i-i$
for $i=1,2,\ldots,p$ are the arm and leg lengths of $F^\lambda$ with respect to its main diagonal. 
It follows that $\alpha$ and $\beta$ are strict partitions in the sense that
$\alpha_1>\alpha_2>\cdots>\alpha_p\geq0$ and $\beta_1>\beta_2>\cdots>\beta_p\geq0$.
The parameter $p$ is sometimes called the Durfee size or the Frobenius rank of $\lambda$.

Given a second partition $\mu=(\mu_1,\mu_2,\ldots,\mu_m)$ of length $\ell(\mu)=m$ such that $F^\mu\subseteq F^\lambda$,
that is $\mu_i\leq\lambda_i$ for all $i=1,2,\ldots,m$, we write $\mu\subseteq\lambda$ and the pair $\lambda$ and $\mu$ 
specify a skew Young diagram $F^{\lambda/\mu}$ of shape $\lambda/\mu$ obtained by deleting from $F^\lambda$ all 
the boxes of $F^\mu$. In Frobenius notation, it follows in such a case 
that if $\mu=(\gamma|\delta)$ has Frobenius rank $q$ with $\gamma=(\gamma_1,\gamma_2,\ldots,\gamma_q)$ and $\delta=(\delta_1,\delta_2,\ldots,\delta_q)$, 
then $q\leq p$ and both $\gamma_i\leq\alpha_i$ and $\delta_i\leq\beta_i$ for all $i=1,2,\ldots,q$. The content
of each box of a skew Young diagram $F^{\lambda/\mu}$ is that of the corresponding box in $F^\lambda$. However it 
should be noted that even if $F^{\lambda/\mu}$ and $F^{\kappa/\nu}$ are of the same skew shape their content may be 
different, but may be shifted with respect to one another by the same amount for each box.

Now we define the tableaux involved in the definition of Schur functions and their generalisations.
Given a skew Young diagram $F^{\lambda/\mu}$, let ${\cal T}^{\lambda/\mu}$ be the set of all semistandard skew 
Young tableaux $T$ obtained by filling each box of $F^{\lambda/\mu}$ 
with an entry from the alphabet $\{1<2<\cdots<n\}$ in such a way that the entries are weakly increasing across rows from left to right 
and strictly increasing down columns from top to bottom 

For $\x=(x_1,x_2,\ldots,x_n)$ the Schur function $s_{\lambda/\mu}(\x)$ may then be defined by
\begin{equation}\label{eqn-sfn}
   s_{\lambda/\mu}(\x)= \sum_{T\in{\cal T}^{\lambda/\mu}} \prod_{(i,j)\in F^{\lambda/\mu}} \wgt(t_{ij}), 
\end{equation}
where $t_{ij}$ is the entry of $T$ at position $(i,j)$ and $\wgt(t_{ij})=x_k$ if $t_{ij}=k$ for any $k\in\{1,2,\ldots,n\}$.

Recall from Section \ref{Sec-introduction}  
our tableaux version of Macdonald's ninth variation of Schur functions~\cite{Mac92},
differing somewhat from that of Nakagawa {\it et al}~\cite{NNSY}, is defined by
introducing weight parameters $\X=(x_{kc})$ for $k\in \{1,2,\ldots,n\}$ and $c\in\Z$ so as to produce
\begin{equation}\label{eqn-sfnX}
   s_{\lambda/\mu}(\X)= \sum_{T\in{\cal T}^{\lambda/\mu}} \prod_{(i,j)\in F^{\lambda/\mu}} \wgt(t_{ij}),
\end{equation}
where now $\wgt(t_{ij})=x_{kc}$ with $k=t_{ij}$ and $c=j-i$.

A partition $\lambda=(\lambda_1,\lambda_2,\ldots,\lambda_m)$ is said to be strict if $\lambda_1>\lambda_2>\cdots>\lambda_m\geq0$. Then,
in an entirely analogous manner, for any pair of strict partitions $\lambda$ and $\mu$ such that $\mu\subseteq\lambda$, let ${\cal ST}^{\lambda/\mu}$ be the set of all primed shifted  tableaux $T$ obtained by filling each box of the skew shifted Young diagam $SF^{\lambda/\mu}$ with an entry from the alphabet $\{1'<1<2'<2<\cdots<n'<n\}$ in such a way that the entries are weakly increasing across rows from left to right and down columns from top to bottom, with no two identical unprimed entries $k$ in any column, and no two identical primed entries $k'$ in any row. 

 The $Q$-function $Q_{\lambda/\mu}(\x)$ may then be defined by
\begin{equation}\label{eqn-Qfn}
   Q_{\lambda/\mu}(\x)= \sum_{T\in{\cal ST}^{\lambda/\mu}} \prod_{(i,j)\in SF^{\lambda/\mu}} \wgt(t_{ij}), 
\end{equation}
where $t_{ij}$ is the entry of $T$ at position $(i,j)$ and $\wgt(t_{ij})=x_k$ if $t_{ij}=k$ or $k'$
for any $k\in\{1,2,\ldots,n\}$.

Recall again from Section \ref{Sec-introduction} our generalisation in the spirit of Macdonald's ninth variation, but this time for $Q$-functions with
weights parametrised by $\X=(x_{kc})$ and $\Y=(y_{kc})$ for $k\in \{1,2,\ldots,n\}$ and $c\in\Z$, takes the form
\begin{equation}\label{eqn-QfnXY}
   Q_{\lambda/\mu}(\X,\!\Y)= \sum_{T\in{\cal ST}^{\lambda/\mu}} \prod_{(i,j)\in SF^{\lambda/\mu}} \wgt(t_{ij}), 
\end{equation}
where now $\wgt(t_{ij})=x_{kc}$ if $t_{ij}=k$ and $\wgt(t_{ij})=y_{kc}$ if $t_{ij}=k'$ with  
$k\in\{1,2,\ldots,n\}$ and $c=j-i$ in both cases.

\subsection{Outside decompositions and cutting strips for skew diagrams}~\label{Sec-outside-decompositions}

A {\em strip} is a skew Young diagram which is edgewise connected and contains no $2 \times 2$ block 
of boxes, that is to say, it has a single box on each of its diagonals. 
The bottommost and leftmost box in a strip $\theta$ is its {\em starting box} and the topmost and rightmost box is its {\em ending box}. 
These are the boxes of lowest and highest content, respectively.  
It is possible (typically in many different ways) to decompose any skew
Young diagram into a disjoint union of strips. If each constituent strip $\theta_p$ starts 
on the left or bottom edge of the diagram and finishes 
on the top or right edge of the diagram, the union $\Theta=(\theta_1,\theta_2,\ldots,\theta_s)$ 
is known as an {\em outside decomposition}, a notion first introduced and exploited in~\cite{HG}.
The content of each box of each strip $\theta_p$ in an outside decomposition of $F^{\lambda/\mu}$ is inherited from 
that of its position in $F^{\lambda/\mu}$. 

As pointed out by Chen et al.~\cite{CYY}, a convenient way of generating outside decompositions exploits the notion of a 
{\em cutting strip}. A cutting strip, $\phi$, 
associated with a skew Young diagram $F^{\lambda/\mu}$ is any strip whose boxes have content ranging
over that of all boxes of $F^{\lambda/\mu}$. 
A cutting strip is not necessarily a subdiagram of $F^{\lambda/\mu}$ or even of $F^{\lambda}$. However, each cutting strip $\phi$ may be used 
to define an associated outside decomposition $\Theta=(\theta_1,\theta_2,\ldots,\theta_s)$ of $F^{\lambda/\mu}$ 
by superposing diagonally translated copies of $\phi$ onto $F^{\lambda/\mu}$,
thereby dividing it into a mutually disjoint set of strips. Each strip $\theta_p$ is thus some connected component of the intersection of one 
of these copies of $\phi$ with $F^{\lambda/\mu}$ and shares exactly the same 
shape and content as some portion, say $\phi_{ab}$, of $\phi$ that extends from a box of content $a$ to one of content $b$ with $a\leq b+1$. 
In such a case we write $\theta_p\simeq\phi_{ab}$.

In the case $\lambda=(5,4,4,2)$ and $\mu=(3,2)$ this is all illustrated below
for one particular choice of cutting strip $\phi$:
\begin{equation}\label{eqn-cs}
\begin{array}{cc}
\mbox{Content of $F^{\lambda/\mu}$}&\mbox{Content of cutting strip $\phi$}\cr
\vcenter{\hbox{
\begin{tikzpicture}[x={(0in,-0.2in)},y={(0.2in,0in)}] 
\foreach \j in {1,...,5} \draw[thick] (1,\j) rectangle +(-1,-1);
\foreach \j in {1,...,4} \draw[thick] (2,\j) rectangle +(-1,-1);
\foreach \j in {1,...,4} \draw[thick] (3,\j) rectangle +(-1,-1);
\foreach \j in {1,...,2} \draw[thick] (4,\j) rectangle +(-1,-1);
\foreach \j in {1,...,3} \draw(1-0.5,\j-0.5)node{$\ast$};
\foreach \j in {1,...,2} \draw(2-0.5,\j-0.5)node{$\ast$};
\draw(1-0.5,4-0.5) node {$3$};\draw(1-0.5,5-0.5) node {$4$};
\draw(2-0.5,3-0.5) node {$1$};\draw(2-0.5,4-0.5) node {$2$};
\draw(3-0.5,1-0.5) node {$\ov2$};\draw(3-0.5,2-0.5) node {$\ov1$};\draw(3-0.5,3-0.5) node {$0$};\draw(3-0.5,4-0.5) node {$1$};
\draw(4-0.5,1-0.5) node {$\ov3$};\draw(4-0.5,2-0.5) node {$\ov2$};
\end{tikzpicture}
}}
&
\vcenter{\hbox{
\begin{tikzpicture}[x={(0in,-0.2in)},y={(0.2in,0in)}] 
\foreach \j in {5,6} \draw[thick] (1,\j) rectangle +(-1,-1);
\foreach \j in {2,...,5} \draw[thick] (2,\j) rectangle +(-1,-1);
\foreach \j in {1,...,2} \draw[thick] (3,\j) rectangle +(-1,-1);
\draw(1-0.5,5-0.5) node {$\green{3}$}; \draw(1-0.5,6-0.5) node {$\brown{4}$};
\draw(2-0.5,2-0.5) node {$\green{\ov1}$}; \draw(2-0.5,3-0.5) node {$\brown{0}$}; \draw(2-0.5,4-0.5) node {$\brown{1}$}; \draw(2-0.5,5-0.5) node {$\brown{2}$};
\draw(3-0.5,1-0.5) node {$\green{\ov3}$}; \draw(3-0.5,2-0.5) node {$\brown{\ov2}$};
\end{tikzpicture}
}}
\cr\cr
\mbox{Superposition of cutting strips}&\mbox{Outside decomposition}\cr\cr
\vcenter{\hbox{
\begin{tikzpicture}[x={(0in,-0.2in)},y={(0.2in,0in)}] 
\foreach \j in {1,...,5} \draw[thick] (1,\j) rectangle +(-1,-1);
\foreach \j in {1,...,4} \draw[thick] (2,\j) rectangle +(-1,-1);
\foreach \j in {1,...,4} \draw[thick] (3,\j) rectangle +(-1,-1);
\foreach \j in {1,...,2} \draw[thick] (4,\j) rectangle +(-1,-1);
\fill[blue!40!white](4.5,0.3)--(3.3,0.3)--(3.3,1.3)--(2.3,1.3)--(2.3,4.3)--(1.3,4.3)--(1.3,6.5)
--(1.7,6.5)--(1.7,4.7)--(2.7,4.7)--(2.7,1.7)--(3.7,1.7)--(3.7,0.7)--(4.5,0.7)--cycle;
\fill[red!40!white](3.5,-0.7)--(2.3,-0.7)--(2.3,0.3)--(1.3,0.3)--(1.3,3.3)--(0.3,3.3)--(0.3,5.5)
--(0.7,5.5)--(0.7,3.7)--(1.7,3.7)--(1.7,0.7)--(2.7,0.7)--(2.7,-0.3)--(3.5,-0.3)--cycle;
\foreach \j in {1,...,3} \draw(1-0.5,\j-0.5)node{$\ast$};
\foreach \j in {1,...,2} \draw(2-0.5,\j-0.5)node{$\ast$};
\end{tikzpicture}
}}
&
\vcenter{\hbox{
\begin{tikzpicture}[x={(0in,-0.2in)},y={(0.2in,0in)}] 
\foreach \j in {1,...,5} \draw[thick] (1,\j) rectangle +(-1,-1);
\foreach \j in {1,...,4} \draw[thick] (2,\j) rectangle +(-1,-1);
\foreach \j in {1,...,4} \draw[thick] (3,\j) rectangle +(-1,-1);
\foreach \j in {1,...,2} \draw[thick] (4,\j) rectangle +(-1,-1);
\foreach \j in {1,...,3} \draw(1-0.5,\j-0.5)node{$\ast$};
\foreach \j in {1,...,2} \draw(2-0.5,\j-0.5)node{$\ast$};
\draw(4.7,0.6)node{$\blue{\theta_1}$};
\draw(2.5,-0.6)node{$\red{\theta_2}$};
\draw(1.5,1.4)node[fill=white,inner sep=1pt]{$\cyan{\theta_3}$};
\draw[blue,ultra thick](4.2,0.5)--(3.5,0.5)--(3.5,1.5)--(2.5,1.5)--(2.5,4.2);
\draw[red,ultra thick](2.5,-0.2)--(2.5,0.5)--(1.8,0.5);
\draw[cyan,ultra thick](1.5,1.8)--(1.5,3.5)--(0.5,3.5)--(0.5,5.2);
\end{tikzpicture}
}}
\end{array}
\end{equation}
It follows in this example that the strips $\theta_p$ and their content are given, along with their identification 
with portions of $\phi$, by:
\begin{equation}\label{eqn-strips}
\theta_1=
\vcenter{\hbox{
\begin{tikzpicture}[x={(0in,-0.2in)},y={(0.2in,0in)}] 
\foreach \j in {2,...,4} \draw[thick] (3,\j) rectangle +(-1,-1);
\foreach \j in {1,...,2} \draw[thick] (4,\j) rectangle +(-1,-1);
\draw(4-0.5,1-0.5) node {$\blue{\ov3}$}; \draw(4-0.5,2-0.5) node {$\blue{\ov2}$};
\draw(3-0.5,2-0.5) node {$\blue{\ov1}$}; \draw(3-0.5,3-0.5) node {$\blue{0}$}; 
\draw(3-0.5,4-0.5) node {$\blue{1}$};
\end{tikzpicture}
}}\simeq \phi_{\ov31};\qquad
\theta_2=
\vcenter{\hbox{
\begin{tikzpicture}[x={(0in,-0.2in)},y={(0.2in,0in)}] 
\foreach \j in {1,...,1} \draw[thick] (3,\j) rectangle +(-1,-1);
\draw(3-0.5,1-0.5) node {$\red{\ov2}$};
\end{tikzpicture}
}}\simeq \phi_{\ov2\ov2};\qquad
\theta_3=
\vcenter{\hbox{
\begin{tikzpicture}[x={(0in,-0.2in)},y={(0.2in,0in)}] 
\foreach \j in {4,...,5} \draw[thick] (1,\j) rectangle +(-1,-1);
\foreach \j in {3,...,4} \draw[thick] (2,\j) rectangle +(-1,-1);
\draw(2-0.5,3-0.5) node {$\cyan{1}$}; \draw(2-0.5,4-0.5) node {$\cyan{2}$};
\draw(1-0.5,4-0.5) node {$\cyan{3}$}; \draw(1-0.5,5-0.5) node {$\cyan{4}$};
\end{tikzpicture}
}}\simeq \phi_{14}
\end{equation}

It should be added that an arbitrary outside decomposition $\Theta$ of $F^{\lambda/\mu}$ may also contain 
null strips~\cite{HG}. They each consist of a single edge through which passes one of the diagonal translations 
of $\phi$. If this edge constituting a null strip $\theta_p$ lies between boxes of content $a-1$ and $a$, 
then we make the identification $\theta_p\simeq\phi_{a,a-1}$. 

We can combine any two strips $\theta_p\simeq\phi_{ab}$ and $\theta_q\simeq\phi_{cd}$ from an outside decomposition to create a further strip. 
This is done via the $\#$ operation of Hamel and Goulden~\cite{HG}, which can be defined so that 
\begin{equation}\label{eqn-hash}
(\theta_p\#\theta_q)=\phi_{ad} \quad\mbox{where}\quad  \phi_{ad} \begin{cases} 
                        \mbox{is well defined for $a\leq d$};\cr
                        \mbox{is defined to be a null strip if $a=d+1$};\cr
                        \mbox{is empty if $a>d+1$}.\cr
                             \end{cases}
\end{equation} 

By way of illustration in the case of the example in (\ref{eqn-cs}) with outside decomposition strips as in (\ref{eqn-strips}) 
the relevant combinations under the $\#$ operation products can be tabulated as follows:
\begin{equation}
\begin{array}{|c|ccc|}
\hline
(\theta_p\#\theta_q)&\theta_1\simeq\phi_{\ov31}&\theta_2\simeq\phi_{\ov2\ov2}&\theta_3\simeq\phi_{14}\cr
\hline
\theta_1\simeq\phi_{\ov31}    &\phi_{\ov31}&\phi_{\ov3\ov2}&\phi_{\ov34}\cr
\theta_2\simeq\phi_{\ov2\ov2}    &\phi_{\ov21}&\phi_{\ov2\ov2}&\phi_{\ov24}\cr
\theta_3\simeq\phi_{14}  &\phi_{11}&\phi_{1\ov2}&\phi_{14}\cr
\hline
\end{array}
\end{equation}
or, more pictorially,
\begin{equation}\label{eqn-pict}
\begin{array}{|c|ccc|}
\hline
(\theta_p\#\theta_q)&\theta_1&\theta_2&\theta_3\cr
\hline
&&&\cr
\theta_1&
\vcenter{\hbox{
\begin{tikzpicture}[x={(0in,-0.2in)},y={(0.2in,0in)}] 
\foreach \j in {2,...,4} \draw[thick] (2,\j) rectangle +(-1,-1);
\foreach \j in {1,...,2} \draw[thick] (3,\j) rectangle +(-1,-1);
\draw(3-0.5,1-0.5) node {$\green{\ov3}$}; \draw(3-0.5,2-0.5) node {$\brown{\ov2}$};
\draw(2-0.5,2-0.5) node {$\green{\ov1}$}; \draw(2-0.5,3-0.5) node {$\brown{0}$}; \draw(2-0.5,4-0.5) node {$\brown{1}$}; 
\end{tikzpicture}
}}   
&
\vcenter{\hbox{
\begin{tikzpicture}[x={(0in,-0.2in)},y={(0.2in,0in)}] 
\foreach \j in {1,...,2} \draw[thick] (3,\j) rectangle +(-1,-1);
\draw(3-0.5,1-0.5) node {$\green{\ov3}$}; \draw(3-0.5,2-0.5) node {$\brown{\ov2}$};
\end{tikzpicture}
}}   
&
\vcenter{\hbox{
\begin{tikzpicture}[x={(0in,-0.2in)},y={(0.2in,0in)}] 
\foreach \j in {5,6} \draw[thick] (1,\j) rectangle +(-1,-1);
\foreach \j in {2,...,5} \draw[thick] (2,\j) rectangle +(-1,-1);
\foreach \j in {1,...,2} \draw[thick] (3,\j) rectangle +(-1,-1);
\draw(3-0.5,1-0.5) node {$\green{\ov3}$}; \draw(3-0.5,2-0.5) node {$\brown{\ov2}$};
\draw(2-0.5,2-0.5) node {$\green{\ov1}$}; \draw(2-0.5,3-0.5) node {$\brown{0}$}; 
\draw(2-0.5,4-0.5) node {$\brown{1}$}; \draw(2-0.5,5-0.5) node {$\brown{2}$};
\draw(1-0.5,5-0.5) node {$\green{3}$}; \draw(1-0.5,6-0.5) node {$\brown{4}$};
\end{tikzpicture}
}}\cr
\theta_2&
\vcenter{\hbox{
\begin{tikzpicture}[x={(0in,-0.2in)},y={(0.2in,0in)}] 
\foreach \j in {2,...,4} \draw[thick] (2,\j) rectangle +(-1,-1);
\foreach \j in {2,...,2} \draw[thick] (3,\j) rectangle +(-1,-1);
\draw(3-0.5,2-0.5) node {$\brown{\ov2}$};
\draw(2-0.5,2-0.5) node {$\green{\ov1}$}; \draw(2-0.5,3-0.5) node {$\brown{0}$}; 
\draw(2-0.5,4-0.5) node {$\brown{1}$}; 
\end{tikzpicture}
}}   
&
\vcenter{\hbox{
\begin{tikzpicture}[x={(0in,-0.2in)},y={(0.2in,0in)}] 
\foreach \j in {2,...,2} \draw[thick] (3,\j) rectangle +(-1,-1);
\draw(3-0.5,2-0.5) node {$\brown{\ov2}$};
\end{tikzpicture}
}}   
&
\vcenter{\hbox{
\begin{tikzpicture}[x={(0in,-0.2in)},y={(0.2in,0in)}] 
\foreach \j in {5,6} \draw[thick] (1,\j) rectangle +(-1,-1);
\foreach \j in {2,...,5} \draw[thick] (2,\j) rectangle +(-1,-1);
\foreach \j in {2,...,2} \draw[thick] (3,\j) rectangle +(-1,-1);
\draw(3-0.5,2-0.5) node {$\brown{\ov2}$};
\draw(2-0.5,2-0.5) node {$\green{\ov1}$}; \draw(2-0.5,3-0.5) node {$\brown{0}$}; 
\draw(2-0.5,4-0.5) node {$\brown{1}$}; \draw(2-0.5,5-0.5) node {$\brown{2}$};
\draw(1-0.5,5-0.5) node {$\green{3}$}; \draw(1-0.5,6-0.5) node {$\brown{4}$};
\end{tikzpicture}
}}\cr
\theta_3&
\vcenter{\hbox{
\begin{tikzpicture}[x={(0in,-0.2in)},y={(0.2in,0in)}] 
\foreach \j in {4,...,4} \draw[thick] (2,\j) rectangle +(-1,-1);
\draw(2-0.5,4-0.5) node {$\brown{1}$}; 
\end{tikzpicture}
}}   
&
\emptyset
&
\vcenter{\hbox{
\begin{tikzpicture}[x={(0in,-0.2in)},y={(0.2in,0in)}] 
\foreach \j in {5,6} \draw[thick] (1,\j) rectangle +(-1,-1);
\foreach \j in {4,...,5} \draw[thick] (2,\j) rectangle +(-1,-1);
\draw(2-0.5,4-0.5) node {$\brown{1}$}; \draw(2-0.5,5-0.5) node {$\brown{2}$};
\draw(1-0.5,5-0.5) node {$\green{3}$}; \draw(1-0.5,6-0.5) node {$\brown{4}$};
\end{tikzpicture}
}} \cr 
&&&\cr
\hline
\end{array}
\end{equation}

There remains one more definition, namely that of a shift parameter.
By construction each diagram $(\theta_p\#\theta_q)\simeq\phi_{ad}$ is either a strip, a single edge 
(not illustrated in this example, but see~\cite{HG}), or empty. In the case of a strip, $\phi_{ad}$
is a skew diagram of the form $F^{\nu/\kappa}$ for some partitions $\nu$ and $\kappa$. These are not uniquely
determined. However, this may be rectified by choosing $\nu$ to be the unique partition such that
the outer rim of $F^\nu$ has shape $\phi_{ad}$. 
It will be seen that content of the boxes of $(\theta_p\#\theta_q)\simeq\phi_{ad}$
do not coincide with the conventionally defined content of the same boxes of $F^{\nu/\kappa}$, but instead they are shifted
by a fixed amount that we denote by $m(p,q)$ and refer to as the shift parameter of $(\theta_p\#\theta_q)$.
This is the value of the content of the box of $\phi_{ad}$ that lies at its intersection with the main diagonal of $F^\nu$.

It can be seen from (\ref{eqn-pict}) that in the case of this example the shift parameters are given by:
\begin{equation}\label{eqn-mpq}
\begin{array}{|c|ccc|}
\hline
m(p,q)&\theta_1&\theta_2&\theta_3\cr
\hline
&&&\cr
\theta_1&\ov2&\ov3&\ov1\cr
\theta_2&\ov1&\ov2&0\cr
\theta_3&1&-&2\cr
\hline
\end{array}
\end{equation}

\subsection{Outside decompositions and cutting strips for skew shifted diagrams}\label{subsec-shifted}
Precisely the same use of a cutting strip may be applied to the construction of outside decompositions of skew 
Young diagrams $SF^{\lambda/\mu}$. This is illustrated below in the case $\lambda=(9,6,4,2)$ and $\mu=(4,3)$
for one particular choice of cutting strip.
\begin{equation}\label{eqn-csQ}
\begin{array}{cc}
\mbox{Content of $SF^{\lambda/\mu}$}&\mbox{Cutting strip $\phi$}
\cr\cr
\vcenter{\hbox{
\begin{tikzpicture}[x={(0in,-0.2in)},y={(0.2in,0in)}] 
\foreach \j in {1,...,9} \draw[thick] (1,\j) rectangle +(-1,-1);
\foreach \j in {2,...,7} \draw[thick] (2,\j) rectangle +(-1,-1);
\foreach \j in {3,...,6} \draw[thick] (3,\j) rectangle +(-1,-1);
\foreach \j in {4,...,5} \draw[thick] (4,\j) rectangle +(-1,-1);
\foreach \j in {1,...,4} \draw(0.5,\j-0.5)node{$\ast$};
\foreach \j in {2,...,4} \draw(1.5,\j-0.5)node{$\ast$};
\foreach \j in {4,...,8} \draw (1-0.5,\j+0.5)node{$\j$}; 
\foreach \j in {3,...,5} \draw (2-0.5,\j+0.5+1)node{$\j$};
\foreach \j in {0,...,3} \draw (3-0.5,\j+0.5+2)node{$\j$};
\foreach \j in {0,...,1} \draw (4-0.5,\j+0.5+3)node{$\j$};
\end{tikzpicture}
}}
&
\vcenter{\hbox{
\begin{tikzpicture}[x={(0in,-0.2in)},y={(0.2in,0in)}] 
\draw[thick] (1,9) rectangle +(-1,-1);\draw(1-0.5,9-0.5)node{$\brown8$};
\draw[thick] (1,8) rectangle +(-1,-1);\draw(1-0.5,8-0.5)node{$\brown7$};
\draw[thick] (1,7) rectangle +(-1,-1);\draw(1-0.5,7-0.5)node{$\brown6$};
\draw[thick] (1,6) rectangle +(-1,-1);\draw(1-0.5,6-0.5)node{$\brown5$};
\draw[thick] (1,5) rectangle +(-1,-1);\draw(1-0.5,5-0.5)node{$\green4$};
\draw[thick] (2,5) rectangle +(-1,-1);\draw(2-0.5,5-0.5)node{$\green3$};
\draw[thick] (3,5) rectangle +(-1,-1);\draw(3-0.5,5-0.5)node{$\brown2$};
\draw[thick] (3,4) rectangle +(-1,-1);\draw(3-0.5,4-0.5)node{$\brown1$};
\draw[thick] (3,3) rectangle +(-1,-1);\draw(3-0.5,3-0.5)node{$\brown0$};
\end{tikzpicture}
}}
\cr\cr
\mbox{Superposition of cutting strips}&\mbox{Outside decomposition}
\cr\cr
\vcenter{\hbox{
\begin{tikzpicture}[x={(0in,-0.2in)},y={(0.2in,0in)}] 
\foreach \j in {1,...,9} \draw[thick] (1,\j) rectangle +(-1,-1);
\foreach \j in {2,...,7} \draw[thick] (2,\j) rectangle +(-1,-1);
\foreach \j in {3,...,6} \draw[thick] (3,\j) rectangle +(-1,-1);
\foreach \j in {4,...,5} \draw[thick] (4,\j) rectangle +(-1,-1);
\fill[blue!60!white](2.3,1.7)--(2.3,4.3)--(0.3,4.3)--(0.3,9.3)--(0.7,9.3)--(0.7,4.7)--(2.7,4.7)--(2.7,1.7)--cycle;
\fill[red!60!white](3.3,2.7)--(3.3,5.3)--(1.3,5.3)--(1.3,9.3)--(1.7,9.3)--(1.7,5.7)--(3.7,5.7)--(3.7,2.7)--cycle;
\foreach \j in {1,...,4} \draw[thick](0.5,\j-0.5)node{$\ast$};
\foreach \j in {2,...,4} \draw[thick](1.5,\j-0.5)node{$\ast$};
\end{tikzpicture}
}}
&
\begin{array}{c}
\vcenter{\hbox{
\begin{tikzpicture}[x={(0in,-0.2in)},y={(0.2in,0in)}] 
\foreach \j in {1,...,9} \draw[thick] (1,\j) rectangle +(-1,-1);
\foreach \j in {2,...,7} \draw[thick] (2,\j) rectangle +(-1,-1);
\foreach \j in {3,...,6} \draw[thick] (3,\j) rectangle +(-1,-1);
\foreach \j in {4,...,5} \draw[thick] (4,\j) rectangle +(-1,-1);
\foreach \j in {1,...,4} \draw(0.5,\j-0.5)node{$\ast$};
\foreach \j in {2,...,4} \draw(1.5,\j-0.5)node{$\ast$};
\draw(2.5,1.3)node{$\theta_1$};
\draw[draw=blue,ultra thick] (2.5,1.8)--(2.5,4.5)--(0.5,4.5)--(0.5,9.2);
\draw(3.5,2.3)node{$\theta_2$};
\draw[draw=red, ultra thick] (3.5,2.8)--(3.5,5.2);
\draw(1.5,7.7)node{$\theta_3$};
\draw[draw=cyan,ultra thick] (3.2,5.5)--(1.5,5.5)--(1.5,7.2);
\end{tikzpicture}
}}\cr\cr
\end{array}
\end{array}
\end{equation}
It can be seen that the individual strips $\theta_p$ are given, along with their content 
and identification with subdiagrams of the cutting strip $\phi$, by
\begin{equation}\label{eqn-stripsQ}
\theta_1=\vcenter{\hbox{
\begin{tikzpicture}[x={(0in,-0.2in)},y={(0.2in,0in)}] 
\foreach \j in{6,7,8,9} \draw[thick] (1,\j) rectangle +(-1,-1);
\foreach \i in{1,2,3} \draw[thick] (\i,5) rectangle +(-1,-1);
\foreach \j in{3,4} \draw[thick] (3,\j) rectangle +(-1,-1);
\draw(1-0.5,9-0.5)node{$\blue8$};
\draw(1-0.5,8-0.5)node{$\blue7$};
\draw(1-0.5,7-0.5)node{$\blue6$};
\draw(1-0.5,6-0.5)node{$\blue5$};
\draw(1-0.5,5-0.5)node{$\blue4$};
\draw(2-0.5,5-0.5)node{$\blue3$};
\draw(3-0.5,5-0.5)node{$\blue2$};
\draw(3-0.5,4-0.5)node{$\blue1$};
\draw(3-0.5,3-0.5)node{$\blue0$};
\end{tikzpicture}
}}\simeq\phi_{08}
\quad
\theta_2=\vcenter{\hbox{
\begin{tikzpicture}[x={(0in,-0.2in)},y={(0.2in,0in)}] 
\draw[thick] (4,2) rectangle +(-1,-1);\draw(4-0.5,2-0.5)node{$\red1$};
\draw[thick] (4,1) rectangle +(-1,-1);\draw(4-0.5,1-0.5)node{$\red0$};
\end{tikzpicture}
}}\simeq\phi_{01}
\quad
\theta_3=\vcenter{\hbox{
\begin{tikzpicture}[x={(0in,-0.2in)},y={(0.2in,0in)}] 
\draw[thick] (2,4) rectangle +(-1,-1);\draw(2-0.5,4-0.5)node{$\cyan5$};
\draw[thick] (2,3) rectangle +(-1,-1);\draw(2-0.5,3-0.5)node{$\cyan4$};
\draw[thick] (3,3) rectangle +(-1,-1);\draw(3-0.5,3-0.5)node{$\cyan3$};
\end{tikzpicture}
}}\simeq\phi_{35}
\end{equation}

Thus a given cutting strip $\phi$ may be used to specify a corresponding outside decomposition
$\theta=(\theta_1,\theta_2,\ldots,\theta_s)$, with each constituent strip 
identifiable in a content preserving manner with a substrip of $\phi$. As before,
if $\theta_p\simeq\phi_{ab}$ and $\theta_q\simeq\phi_{cd}$ then $(\theta_p\#\theta_q)=\phi_{ad}$. 
In order to extend the strips of our outside decomposition to the main diagonal, let $\rho$ be a strip consisting of a 
single box of content $0$ so that $\rho\simeq\phi_{00}$. This allows us
to define $\ov{\theta_p}\simeq(\rho\#\theta_p)=\phi_{0b}$ and $\ov{\theta_q}\simeq(\rho\#\theta_q)=\phi_{0d}$.
Combining these, let $(\ov{\theta}_p,\ov{\theta}_q)$ be formed by juxtaposing $\ov{\theta}_p$ and $\ov{\theta}_q$ 
with their boxes of content $0$ lying on the main diagonal with that of $\ov{\theta}_p$ immediately above and to the left of $\ov{\theta}_q$.
It will be noted that $(\ov{\theta}_p,\ov{\theta}_q)$ will be a standard skew shifted Young diagram if $b>d$, but not if $b\leq d$. In particular 
$(\ov{\theta}_p,\ov{\theta}_p)$ is a non-standard skew shifted Young diagram.

This is illustrated in our running example based on (\ref{eqn-csQ}) and (\ref{eqn-stripsQ}) by
\begin{equation}\label{eqn-ovp}
\begin{array}{|c|c|}
\hline
&\cr
\theta_p~~\mbox{and}~~\ov{\theta}_p=(\rho\#\theta_p)&(\theta_3\#\theta_p)\cr
\hline
&\cr
\theta_1=
\vcenter{\hbox{
\begin{tikzpicture}[x={(0in,-0.2in)},y={(0.2in,0in)}] 
\draw[thick] (1,9) rectangle +(-1,-1);\draw(1-0.5,9-0.5)node{$\blue8$};
\draw[thick] (1,8) rectangle +(-1,-1);\draw(1-0.5,8-0.5)node{$\blue7$};
\draw[thick] (1,7) rectangle +(-1,-1);\draw(1-0.5,7-0.5)node{$\blue6$};
\draw[thick] (1,6) rectangle +(-1,-1);\draw(1-0.5,6-0.5)node{$\blue5$};
\draw[thick] (1,5) rectangle +(-1,-1);\draw(1-0.5,5-0.5)node{$\blue4$};
\draw[thick] (2,5) rectangle +(-1,-1);\draw(2-0.5,5-0.5)node{$\blue3$};
\draw[thick] (3,5) rectangle +(-1,-1);\draw(3-0.5,5-0.5)node{$\blue2$};
\draw[thick] (3,4) rectangle +(-1,-1);\draw(3-0.5,4-0.5)node{$\blue1$};
\draw[thick] (3,3) rectangle +(-1,-1);\draw(3-0.5,3-0.5)node{$\blue0$};
\end{tikzpicture}
}}=\ov{\theta}_1
&
\vcenter{\hbox{
\begin{tikzpicture}[x={(0in,-0.2in)},y={(0.2in,0in)}] 
\draw[thick] (1,9) rectangle +(-1,-1);\draw(1-0.5,9-0.5)node{$\blue8$};
\draw[thick] (1,8) rectangle +(-1,-1);\draw(1-0.5,8-0.5)node{$\blue7$};
\draw[thick] (1,7) rectangle +(-1,-1);\draw(1-0.5,7-0.5)node{$\blue6$};
\draw[thick] (1,6) rectangle +(-1,-1);\draw(1-0.5,6-0.5)node{$\blue5$};
\draw[thick] (1,5) rectangle +(-1,-1);\draw(1-0.5,5-0.5)node{$\blue4$};
\draw[thick] (2,5) rectangle +(-1,-1);\draw(2-0.5,5-0.5)node{$\blue3$};
\end{tikzpicture}
}}
\cr&\cr
\theta_2=
\vcenter{\hbox{
\begin{tikzpicture}[x={(0in,-0.2in)},y={(0.2in,0in)}] 
\draw[thick] (3,4) rectangle +(-1,-1);\draw(3-0.5,4-0.5)node{$\red1$};
\draw[thick] (3,3) rectangle +(-1,-1);\draw(3-0.5,3-0.5)node{$\red0$};
\end{tikzpicture}
}}=\ov{\theta}_2
&
\emptyset
\cr&\cr
\theta_3=
\vcenter{\hbox{
\begin{tikzpicture}[x={(0in,-0.2in)},y={(0.2in,0in)}] 
\draw[thick] (1,6) rectangle +(-1,-1);\draw(1-0.5,6-0.5)node{$\cyan5$};
\draw[thick] (1,5) rectangle +(-1,-1);\draw(1-0.5,5-0.5)node{$\cyan4$};
\draw[thick] (2,5) rectangle +(-1,-1);\draw(2-0.5,5-0.5)node{$\cyan3$};
\end{tikzpicture}
}}
\quad
\ov{\theta}_3=
\vcenter{\hbox{
\begin{tikzpicture}[x={(0in,-0.2in)},y={(0.2in,0in)}] 
\draw[thick] (1,6) rectangle +(-1,-1);\draw(1-0.5,6-0.5)node{$\cyan5$};
\draw[thick] (1,5) rectangle +(-1,-1);\draw(1-0.5,5-0.5)node{$\cyan4$};
\draw[thick] (2,5) rectangle +(-1,-1);\draw(2-0.5,5-0.5)node{$\cyan3$};
\draw[thick] (3,5) rectangle +(-1,-1);\draw(3-0.5,5-0.5)node{$\cyan2$};
\draw[thick] (3,4) rectangle +(-1,-1);\draw(3-0.5,4-0.5)node{$\cyan1$};
\draw[thick] (3,3) rectangle +(-1,-1);\draw(3-0.5,3-0.5)node{$\cyan0$};
\end{tikzpicture}
}}
&
\vcenter{\hbox{
\begin{tikzpicture}[x={(0in,-0.2in)},y={(0.2in,0in)}] 
\draw[thick] (1,6) rectangle +(-1,-1);\draw(1-0.5,6-0.5)node{$\cyan5$};
\draw[thick] (1,5) rectangle +(-1,-1);\draw(1-0.5,5-0.5)node{$\cyan4$};
\draw[thick] (2,5) rectangle +(-1,-1);\draw(2-0.5,5-0.5)node{$\cyan3$};
\end{tikzpicture}
}}\cr
&\cr
\hline
\end{array}
\end{equation}
and
\begin{equation}\label{eqn-ovpq}
\begin{array}{|c|ccc|}
\hline
&&&\cr
(\ov{\theta}_p,\ov{\theta}_q)&\ov{\theta}_1&\ov{\theta}_2&\ov{\theta}_3\cr
\hline
&&&\cr
\ov{\theta}_1&
\vcenter{\hbox{
\begin{tikzpicture}[x={(0in,-0.2in)},y={(0.2in,0in)}] 
\draw[thick] (1,9) rectangle +(-1,-1);\draw(1-0.5,9-0.5)node{$\blue8$};
\draw[thick] (1,8) rectangle +(-1,-1);\draw(1-0.5,8-0.5)node{$\blue7$};
\draw[thick] (1,7) rectangle +(-1,-1);\draw(1-0.5,7-0.5)node{$\blue6$};
\draw[thick] (1,6) rectangle +(-1,-1);\draw(1-0.5,6-0.5)node{$\blue5$};
\draw[thick] (1,5) rectangle +(-1,-1);\draw(1-0.5,5-0.5)node{$\blue4$};
\draw[thick] (2,5) rectangle +(-1,-1);\draw(2-0.5,5-0.5)node{$\blue3$};
\draw[thick] (3,5) rectangle +(-1,-1);\draw(3-0.5,5-0.5)node{$\blue2$};
\draw[thick] (3,4) rectangle +(-1,-1);\draw(3-0.5,4-0.5)node{$\blue1$};
\draw[thick] (3,3) rectangle +(-1,-1);\draw(3-0.5,3-0.5)node{$\blue0$};
\draw[thick] (1+1,9+1) rectangle +(-1,-1);\draw(1+0.5,9+0.5)node{$\blue8$};
\draw[thick] (1+1,8+1) rectangle +(-1,-1);\draw(1+0.5,8+0.5)node{$\blue7$};
\draw[thick] (1+1,7+1) rectangle +(-1,-1);\draw(1+0.5,7+0.5)node{$\blue6$};
\draw[thick] (1+1,6+1) rectangle +(-1,-1);\draw(1+0.5,6+0.5)node{$\blue5$};
\draw[thick] (1+1,5+1) rectangle +(-1,-1);\draw(1+0.5,5+0.5)node{$\blue4$};
\draw[thick] (2+1,5+1) rectangle +(-1,-1);\draw(2+0.5,5+0.5)node{$\blue3$};
\draw[thick] (3+1,5+1) rectangle +(-1,-1);\draw(3+0.5,5+0.5)node{$\blue2$};
\draw[thick] (3+1,4+1) rectangle +(-1,-1);\draw(3+0.5,4+0.5)node{$\blue1$};
\draw[thick] (3+1,3+1) rectangle +(-1,-1);\draw(3+0.5,3+0.5)node{$\blue0$};
\end{tikzpicture}
}}
&
\vcenter{\hbox{
\begin{tikzpicture}[x={(0in,-0.2in)},y={(0.2in,0in)}] 
\draw[thick] (1,9) rectangle +(-1,-1);\draw(1-0.5,9-0.5)node{$\blue8$};
\draw[thick] (1,8) rectangle +(-1,-1);\draw(1-0.5,8-0.5)node{$\blue7$};
\draw[thick] (1,7) rectangle +(-1,-1);\draw(1-0.5,7-0.5)node{$\blue6$};
\draw[thick] (1,6) rectangle +(-1,-1);\draw(1-0.5,6-0.5)node{$\blue5$};
\draw[thick] (1,5) rectangle +(-1,-1);\draw(1-0.5,5-0.5)node{$\blue4$};
\draw[thick] (2,5) rectangle +(-1,-1);\draw(2-0.5,5-0.5)node{$\blue3$};
\draw[thick] (3,5) rectangle +(-1,-1);\draw(3-0.5,5-0.5)node{$\blue2$};
\draw[thick] (3,4) rectangle +(-1,-1);\draw(3-0.5,4-0.5)node{$\blue1$};
\draw[thick] (3,3) rectangle +(-1,-1);\draw(3-0.5,3-0.5)node{$\blue0$};
\draw[thick] (3+1,4+1) rectangle +(-1,-1);\draw(3+1-0.5,4+1-0.5)node{$\red1$};
\draw[thick] (3+1,3+1) rectangle +(-1,-1);\draw(3+1-0.5,3+1-0.5)node{$\red0$};
\end{tikzpicture}
}}
&
\vcenter{\hbox{
\begin{tikzpicture}[x={(0in,-0.2in)},y={(0.2in,0in)}] 
\draw[thick] (1,9) rectangle +(-1,-1);\draw(1-0.5,9-0.5)node{$\blue8$};
\draw[thick] (1,8) rectangle +(-1,-1);\draw(1-0.5,8-0.5)node{$\blue7$};
\draw[thick] (1,7) rectangle +(-1,-1);\draw(1-0.5,7-0.5)node{$\blue6$};
\draw[thick] (1,6) rectangle +(-1,-1);\draw(1-0.5,6-0.5)node{$\blue5$};
\draw[thick] (1,5) rectangle +(-1,-1);\draw(1-0.5,5-0.5)node{$\blue4$};
\draw[thick] (2,5) rectangle +(-1,-1);\draw(2-0.5,5-0.5)node{$\blue3$};
\draw[thick] (3,5) rectangle +(-1,-1);\draw(3-0.5,5-0.5)node{$\blue2$};
\draw[thick] (3,4) rectangle +(-1,-1);\draw(3-0.5,4-0.5)node{$\blue1$};
\draw[thick] (3,3) rectangle +(-1,-1);\draw(3-0.5,3-0.5)node{$\blue0$};
\draw[thick] (1+1,6+1) rectangle +(-1,-1);\draw(1+1-0.5,6+1-0.5)node{$\cyan5$};
\draw[thick] (1+1,5+1) rectangle +(-1,-1);\draw(1+1-0.5,5+1-0.5)node{$\cyan4$};
\draw[thick] (2+1,5+1) rectangle +(-1,-1);\draw(2+1-0.5,5+1-0.5)node{$\cyan3$};
\draw[thick] (3+1,5+1) rectangle +(-1,-1);\draw(3+1-0.5,5+1-0.5)node{$\cyan2$};
\draw[thick] (3+1,4+1) rectangle +(-1,-1);\draw(3+1-0.5,4+1-0.5)node{$\cyan1$};
\draw[thick] (3+1,3+1) rectangle +(-1,-1);\draw(3+1-0.5,3+1-0.5)node{$\cyan0$};
\end{tikzpicture}
}}
\cr
&&&\cr
\ov{\theta}_2&
\vcenter{\hbox{
\begin{tikzpicture}[x={(0in,-0.2in)},y={(0.2in,0in)}] 
\draw[thick] (3-1,4-1) rectangle +(-1,-1);\draw(3-1-0.5,4-1-0.5)node{$\red1$};
\draw[thick] (3-1,3-1) rectangle +(-1,-1);\draw(3-1-0.5,3-1-0.5)node{$\red0$};
\draw[thick] (1,9) rectangle +(-1,-1);\draw(1-0.5,9-0.5)node{$\blue8$};
\draw[thick] (1,8) rectangle +(-1,-1);\draw(1-0.5,8-0.5)node{$\blue7$};
\draw[thick] (1,7) rectangle +(-1,-1);\draw(1-0.5,7-0.5)node{$\blue6$};
\draw[thick] (1,6) rectangle +(-1,-1);\draw(1-0.5,6-0.5)node{$\blue5$};
\draw[thick] (1,5) rectangle +(-1,-1);\draw(1-0.5,5-0.5)node{$\blue4$};
\draw[thick] (2,5) rectangle +(-1,-1);\draw(2-0.5,5-0.5)node{$\blue3$};
\draw[thick] (3,5) rectangle +(-1,-1);\draw(3-0.5,5-0.5)node{$\blue2$};
\draw[thick] (3,4) rectangle +(-1,-1);\draw(3-0.5,4-0.5)node{$\blue1$};
\draw[thick] (3,3) rectangle +(-1,-1);\draw(3-0.5,3-0.5)node{$\blue0$};
\end{tikzpicture}
}}
&
\vcenter{\hbox{
\begin{tikzpicture}[x={(0in,-0.2in)},y={(0.2in,0in)}] 
\draw[thick] (3-1,4-1) rectangle +(-1,-1);\draw(3-1-0.5,4-1-0.5)node{$\red1$};
\draw[thick] (3-1,3-1) rectangle +(-1,-1);\draw(3-1-0.5,3-1-0.5)node{$\red0$};
\draw[thick] (3,4) rectangle +(-1,-1);\draw(3-0.5,4-0.5)node{$\red1$};
\draw[thick] (3,3) rectangle +(-1,-1);\draw(3-0.5,3-0.5)node{$\red0$};
\end{tikzpicture}
}}
&
\vcenter{\hbox{
\begin{tikzpicture}[x={(0in,-0.2in)},y={(0.2in,0in)}] 
\draw[thick] (3,4) rectangle +(-1,-1);\draw(3-0.5,4-0.5)node{$\red1$};
\draw[thick] (3,3) rectangle +(-1,-1);\draw(3-0.5,3-0.5)node{$\red0$};
\draw[thick] (1+1,6+1) rectangle +(-1,-1);\draw(1+1-0.5,6+1-0.5)node{$\cyan5$};
\draw[thick] (1+1,5+1) rectangle +(-1,-1);\draw(1+1-0.5,5+1-0.5)node{$\cyan4$};
\draw[thick] (2+1,5+1) rectangle +(-1,-1);\draw(2+1-0.5,5+1-0.5)node{$\cyan3$};
\draw[thick] (3+1,5+1) rectangle +(-1,-1);\draw(3+1-0.5,5+1-0.5)node{$\cyan2$};
\draw[thick] (3+1,4+1) rectangle +(-1,-1);\draw(3+1-0.5,4+1-0.5)node{$\cyan1$};
\draw[thick] (3+1,3+1) rectangle +(-1,-1);\draw(3+1-0.5,3+1-0.5)node{$\cyan0$};
\end{tikzpicture}
}}
\cr
&&&\cr
\ov{\theta}_3&
\vcenter{\hbox{
\begin{tikzpicture}[x={(0in,-0.2in)},y={(0.2in,0in)}] 
\draw[thick] (1-1,6-1) rectangle +(-1,-1);\draw(1-1-0.5,6-1-0.5)node{$\cyan5$};
\draw[thick] (1-1,5-1) rectangle +(-1,-1);\draw(1-1-0.5,5-1-0.5)node{$\cyan4$};
\draw[thick] (2-1,5-1) rectangle +(-1,-1);\draw(2-1-0.5,5-1-0.5)node{$\cyan3$};
\draw[thick] (3-1,5-1) rectangle +(-1,-1);\draw(3-1-0.5,5-1-0.5)node{$\cyan2$};
\draw[thick] (3-1,4-1) rectangle +(-1,-1);\draw(3-1-0.5,4-1-0.5)node{$\cyan1$};
\draw[thick] (3-1,3-1) rectangle +(-1,-1);\draw(3-1-0.5,3-1-0.5)node{$\cyan0$};
\draw[thick] (1,9) rectangle +(-1,-1);\draw(1-0.5,9-0.5)node{$\blue8$};
\draw[thick] (1,8) rectangle +(-1,-1);\draw(1-0.5,8-0.5)node{$\blue7$};
\draw[thick] (1,7) rectangle +(-1,-1);\draw(1-0.5,7-0.5)node{$\blue6$};
\draw[thick] (1,6) rectangle +(-1,-1);\draw(1-0.5,6-0.5)node{$\blue5$};
\draw[thick] (1,5) rectangle +(-1,-1);\draw(1-0.5,5-0.5)node{$\blue4$};
\draw[thick] (2,5) rectangle +(-1,-1);\draw(2-0.5,5-0.5)node{$\blue3$};
\draw[thick] (3,5) rectangle +(-1,-1);\draw(3-0.5,5-0.5)node{$\blue2$};
\draw[thick] (3,4) rectangle +(-1,-1);\draw(3-0.5,4-0.5)node{$\blue1$};
\draw[thick] (3,3) rectangle +(-1,-1);\draw(3-0.5,3-0.5)node{$\blue0$};
\end{tikzpicture}
}}
&
\vcenter{\hbox{
\begin{tikzpicture}[x={(0in,-0.2in)},y={(0.2in,0in)}] 
\draw[thick] (1-1,6-1) rectangle +(-1,-1);\draw(1-1-0.5,6-1-0.5)node{$\cyan5$};
\draw[thick] (1-1,5-1) rectangle +(-1,-1);\draw(1-1-0.5,5-1-0.5)node{$\cyan4$};
\draw[thick] (2-1,5-1) rectangle +(-1,-1);\draw(2-1-0.5,5-1-0.5)node{$\cyan3$};
\draw[thick] (3-1,5-1) rectangle +(-1,-1);\draw(3-1-0.5,5-1-0.5)node{$\cyan2$};
\draw[thick] (3-1,4-1) rectangle +(-1,-1);\draw(3-1-0.5,4-1-0.5)node{$\cyan1$};
\draw[thick] (3-1,3-1) rectangle +(-1,-1);\draw(3-1-0.5,3-1-0.5)node{$\cyan0$};
\draw[thick] (3,4) rectangle +(-1,-1);\draw(3-0.5,4-0.5)node{$\red1$};
\draw[thick] (3,3) rectangle +(-1,-1);\draw(3-0.5,3-0.5)node{$\red0$};
\end{tikzpicture}
}}
&
\vcenter{\hbox{
\begin{tikzpicture}[x={(0in,-0.2in)},y={(0.2in,0in)}] 
\draw[thick] (1-1,6-1) rectangle +(-1,-1);\draw(1-1-0.5,6-1-0.5)node{$\cyan5$};
\draw[thick] (1-1,5-1) rectangle +(-1,-1);\draw(1-1-0.5,5-1-0.5)node{$\cyan4$};
\draw[thick] (2-1,5-1) rectangle +(-1,-1);\draw(2-1-0.5,5-1-0.5)node{$\cyan3$};
\draw[thick] (3-1,5-1) rectangle +(-1,-1);\draw(3-1-0.5,5-1-0.5)node{$\cyan2$};
\draw[thick] (3-1,4-1) rectangle +(-1,-1);\draw(3-1-0.5,4-1-0.5)node{$\cyan1$};
\draw[thick] (3-1,3-1) rectangle +(-1,-1);\draw(3-1-0.5,3-1-0.5)node{$\cyan0$};
\draw[thick] (1,6) rectangle +(-1,-1);\draw(1-0.5,6-0.5)node{$\cyan5$};
\draw[thick] (1,5) rectangle +(-1,-1);\draw(1-0.5,5-0.5)node{$\cyan4$};
\draw[thick] (2,5) rectangle +(-1,-1);\draw(2-0.5,5-0.5)node{$\cyan3$};
\draw[thick] (3,5) rectangle +(-1,-1);\draw(3-0.5,5-0.5)node{$\cyan2$};
\draw[thick] (3,4) rectangle +(-1,-1);\draw(3-0.5,4-0.5)node{$\cyan1$};
\draw[thick] (3,3) rectangle +(-1,-1);\draw(3-0.5,3-0.5)node{$\cyan0$};
\end{tikzpicture}
}}
\cr
&&&\cr
\hline
\end{array}
\end{equation}

It should be noted here that the double strip diagrams $(\ov{\theta}_p,\ov{\theta}_q)$ coincide with 
standard skew shifted diagrams of shape $SF^{\kappa/\nu}$ with $\kappa$ and $\nu$ both strict partitions
provided that $\ov{\theta}_p\simeq\phi_{0b}$ and $\ov{\theta}_q\simeq\phi_{0d}$ with $b>d$. 
In such a case $\kappa$ may be fixed by ensuring that the main diagonal of $SF^{\kappa}$ passes through
the boxes of content $0$ in both $\phi_{0b}$ and $\phi_{0d}$. The fact that these two strips are anchored
to the main diagonal ensures that the content of each box of $(\ov{\theta}_p,\ov{\theta}_q)$
coincides with the content of the corresponding box of $SF^{\kappa/\nu}$, without the necessity of invoking any 
shift parameter. The same is not true of the single strip diagrams 
$(\theta_k\#\theta_p)\simeq\phi_{ab}$ with $b\geq a>0$
for which it is necessary to define a shift parameter $m(k,p)=a$. This is in the spirit of
our previous definition of a shift parameter, in the sense that $\phi_{ab}$ is to
be interpreted as a skew shifted diagram of shape $SF^{\kappa/\nu}$ with $\kappa$ chosen so that
the main diagonal of $SF^\kappa$ intersects the lower left hand end of $\phi_{ab}$ in a box of content $a$.

This completes the combinatorial machinery relevant to the statement of Theorems~\ref{The-HG}
and~\ref{The-Ham} and their generalisations which are the subjects of the next section.

\section{Main Theorems}\label{Sec-main-theorems}

\subsection{Generalisation of Hamel and Goulden's Theorem~\ref{The-HG}}

\begin{Theorem}\label{The-HGX} 
For partitions $\lambda$ and $\mu$ such that $\mu\subseteq\lambda$, let 
$\Theta=(\theta_1,\theta_2,\ldots,\theta_s)$ be the outside decomposition of $F^{\lambda/\mu}$
corresponding to a cutting strip $\phi$. Let $\tau$ be a shift operator whose action on $\X=(x_{kc})$
with $k\in\{1,2,\ldots,n\}$ and $c\in\Z$ is such that $\tau^m\X=(x_{k,m+c})$ for any $m\in\Z$.

Then the generalised Schur functions of (\ref{eqn-sfnX}) satisfy the identity
\begin{equation}\label{eqn-HGX}
    s_{\lambda/\mu}(\X) = \left|\,  s_{(\theta_p\#\theta_q)}(\tau^{m(p,q)}\X) \,\right|_{1\leq p\leq q\leq s}  \,,
\end{equation}
where $m(p,q)$ is the value of the content of the box of $\phi_{ad}$ that lies at the intersection of the main diagonal of the diagram having $\phi_{ad}$ as its outer rim, and
where if $(\theta_p\#\theta_q)=\phi_{ad}$ then
\begin{equation}\label{eqn-phashq}
   s_{(\theta_p\#\theta_q)}(\tau^{m(p,q)}\X)= \begin{cases} 
	                                               s_{\phi_{ad}}(\tau^{m(p,q)}\X)&\mbox{if $a\leq d$};\cr
		                                             1&\mbox{if $a=d+1$};\cr
																							   0&\mbox{if $a>d+1$}.\cr
															                \end{cases}
\end{equation}																	
\end{Theorem}

\noindent{\bf Proof}: In order to prove Theorem~\ref{The-HG} Hamel and Goulden~\cite{HG} established a bijection between
semistandard skew tableaux $T\in{\cal T}^{\lambda/\mu}$ and $L\in{\cal L}^\phi$, the set of $s$-tuples of non-intersecting 
paths on a rectangular lattice with directed edges determined by the cutting strip $\phi$ and paths determined by the entries 
of $T$ in the constituent strips $\theta_p$ of the corresponding outside decomposition $\Theta=(\theta_1,\theta_2,\ldots,\theta_s)$ 
of $F^{\lambda/\mu}$. They showed furthermore that this same map takes non-semistandard tableaux to $s$-tuples of paths involving 
at least one intersection. The same bijection can be used here. 

It can be exemplified as follows. 
In the case $\lambda=(5,4,4,2)$ and $\mu=(3,2)$ a typical tableau $T\in{\cal T}^{\lambda/\mu}$ with entries $k$ from $\{1<2<3<4\}$
is as shown below on the right alongside a reminder of the content $c$ of the boxes of $F^{\lambda/\mu}$ on the left.
\begin{equation}\label{eqn-FT}
\begin{array}{cc}
\vcenter{\hbox{
\begin{tikzpicture}[x={(0in,-0.2in)},y={(0.2in,0in)}] 
\foreach \j in {1,...,5} \draw[thick] (1,\j) rectangle +(-1,-1);
\foreach \j in {1,...,4} \draw[thick] (2,\j) rectangle +(-1,-1);
\foreach \j in {1,...,4} \draw[thick] (3,\j) rectangle +(-1,-1);
\foreach \j in {1,...,2} \draw[thick] (4,\j) rectangle +(-1,-1);
\foreach \j in {1,...,3} \draw(1-0.5,\j-0.5)node{$\ast$};
\foreach \j in {1,...,2} \draw(2-0.5,\j-0.5)node{$\ast$};
\draw(1-0.5,4-0.5) node {$3$};\draw(1-0.5,5-0.5) node {$4$};
\draw(2-0.5,3-0.5) node {$1$};\draw(2-0.5,4-0.5) node {$2$};
\draw(3-0.5,1-0.5) node {$\ov2$};\draw(3-0.5,2-0.5) node {$\ov1$};\draw(3-0.5,3-0.5) node {$0$};\draw(3-0.5,4-0.5) node {$1$};
\draw(4-0.5,1-0.5) node {$\ov3$};\draw(4-0.5,2-0.5) node {$\ov2$};
\end{tikzpicture}
}}
&
\vcenter{\hbox{
\begin{tikzpicture}[x={(0in,-0.2in)},y={(0.2in,0in)}] 
\foreach \j in {1,...,5} \draw[thick] (1,\j) rectangle +(-1,-1);
\foreach \j in {1,...,4} \draw[thick] (2,\j) rectangle +(-1,-1);
\foreach \j in {1,...,4} \draw[thick] (3,\j) rectangle +(-1,-1);
\foreach \j in {1,...,2} \draw[thick] (4,\j) rectangle +(-1,-1);
\foreach \j in {1,...,3} \draw(1-0.5,\j-0.5)node{$\ast$};
\foreach \j in {1,...,2} \draw(2-0.5,\j-0.5)node{$\ast$};
\draw(1-0.5,4-0.5)node{$\cyan1$};   \draw(1-0.5,5-0.5) node {$\cyan3$};
\draw(2-0.5,3-0.5)node{$\cyan1$};    \draw(2-0.5,4-0.5) node {$\cyan2$};
\draw(3-0.5,1-0.5)node{$\red1$};\draw(3-0.5,2-0.5) node {$\blue2$};\draw(3-0.5,3-0.5) node {$\blue2$};\draw(3-0.5,4-0.5) node {$\blue4$};
\draw(4-0.5,1-0.5)node{$\blue3$};\draw(4-0.5,2-0.5) node {$\blue3$};
\end{tikzpicture}
}}\cr
\end{array}
\end{equation}

The required lattice is set up on a rectangular grid with coordinates $(k,c)$ and directed edges determined by the chosen cutting strip.
In our example with  
\begin{equation}
  \phi = \vcenter{\hbox{
\begin{tikzpicture}[x={(0in,-0.2in)},y={(0.2in,0in)}] 
\foreach \j in {5,6} \draw[thick] (1,\j) rectangle +(-1,-1);
\foreach \j in {2,...,5} \draw[thick] (2,\j) rectangle +(-1,-1);
\foreach \j in {1,...,2} \draw[thick] (3,\j) rectangle +(-1,-1);
\draw(1-0.5,5-0.5) node {$\green{3}$}; \draw(1-0.5,6-0.5) node {$\brown{4}$};
\draw(2-0.5,2-0.5) node {$\green{\ov1}$}; \draw(2-0.5,3-0.5) node {$\brown{0}$}; \draw(2-0.5,4-0.5) node {$\brown{1}$}; \draw(2-0.5,5-0.5) node {$\brown{2}$};
\draw(3-0.5,1-0.5) node {$\green{\ov3}$}; \draw(3-0.5,2-0.5) node {$\brown{\ov2}$};
\end{tikzpicture}
}}
\end{equation}
we have
\begin{equation}
\vcenter{\hbox{
\begin{tikzpicture}[x={(0in,-0.4in)},y={(0.4in,0in)}] 
\foreach \i in {1,...,4} \draw(\i,-5)node{$k=\i$};
\draw(-0.5,-5)node{$c=$};
\foreach \j in {1,3} \draw(-0.6,-\j)node{$\green{\ov\j}$};
\foreach \j in {2,4} \draw(-0.6,-\j)node{$\brown{\ov\j}$};
\foreach \j in {3} \draw(-0.5,\j)node{$\green{\j}$};
\foreach \j in {0,1,2,4} \draw(-0.5,\j)node{$\brown{\j}$};
\foreach \i in {1,...,4} \foreach \j in {-2,0,1,2,4} \draw[-latex](\i,\j-1)--(\i,\j-0.4);
\foreach \i in {1,...,4} \foreach \j in {-2,0,1,2,4} \draw(\i,\j-0.4)--(\i,\j);
\foreach \i in {1,...,4} \foreach \j in {-3,-1,3} \draw[-latex](\i+1,\j-1)--(\i+0.4,\j-0.4); 
\foreach \i in {1,...,4} \foreach \j in {-3,-1,3} \draw(\i+0.4,\j-0.4)--(\i,\j);
\foreach \j in {-3,-1,0,1,3,4} \foreach \i in {1,...,5} \draw[-latex](\i-1,\j)--(\i-0.4,\j);
\foreach \j in {-3,-1,0,1,3,4} \foreach \i in {1,...,5}\draw(\i-0.4,\j)--(\i,\j);
\foreach \j in {-4,-2,2} \foreach \i in {0,...,4} \draw[-latex](\i+1,\j)--(\i+0.4,\j);
\foreach \j in {-4,-2,2} \foreach \i in {0,...,4} \draw(\i+0.4,\j)--(\i,\j);
\foreach \i in {0,...,5} \foreach \j in {-4,...,4} \draw(\i,\j)node{$\ssc\bullet$};
\end{tikzpicture}
}}
\end{equation}
This is constructed by reading $\phi$ from bottom left to top right, then if a box of content $c$ in $\phi$ is approached either from the left or from below 
then the edges from column $c-1$ to $c$ are either horizontal from left to right or diagonal from lower left to upper right, and the vertical edges
in column $c$ are up or down, respectively. In the case of the column corresponding to the lowest value of the content in $\phi$ the approach could be horizontal or diagonal, and we have chosen the latter.

The start and end points of the lattice paths $P_p$ and $Q_p$ corresponding to the entries of $\theta_p$ in $T$ are determined~\cite{HG}
by the manner in which $\theta_p\simeq\phi_{ab}$, when followed from bottom left to top right, enters and leaves $F^{\lambda/\mu}$:
\begin{itemize}
\item if $\theta_p\simeq\phi_{ab}$ enters from the left then $P_p=(0,a-1)$;
\item if $\theta_p\simeq\phi_{ab}$ enters from below then $P_p=(n+1,a-1)$;
\item if $\theta_p\simeq\phi_{ab}$ leaves to the right then $Q_p=(n+1,b)$;
\item if $\theta_p\simeq\phi_{ab}$ leaves upwards then $Q_p=(0,b)$.
\end{itemize}

The map from tableau $T$ of shape $F^{\lambda/\mu}$ to $s$-tuples of lattice paths between end points determined by
the outside decomposition $\Theta$ maps each entry $k$ in a box of content $c$
to a horizontal or diagonal edge terminating at its right hand end at $(k,c)$. Each path is completed so as to be continuous
through the addition of vertical edges, as illustrated below, where for the sake of clarity the directions 
of edges have been omitted:
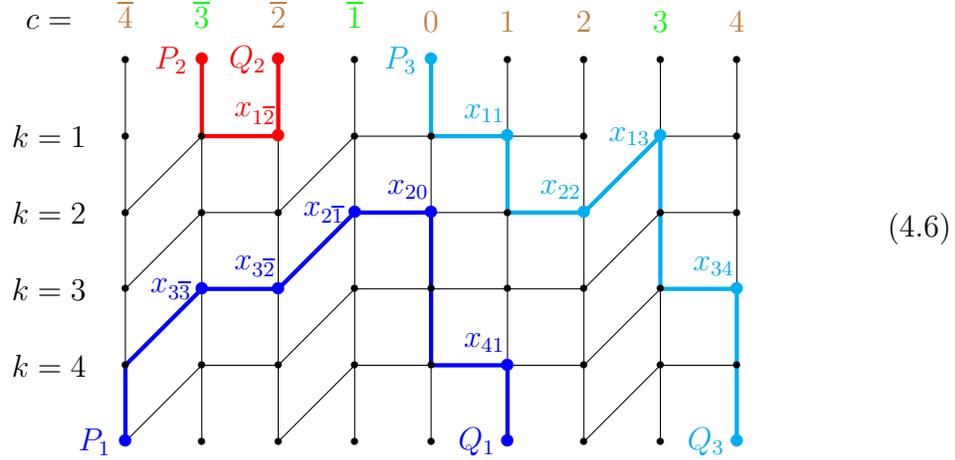
\begin{equation}
\vcenter{\hbox{
\begin{tikzpicture}[x={(0in,-0.4in)},y={(0.4in,0in)}] 
\foreach \i in {1,...,4} \draw(\i,-5)node{$k=\i$};
\draw(-0.5,-5)node{$c=$};
\foreach \j in {1,3} \draw(-0.6,-\j)node{$\green{\ov\j}$};
\foreach \j in {2,4} \draw(-0.6,-\j)node{$\brown{\ov\j}$};
\foreach \j in {3} \draw(-0.5,\j)node{$\green{\j}$};
\foreach \j in {0,1,2,4} \draw(-0.5,\j)node{$\brown{\j}$};
\foreach \i in {1,...,4} \foreach \j in {-2,0,1,2,4} \draw(\i,\j-1)--(\i,\j); 
\foreach \i in {1,...,4} \foreach \j in {-3,-1,3} \draw(\i+1,\j-1)--(\i,\j); 
\foreach \i in {1,...,5} \foreach \j in {-3,0,3,4} \draw(\i-1,\j)--(\i,\j);  
\foreach \i in {1,...,5} \foreach \j in {-4,-2,-1,1,2} \draw(\i,\j)--(\i-1,\j);  
\draw[draw=blue,ultra thick] (5,-4)--(4,-4)--(3,-3)--(3,-2)--(2,-1)--(2,0)--(4,0)--(4,1)--(5,1); 
\draw[draw=red,ultra thick] (0,-3)--(1,-3)--(1,-2)--(0,-2);
\draw[draw=cyan,ultra thick] (0,0)--(1,0)--(1,1)--(2,1)--(2,2)--(1,3)--(3,3)--(3,4)--(5,4); 
\foreach \i in {0,...,5} \foreach \j in {-4,...,4} \draw(\i,\j)node{$\ssc\bullet$};
\draw(5,-4.4)node{$\blue{P_1}$}; \draw(5,-4)node{$\blue\bullet$}; \draw(5,0.6)node{$\blue{Q_1}$}; \draw(5,1)node{$\blue\bullet$};
\draw(0,-3.4)node{$\red{P_2}$}; \draw(0,-3)node{$\red\bullet$}; \draw(0,-2.4)node{$\red{Q_2}$}; \draw(0,-2)node{$\red\bullet$};
\draw(0,-0.4)node{$\cyan{P_3}$}; \draw(0,0)node{$\cyan\bullet$}; \draw(5,3.6)node{$\cyan{Q_3}$}; \draw(5,4)node{$\cyan\bullet$};
\draw(1,3)node{$\cyan\bullet$};\draw(3,4)node{$\cyan\bullet$};
\draw(1,1)node{$\cyan\bullet$};\draw(2,2)node{$\cyan\bullet$};
\draw(1,-2)node{$\red\bullet$}; \draw(2,-1)node{$\blue\bullet$}; \draw(2,0)node{$\blue\bullet$}; \draw(4,1)node{$\blue\bullet$};
\draw(3,-3)node{$\blue\bullet$};\draw(3,-2)node{$\blue\bullet$};
\draw(1+0.0,3-0.4)node{$\cyan{x_{13}}$};\draw(3-0.3,4-0.3)node{$\cyan{x_{34}}$};
\draw(1-0.3,1-0.3)node{$\cyan{x_{11}}$};\draw(2-0.3,2-0.3)node{$\cyan{x_{22}}$};
\draw(1-0.3,-2-0.3)node{$\red{x_{1\ov2}}$}; \draw(2+0.0,-1-0.4)node{$\blue{x_{2\ov1}}$}; \draw(2-0.3,0-0.3)node{$\blue{x_{20}}$}; \draw(4-0.3,1-0.3)node{$\blue{x_{41}}$};
\draw(3+0.0,-3-0.4)node{$\blue{x_{3\ov3}}$};\draw(3-0.3,-2-0.3)node{$\blue{x_{3\ov2}}$};
\end{tikzpicture}
}}
\end{equation}

Hamel and Goulden~\cite{HG} have shown not only that this map provides a bijection between semistandard tableaux $T\in{\cal T}^{\lambda/\mu}$
and $s$-tuples of non-intersecting lattice paths $L\in{\cal L}^\phi$, but also that the image of any non-semistandard tableaux is
an $s$-tuple of intersecting lattice paths. By invoking the well-established link between non-intersecting lattice path models and determinants
provided by the Lindstr\"om-Gessel-Viennot Theorem~\cite{Lin,GV}, they completed their proof of Theorem~\ref{The-HG} appropriate 
to the classical Schur functions of (\ref{eqn-sfn}) by mapping the weight $x_k$ of each entry $k$ in $T$ to a corresponding weight $x_k$
of each horizontal or diagonal edge terminating at level $k$ in $L$, with vertical edges all carrying weight $1$. 
This then provides the special case of (\ref{eqn-HGX}) for which $x_{kc}=x_k$, that is to say when the weight of $t_{ij}$ is independent 
of the content $c$ of the box $(i,j)\in F^{\lambda/\mu}$. The extension to the case of weights $x_{kc}$ determined by an arbitrary choice of 
$\X=(x_{kc})$ is accomplished merely by noting that the bijection between tableaux and lattice paths must be made weight preserving in
the sense that the weight assigned to each directed edge is fixed to coincide with the weight of the corresponding tableau entry. 
This completes the proof of Theorem~\ref{The-HGX} once it is recognised that on the left-hand side the assignment of weights
depends on the content of the boxes of each tableaux contributing to $s_{\lambda/\mu}(\X)$ and that these weights are translated 
into weights assigned to the boxes of the tableaux contributing to $s_{(\theta_p\#\theta_q)}(\tau^{m(p,q)}\X)$ in which the content 
of each box is shifted under the outside decomposition from its usual value by the parameter $m(p,q)$. 
\qed

As a result it will be seen by examining the data of (\ref{eqn-pict}) and (\ref{eqn-mpq}) that Theorem~\ref{The-HGX} 
yields the generalised Schur function determinantal identity
\begin{equation}\label{eqn-sfndet-ex}
s_{5442/32}(\X)= \left|\, \begin{array}{ccc} 
      s_{42/1}(\tau^{-2}\X)&s_{2}(\tau^{-3}\X)&s_{652/41}(\tau^{-1}\X) \cr
       s_{31}(\tau^{-1}\X)&s_{1}(\tau^{-2}\X)&s_{541/3}(\X)\cr
			s_{1}(\tau\X)&0&s_{32/1}(\tau^2\X)\cr
\end{array}\,\right| 
\end{equation}
for all $\X=(x_{kc})$.

As claimed in the Introduction, Section \ref{Sec-introduction}, the non-intersecting lattice path model does not apply in the more general non-diagonal tenth generalisation case with a weight function $\X=(x_{k,(i,j)})$.
We illustrate this claim with the following example in which
the three alternative assignments of weights correspond to the first, ninth and tenth variation of Schur functions.
As can be seen, in each of the three pairs involving intersecting paths the weights are identical
in the left and right hand cases for weights of the first and ninth variations, but not the tenth.
This latter case therefore violates the hypothesis of the Lindstr\"om-Gessel-Viennot Theorem,
thereby precluding the possibility of deriving determinantal identities by this route. 
This same argument applies in the case of any skew Schur function and any outside decomposition.
This is because in every case of intersecting paths there is an ambiguity over the assignment of edges beyond 
the lowest rightmost point of intersection to one of two particular paths, and each of these edges will necessarily
carry a different weight $x_{k,(i,j)}$ depending on which path it lies, since its corresponding position
in the pair of distinct non-standard tableaux will have different coordinates $(i,j)$.

\begin{Example} Let $n=2$, $\lambda=(1,1)$ and $\mu=(0)$, for which there exists
a single semi-standard tableau and correspondingly a single pair of non-intersecting paths.
\begin{equation}
\begin{array}{ccc}
\vcenter{\hbox{
\begin{tikzpicture}[x={(0in,-0.2in)},y={(0.2in,0in)}] 
\draw (1,1) rectangle +(-1,-1);
\draw (2,1) rectangle +(-1,-1);
\draw(0.5,0.5)node{$1$};
\draw(1.5,0.5)node{$2$};
\end{tikzpicture}
}}
&
\vcenter{\hbox{
\begin{tikzpicture}[x={(0in,-0.4in)},y={(0.4in,0in)}] 
\draw(1,-1)node{$k=1$};\draw(2,-1)node{$k=2$};
\draw(-0.5,0)node{$c=$};\draw(-0.55,1)node{$\ov1$};\draw(-0.5,2)node{$0$};
\draw(1,-0.4)--(1,2.4); \draw(2,-0.4)--(2,2.4); 
\foreach \i in {1,...,3} \foreach \j in {0,...,2} \draw(\i-1,\j)--(\i,\j); 
\draw[draw=blue,ultra thick] (0,1)--(1,1)--(1,2)--(2,2)--(3,2); 
\draw[draw=red,ultra thick] (0,0)--(2,0)--(2,1)--(3,1);
\foreach \i in {0,...,3} \foreach \j in {0,1,2} \draw(\i,\j)node{$\ssc\bullet$};
\draw(0,0.6)node{$\blue{P_1}$}; \draw(0,-0.4)node{$\red{P_2}$}; \draw(3,2.4)node{$\blue{Q_1}$}; \draw(3,1.4)node{$\red{Q_2}$}; 
\draw(0,1)node{$\blue\bullet$};\draw(1,2)node{$\blue\bullet$};\draw(3,2)node{$\blue\bullet$};
\draw(0,0)node{$\red\bullet$};\draw(2,1)node{$\red\bullet$}; \draw(3,1)node{$\red\bullet$}; 
\draw(1-0.3,2-0.3)node{$\blue{B}$};\draw(2-0.3,1-0.3)node{$\red{C}$};
\end{tikzpicture}
}}
&
\begin{array}{cc}
x_k& x_1 x_2\cr
x_{kc}&x_{10}x_{2\ov1}\cr
x_{k,(i,j)}&x_{1,(1,1)} x_{2,(2,1)}\cr
\end{array}
\end{array}
\end{equation}

The remaining non-standard tableaux and intersecting pairs of paths and their corresponding weights take 
the following form:
\begin{equation}
\begin{array}{cccc}
\vcenter{\hbox{
\begin{tikzpicture}[x={(0in,-0.2in)},y={(0.2in,0in)}] 
\draw (1,1) rectangle +(-1,-1);
\draw (2,1) rectangle +(-1,-1);
\draw(0.5,0.5)node{$1$};
\draw(1.5,0.5)node{$1$};
\end{tikzpicture}
}}
&
\vcenter{\hbox{
\begin{tikzpicture}[x={(0in,-0.4in)},y={(0.4in,0in)}] 
\draw(1,-1)node{$k=1$};\draw(2,-1)node{$k=2$};
\draw(-0.5,0)node{$c=$};\draw(-0.55,1)node{$\ov1$};\draw(-0.5,2)node{$0$};
\draw(1,-0.4)--(1,2.4); \draw(2,-0.4)--(2,2.4); 
\foreach \i in {1,...,3} \foreach \j in {0,...,2} \draw(\i-1,\j)--(\i,\j); 
\draw[draw=blue,ultra thick] (0,1)--(1,1)--(1,2)--(2,2)--(3,2); 
\draw[draw=red,ultra thick] (0,0)--(1,0)--(1,1)--(2,1)--(3,1);
\foreach \i in {0,...,3} \foreach \j in {0,1,2} \draw(\i,\j)node{$\ssc\bullet$};
\draw(0,0.6)node{$\blue{P_1}$}; \draw(0,-0.4)node{$\red{P_2}$}; \draw(3,2.4)node{$\blue{Q_1}$}; \draw(3,1.4)node{$\red{Q_2}$}; 
\draw(0,1)node{$\blue\bullet$};\draw(1,2)node{$\blue\bullet$};\draw(3,2)node{$\blue\bullet$};
\draw(0,0)node{$\red\bullet$};\draw(1,1)node{$\red\bullet$}; \draw(3,1)node{$\red\bullet$}; 
\draw(1-0.3,1-0.3)node{$\red{A}$};\draw(1-0.3,2-0.3)node{$\blue{B}$};
\end{tikzpicture}
}}
&
\vcenter{\hbox{
\begin{tikzpicture}[x={(0in,-0.2in)},y={(0.2in,0in)}] 
\draw (0,0)--(1,0);
\draw (2,1) rectangle +(-1,-1);
\draw (2,2) rectangle +(-1,-1);
\draw(1.5,0.5)node{$1$};
\draw(1.5,1.5)node{$1$};
\end{tikzpicture}
}}
&
\vcenter{\hbox{
\begin{tikzpicture}[x={(0in,-0.4in)},y={(0.4in,0in)}] 
\draw(1,-1)node{$k=1$};\draw(2,-1)node{$k=2$};
\draw(-0.5,0)node{$c=$};\draw(-0.55,1)node{$\ov1$};\draw(-0.5,2)node{$0$};
\draw(1,-0.4)--(1,2.4); \draw(2,-0.4)--(2,2.4); 
\foreach \i in {1,...,3} \foreach \j in {0,...,2} \draw(\i-1,\j)--(\i,\j); 
\draw[draw=blue,ultra thick] (0,1)--(3,1); 
\draw[draw=red,ultra thick] (0,0)--(1,0)--(1,1)--(1,2)--(3,2);
\foreach \i in {0,...,3} \foreach \j in {0,1,2} \draw(\i,\j)node{$\ssc\bullet$};
\draw(0,0.6)node{$\blue{P_1}$}; \draw(0,-0.4)node{$\red{P_2}$}; \draw(3,2.4)node{$\red{Q_1}$}; \draw(3,1.4)node{$\blue{Q_2}$}; 
\draw(0,1)node{$\blue\bullet$};\draw(1,1)node{$\red\bullet$};\draw(3,1)node{$\blue\bullet$};
\draw(0,0)node{$\red\bullet$};\draw(1,2)node{$\red\bullet$}; \draw(3,2)node{$\red\bullet$}; 
\draw(1-0.3,1-0.3)node{$\red{A}$};\draw(1-0.3,2-0.3)node{$\red{B}$};
\end{tikzpicture}
}}
\cr\cr
&
\begin{array}{rl}
x_k& x_1 x_1\cr
x_{kc}&x_{1\ov1}x_{10}\cr
x_{k,(i,j)}&x_{1,(2,1)} x_{1,(1,1)}\cr
\end{array}
&&
\begin{array}{rl}
x_k& x_1 x_1\cr
x_{kc}&x_{1\ov1}x_{10}\cr
x_{k,(i,j)}&x_{1,(2,1)} x_{1,(2,2)}\cr
\end{array}
\end{array}
\end{equation}

\begin{equation}
\begin{array}{cccc}
\vcenter{\hbox{
\begin{tikzpicture}[x={(0in,-0.2in)},y={(0.2in,0in)}] 
\draw (1,1) rectangle +(-1,-1);
\draw (2,1) rectangle +(-1,-1);
\draw(0.5,0.5)node{$2$};
\draw(1.5,0.5)node{$2$};
\end{tikzpicture}
}}
&
\vcenter{\hbox{
\begin{tikzpicture}[x={(0in,-0.4in)},y={(0.4in,0in)}] 
\draw(1,-1)node{$k=1$};\draw(2,-1)node{$k=2$};
\draw(-0.5,0)node{$c=$};\draw(-0.55,1)node{$\ov1$};\draw(-0.5,2)node{$0$};
\draw(1,-0.4)--(1,2.4); \draw(2,-0.4)--(2,2.4); 
\foreach \i in {1,...,3} \foreach \j in {0,...,2} \draw(\i-1,\j)--(\i,\j); 
\draw[draw=blue,ultra thick] (0,1)--(2,1)--(2,2)--(3,2); 
\draw[draw=red,ultra thick] (0,0)--(2,0)--(2,1)--(3,1);
\foreach \i in {0,...,3} \foreach \j in {0,1,2} \draw(\i,\j)node{$\ssc\bullet$};
\draw(0,0.6)node{$\blue{P_1}$}; \draw(0,-0.4)node{$\red{P_2}$}; \draw(3,2.4)node{$\blue{Q_1}$}; \draw(3,1.4)node{$\red{Q_2}$}; 
\draw(0,1)node{$\blue\bullet$};\draw(2,2)node{$\blue\bullet$};\draw(3,2)node{$\blue\bullet$};
\draw(0,0)node{$\red\bullet$};\draw(2,1)node{$\red\bullet$}; \draw(3,1)node{$\red\bullet$}; 
\draw(2-0.3,1-0.3)node{$\red{C}$};\draw(2-0.3,2-0.3)node{$\blue{D}$};
\end{tikzpicture}
}}
&
\vcenter{\hbox{
\begin{tikzpicture}[x={(0in,-0.2in)},y={(0.2in,0in)}] 
\draw (0,0)--(1,0);
\draw (2,1) rectangle +(-1,-1);
\draw (2,2) rectangle +(-1,-1);
\draw(1.5,0.5)node{$2$};
\draw(1.5,1.5)node{$2$};
\end{tikzpicture}
}}
&
\vcenter{\hbox{
\begin{tikzpicture}[x={(0in,-0.4in)},y={(0.4in,0in)}] 
\draw(1,-1)node{$k=1$};\draw(2,-1)node{$k=2$};
\draw(-0.5,0)node{$c=$};\draw(-0.55,1)node{$\ov1$};\draw(-0.5,2)node{$0$};
\draw(1,-0.4)--(1,2.4); \draw(2,-0.4)--(2,2.4); 
\foreach \i in {1,...,3} \foreach \j in {0,...,2} \draw(\i-1,\j)--(\i,\j); 
\draw[draw=blue,ultra thick] (0,1)--(3,1); 
\draw[draw=red,ultra thick] (0,0)--(2,0)--(2,2)--(3,2);
\foreach \i in {0,...,3} \foreach \j in {0,1,2} \draw(\i,\j)node{$\ssc\bullet$};
\draw(0,0.6)node{$\blue{P_1}$}; \draw(0,-0.4)node{$\red{P_2}$}; \draw(3,2.4)node{$\red{Q_1}$}; \draw(3,1.4)node{$\blue{Q_2}$}; 
\draw(0,1)node{$\blue\bullet$};\draw(3,1)node{$\blue\bullet$};
\draw(0,0)node{$\red\bullet$};\draw(2,1)node{$\red\bullet$};\draw(2,2)node{$\red\bullet$}; \draw(3,2)node{$\red\bullet$}; 
\draw(2-0.3,1-0.3)node{$\red{C}$};\draw(2-0.3,2-0.3)node{$\red{D}$};
\end{tikzpicture}
}}
\cr\cr
&
\begin{array}{rl}
x_k& x_2 x_2\cr
x_{kc}&x_{2\ov1}x_{20}\cr
x_{k,(i,j)}&x_{2,(2,1)} x_{2,(1,1)}\cr
\end{array}
&&
\begin{array}{rl}
x_k& x_2 x_2\cr
x_{kc}&x_{2\ov1}x_{20}\cr
x_{k,(i,j)}&x_{2,(2,1)} x_{1,(2,2)}\cr
\end{array}
\end{array}
\end{equation}

\begin{equation}
\begin{array}{cccc}
\vcenter{\hbox{
\begin{tikzpicture}[x={(0in,-0.2in)},y={(0.2in,0in)}] 
\draw (1,1) rectangle +(-1,-1);
\draw (2,1) rectangle +(-1,-1);
\draw(0.5,0.5)node{$2$};
\draw(1.5,0.5)node{$1$};
\end{tikzpicture}
}}
&
\vcenter{\hbox{
\begin{tikzpicture}[x={(0in,-0.4in)},y={(0.4in,0in)}] 
\draw(1,-1)node{$k=1$};\draw(2,-1)node{$k=2$};
\draw(-0.5,0)node{$c=$};\draw(-0.55,1)node{$\ov1$};\draw(-0.5,2)node{$0$};
\draw(1,-0.4)--(1,2.4); \draw(2,-0.4)--(2,2.4); 
\foreach \i in {1,...,3} \foreach \j in {0,...,2} \draw(\i-1,\j)--(\i,\j); 
\draw[draw=blue,ultra thick] (0,1)--(2,1)--(2,2)--(3,2); 
\draw[draw=red,ultra thick] (0,0)--(1,0)--(1,1)--(2,1)--(3,1);
\draw[draw=magenta,ultra thick] (1,1)--(2,1);
\foreach \i in {0,...,3} \foreach \j in {0,1,2} \draw(\i,\j)node{$\ssc\bullet$};
\draw(0,0.6)node{$\blue{P_1}$}; \draw(0,-0.4)node{$\red{P_2}$}; \draw(3,2.4)node{$\blue{Q_1}$}; \draw(3,1.4)node{$\red{Q_2}$}; 
\draw(0,1)node{$\blue\bullet$};\draw(2,2)node{$\blue\bullet$};\draw(3,2)node{$\blue\bullet$};
\draw(0,0)node{$\red\bullet$};\draw(1,1)node{$\red\bullet$}; \draw(3,1)node{$\red\bullet$}; 
\draw(2-0.3,2-0.3)node{$\blue{D}$};\draw(1-0.3,1-0.3)node{$\red{B}$};
\end{tikzpicture}
}}
&
\vcenter{\hbox{
\begin{tikzpicture}[x={(0in,-0.2in)},y={(0.2in,0in)}] 
\draw (0,0)--(1,0);
\draw (2,1) rectangle +(-1,-1);
\draw (2,2) rectangle +(-1,-1);
\draw(1.5,0.5)node{$1$};
\draw(1.5,1.5)node{$2$};
\end{tikzpicture}
}}
&
\vcenter{\hbox{
\begin{tikzpicture}[x={(0in,-0.4in)},y={(0.4in,0in)}] 
\draw(1,-1)node{$k=1$};\draw(2,-1)node{$k=2$};
\draw(-0.5,0)node{$c=$};\draw(-0.55,1)node{$\ov1$};\draw(-0.5,2)node{$0$};
\draw(1,-0.4)--(1,2.4); \draw(2,-0.4)--(2,2.4); 
\foreach \i in {1,...,3} \foreach \j in {0,...,2} \draw(\i-1,\j)--(\i,\j); 
\draw[draw=blue,ultra thick] (0,1)--(3,1); 
\draw[draw=red,ultra thick] (0,0)--(1,0)--(1,1)--(2,1)--(2,2)--(3,2);
\draw[draw=magenta,ultra thick] (1,1)--(2,1);
\foreach \i in {0,...,3} \foreach \j in {0,1,2} \draw(\i,\j)node{$\ssc\bullet$};
\draw(0,0.6)node{$\blue{P_1}$}; \draw(0,-0.4)node{$\red{P_2}$}; \draw(3,2.4)node{$\red{Q_1}$}; \draw(3,1.4)node{$\blue{Q_2}$}; 
\draw(0,1)node{$\blue\bullet$};\draw(1,1)node{$\red\bullet$};\draw(3,1)node{$\blue\bullet$};
\draw(0,0)node{$\red\bullet$};\draw(2,2)node{$\red\bullet$}; \draw(3,2)node{$\red\bullet$}; 
\draw(2-0.3,2-0.3)node{$\red{D}$};\draw(1-0.3,1-0.3)node{$\red{B}$};
\end{tikzpicture}
}}
\cr\cr
&
\begin{array}{rl}
x_k& x_1 x_2\cr
x_{kc}&x_{1\ov1}x_{20}\cr
x_{k,(i,j)}&x_{1,(2,1)} x_{2,(1,1)}\cr
\end{array}
&&
\begin{array}{rl}
x_k& x_1 x_2\cr
x_{kc}&x_{1\ov1}x_{20}\cr
x_{k,(i,j)}&x_{1,(2,1)} x_{2,(2,2)}\cr
\end{array}
\end{array}
\end{equation}

\end{Example}


\subsection{Generalisation of Hamel's Theorem~\ref{The-Ham}}\label{subsec HamX}

\begin{Theorem}\label{The-HamXY}
Let $\lambda$ and $\mu$ be strict partitions with $\mu\subseteq\lambda$, and let $d$ be
the number of boxes on the main diagonal of the skew shifted diagram $F^{\lambda/\mu}$. 
For a given cutting strip $\phi$, let $\theta=(\theta_1,\theta_2,\ldots,\theta_r,\theta_{r+1},\ldots,\theta_s)$ be an outside decomposition
of $F^{\lambda/\mu}$, where $\theta_p$ includes a box on the main diagonal of $F^{\lambda/\mu}$ if and only 
if $1\leq i\leq d$, with $r=d$ if $d$ is even and $r=d+1$ if $d$ is odd, in which case $\theta_r$ is a null strip.
Let $\rho$ be identified with the unique box of the cutting strip $\phi$ having content $0$, and let
$\ov{\theta}_p=(\rho\#\theta_p)$ for all $p=1,2,\ldots,s$. 
Let $\tau$ be a shift operator whose action on $(\X,\!\Y)=((x_{kc}),(y_{lc}))$,with $k,l\in\{1,2,\ldots,n\}$ and $c\in\Z$,
is such that $\tau^m(\X,\!\Y)=((x_{k,m+c}),(y_{l,m+c}))$ for any $m\in\Z$
Then the generalised $Q$-functions of (\ref{eqn-QfnXY}) satisfy the identity
\begin{eqnarray}\label{eqn-HamXY}
    \lefteqn{Q_{\lambda/\mu}(\X,\!\Y)=}\\ \nonumber
   & &  \pf\left( \begin{array}{cc} Q_{(\ov{\theta}_p,\ov{\theta}_q)}(\X,\!\Y)&Q_{(\theta_{s+r+1-k}\#\theta_p)}(\tau^{m(s+r+1-k,p)}(\X,\!\Y))\cr
		                                      -\,^t(Q_{(\theta_{s+r+1-k}\#\theta_p)}(\tau^{m(s+r+1-k,p)}(\X,\!\Y)))&0\cr \end{array} \right) 
\end{eqnarray}
with $1\leq p,q\leq s$ and $r+1\leq k\leq s$, where $(\ov{\theta}_p,\ov{\theta}_q)$ 
is the skew shifted shape formed by abutting $\ov{\theta}_p$ and $\ov{\theta}_q$ with their 
boxes of content $0$ lying on the main diagonal, with that of $\ov{\theta}_p$ immediately above 
and to the left of that of $\ov{\theta}_q$.
Here
\begin{equation}\label{eqn-Qpq}
    Q_{(\ov{\theta}_p,\ov{\theta}_q)}(\X,\!\Y)=- Q_{(\ov{\theta}_q,\ov{\theta}_p)}(\X,\!\Y)
    \qquad \mbox{with}\qquad    Q_{(\ov{\theta}_p,\ov{\theta}_p)}(\X,\Y)=0,
\end{equation}
and if $(\theta_{s+r+1-k}\#\theta_p)=\phi_{ad}$ then
\begin{equation}\label{eqn-Qkp}
       Q_{(\theta_{s+r+1-k}\#\theta_p)}(\tau^{m(s+r+1-k,p)}(\X,\!\Y))= \begin{cases}
                                                    Q_{\phi_{ad}}(\tau^{\a}(\X,\!\Y))&\mbox{if $a\leq d$};\cr
                            1&\mbox{if $a=d+1$};\cr
                            0&\mbox{if $a>d+1$}.\cr
                                                \end{cases}
\end{equation}

\end{Theorem}

\noindent{\bf Proof}: Very much as in the case of Theorem~\ref{The-HG}, in order to prove Theorem~\ref{The-Ham} 
Hamel~\cite{Ham} established a bijection between
standard skew shifted primed tableaux $T\in{\cal ST}^{\lambda/\mu}$ and $L\in{\cal L}^\phi$, the set of $s$-tuples of non-intersecting 
paths on a rectangular lattice with directed edges determined by the cutting strip $\phi$ and paths determined by the entries 
of $T$ in the constituent strips $\theta_p$ of the corresponding outside decomposition $\Theta=(\theta_1,\theta_2,\ldots,\theta_s)$ 
of $F^{\lambda/\mu}$. She showed furthermore that this same map takes non-semistandard skew shifted primed tableaux to $s$-tuples of paths involving 
at least one intersection.  

The use of Hamel's bijection can be exemplified as follows. 
In the case $\lambda=(9,6,4,2)$ and $\mu=(4,3)$ a typical primed shifted skew tableau $T\in{\cal ST}^{\lambda/\mu}$ with entries 
$k$ or $k'$ from $\{1'<1<2'<2<3'<3<4'<4\}$ is as shown below on the right, with the content $c$ of the boxes of $SF^{\lambda/\mu}$ 
shown on the left.
\begin{equation}\label{eqn-SF-T}
\vcenter{\hbox{
\begin{tikzpicture}[x={(0in,-0.2in)},y={(0.2in,0in)}] 
\foreach \j in {1,...,9} \draw[thick] (1,\j) rectangle +(-1,-1);
\foreach \j in {2,...,7} \draw[thick] (2,\j) rectangle +(-1,-1);
\foreach \j in {3,...,6} \draw[thick] (3,\j) rectangle +(-1,-1);
\foreach \j in {4,...,5} \draw[thick] (4,\j) rectangle +(-1,-1);
\foreach \j in {1,...,4} \draw(0.5,\j-0.5)node{$\ast$};
\foreach \j in {2,...,4} \draw(1.5,\j-0.5)node{$\ast$};
\foreach \j in {4,...,8} \draw(1-0.5,\j+0.5)node{$\j$}; 
\foreach \j in {3,...,5} \draw(2-0.5,\j+0.5+1)node{$\j$};
\foreach \j in {0,...,3} \draw(3-0.5,\j+0.5+2)node{$\j$};
\foreach \j in {0,...,1} \draw(4-0.5,\j+0.5+3)node{$\j$};
\end{tikzpicture}
}}
\qquad\qquad
\vcenter{\hbox{
\begin{tikzpicture}[x={(0in,-0.2in)},y={(0.2in,0in)}] 
\foreach \j in {1,...,9} \draw[thick]  (1,\j) rectangle +(-1,-1);
\foreach \j in {2,...,7} \draw[thick]  (2,\j) rectangle +(-1,-1);
\foreach \j in {3,...,6} \draw[thick]  (3,\j) rectangle +(-1,-1);
\foreach \j in {4,...,5} \draw[thick]  (4,\j) rectangle +(-1,-1);
\foreach \j in {1,...,4} \draw(0.5,\j-0.5)node{$\ast$};
\foreach \j in {2,...,4} \draw(1.5,\j-0.5)node{$\ast$};
\draw (1-0.5,4+0.5)node{$\blue{1}$};\draw (1-0.5,5+0.5)node{$\blue{2'}$};\draw (1-0.5,6+0.5)node{$\blue{2}$};\draw (1-0.5,7+0.5)node{$\blue{4'}$};\draw (1-0.5,8+0.5)node{$\blue{4}$}; 
\draw (2-0.5,3+0.5+1)node{$\blue{3'}$}; \draw (2-0.5,4+0.5+1)node{$\cyan{3}$};\draw (2-0.5,5+0.5+1)node{$\cyan{3}$};
\draw (3-0.5,0+0.5+2)node{$\blue{2}$}; \draw (3-0.5,1+0.5+2)node{$\blue{2}$};\draw (3-0.5,2+0.5+2)node{$\blue{3'}$};\draw (3-0.5,3+0.5+2)node{$\cyan{4}$};
\draw (4-0.5,1+0.5+2)node{$\red{4'}$}; \draw (4-0.5,2+0.5+2)node{$\red{4}$};
\end{tikzpicture}
}}
\end{equation}
The required lattice, which differs somewhat from that of~\cite{Ham}, is set up on rectangular grid with directed edges 
determined by the choice of cutting strip. In the case of the cutting strip
\begin{equation}\label{eqn-csQex}
\vcenter{\hbox{
\begin{tikzpicture}[x={(0in,-0.2in)},y={(0.2in,0in)}] 
\draw[thick](1,9) rectangle +(-1,-1);\draw(1-0.5,9-0.5)node{$\brown8$};
\draw[thick](1,8) rectangle +(-1,-1);\draw(1-0.5,8-0.5)node{$\brown7$};
\draw[thick](1,7) rectangle +(-1,-1);\draw(1-0.5,7-0.5)node{$\brown6$};
\draw[thick](1,6) rectangle +(-1,-1);\draw(1-0.5,6-0.5)node{$\brown5$};
\draw[thick](1,5) rectangle +(-1,-1);\draw(1-0.5,5-0.5)node{$\green4$};
\draw[thick](2,5) rectangle +(-1,-1);\draw(2-0.5,5-0.5)node{$\green3$};
\draw[thick](3,5) rectangle +(-1,-1);\draw(3-0.5,5-0.5)node{$\brown2$};
\draw[thick](3,4) rectangle +(-1,-1);\draw(3-0.5,4-0.5)node{$\brown1$};
\draw[thick](3,3) rectangle +(-1,-1);\draw(3-0.5,3-0.5)node{$\brown0$};
\end{tikzpicture}
}}
\end{equation}
and entries no greater than $n=4$ the corresponding lattice is given by
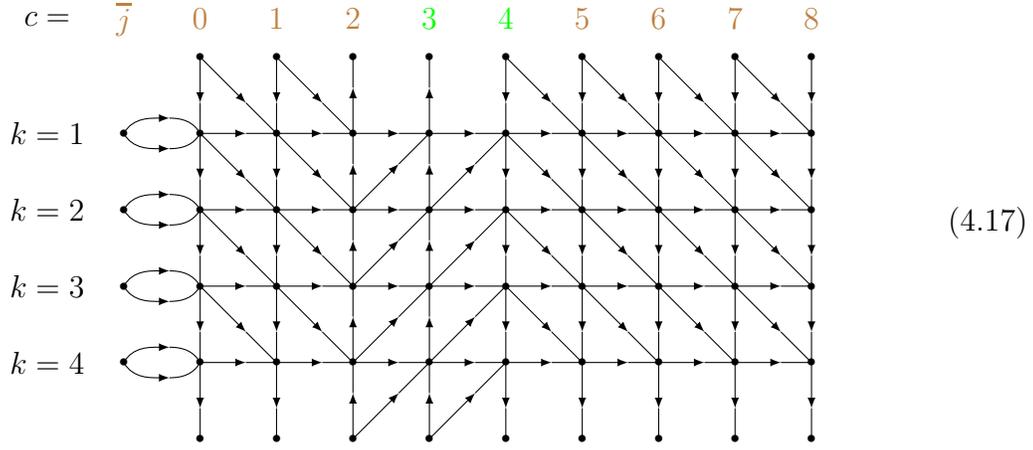
\begin{equation}\label{eqn-LQ}
\vcenter{\hbox{
\begin{tikzpicture}[x={(0in,-0.4in)},y={(0.4in,0in)}] 
\foreach \i in {1,...,4} \draw(\i,-2)node{$k=\i$};
\draw(-0.5,-2)node{$c=$};
\foreach \j in {1} \draw(-0.5,-\j)node{$\brown{\ov{j}}$};
\foreach \j in {0,...,2} \draw(-0.5,\j)node{$\brown{\j}$};
\foreach \j in {3,4} \draw(-0.5,\j)node{$\green{\j}$};
\foreach \j in {5,...,8} \draw(-0.5,\j)node{$\brown{\j}$};
\foreach \i in {1,...,4} \foreach \j in {1,...,8} \draw[-latex](\i,\j-1)--(\i,\j-0.4);
\foreach \i in {1,...,4} \foreach \j in {1,...,8} \draw(\i,\j-0.4)--(\i,\j);
\foreach \i in {1,...,4} \foreach \j in {1,2,5,6,7,8} \draw[-latex](\i-1,\j-1)--(\i-0.4,\j-0.4);
\foreach \i in {1,...,4} \foreach \j in {1,2,5,6,7,8} \draw(\i-0.4,\j-0.4)--(\i,\j);
\foreach \i in {1,...,4} \foreach \j in {3,4} \draw[-latex](\i+1,\j-1)--(\i+0.4,\j-0.4);
\foreach \i in {1,...,4} \foreach \j in {3,4} \draw(\i+0.4,\j-0.4)--(\i,\j);
\foreach \j in {0,1,4,5,6,7,8} \foreach \i in {1,...,5} \draw[-latex](\i-1,\j)--(\i-0.4,\j);
\foreach \j in {0,1,4,5,6,7,8} \foreach \i in {1,...,5}\draw(\i-0.4,\j)--(\i,\j);
\foreach \j in {2,3} \foreach \i in {0,...,4} \draw[-latex](\i+1,\j)--(\i+0.4,\j);
\foreach \j in {2,3} \foreach \i in {0,...,4} \draw(\i+0.4,\j)--(\i,\j);
\foreach \i in {1,...,4} \draw[-latex](\i,-1)to[out=60,in=180](\i-0.2,-0.4);
\foreach \i in {1,...,4} \draw(\i-0.2,-0.4)to[out=0,in=120](\i,0);
\foreach \i in {1,...,4} \draw[-latex](\i,-1)to[out=-60,in=-180](\i+0.2,-0.4);
\foreach \i in {1,...,4} \draw(\i+0.2,-0.4)to[out=0,in=-120](\i,0);
\foreach \i in {0,...,5} \foreach \j in {0,...,8} \draw(\i,\j)node{$\ssc\bullet$};
\foreach \i in {1,...,4} \draw(\i,-1)node{$\ssc\bullet$};
\end{tikzpicture}
}}
\end{equation}
The orientation of horizontal edges is always left to right, while diagonal edges are either upwards or downwards
to the right according as their end point is in column $c$ and the box of content $c$ in $\phi$
is approached from either below or the left, respectively, in reading $\phi$ as usual from bottom left to top right.
The curved edges at level $k$ from $c=\ov1$ to $c=0$ are to accommodate
entries $k$ and $k'$ that are allowed in boxes of content $0$ of $T$.

The start and end points, $P_p$ and $Q_p$, of the lattice paths corresponding to each component $\theta_p$ of the outside 
decomposition are determined as follows:
\begin{itemize}
\item if $\theta_p\simeq\phi_{0b}$ enters from the left then $P_p=(k,-1)$ for some $k\in\{1,2,\ldots,n\}$;
\item if $\theta_p\simeq\phi_{ab}$ enters from the left and $a>0$ then $P_p=(0,a-1)$;
\item if $\theta_p\simeq\phi_{ab}$ enters from below then $P_p=(n+1,a-1)$;
\item if $\theta_p\simeq\phi_{ab}$ leaves to the right then $Q_p=(n+1,b)$;
\item if $\theta_p\simeq\phi_{ab}$ leaves upwards then $Q_p=(0,b)$.
\end{itemize}

An example of a $3$-tuple of non-intersecting paths may be obtained by noting the entries $k$ or $k'$ in the 
the primed shifted tableaux $T$ of (\ref{eqn-SF-T}) in the boxes of $\theta_1$, $\theta_2$ and $\theta_3$
of content $c$ as specified in (\ref{eqn-stripsQ}). This takes the form
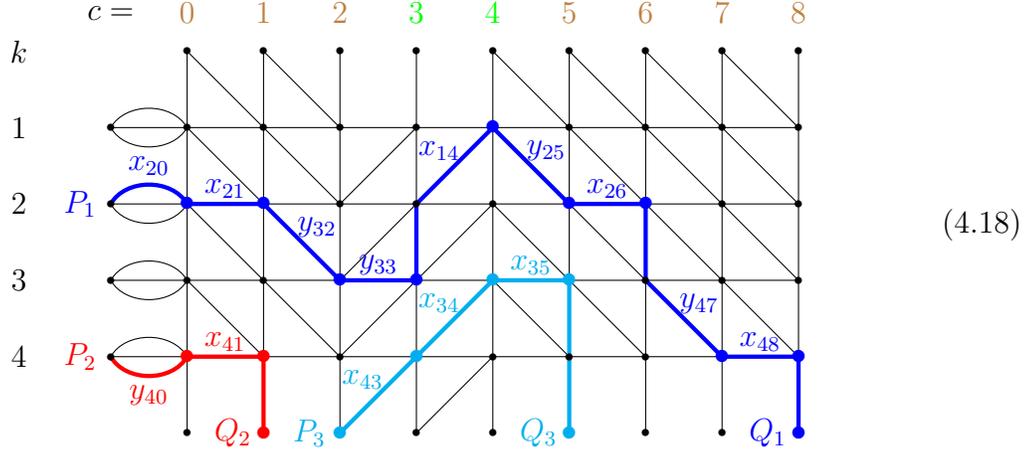
\begin{equation}
\vcenter{\hbox{
\begin{tikzpicture}[x={(0in,-0.4in)},y={(0.4in,0in)}] 
\draw(0,-2.2)node{$k$};
\foreach \k in {1,...,4} \draw(\k,-2.2)node{$\k$};
\draw(-0.5,-1)node{$c=$};
\foreach \j in {0,...,2} \draw(-0.5,\j)node{$\brown{\j}$};
\foreach \j in {3,4} \draw(-0.5,\j)node{$\green{\j}$};
\foreach \j in {5,...,8} \draw(-0.5,\j)node{$\brown{\j}$};
\foreach \i in {1,...,4} \foreach \j in {0,...,8} \draw(\i,\j-1)--(\i,\j);
\foreach \i in {1,...,4} \foreach \j in {1,2,5,6,7,8} \draw(\i-1,\j-1)--(\i,\j);
\foreach \i in {1,...,4} \foreach \j in {3,4} \draw(\i+1,\j-1)--(\i,\j);
\foreach \j in {0,...,8} \draw(0,\j)--(5,\j);
\foreach \i in {1,...,4} \draw(\i,-1)to[out=60,in=120](\i,0);
\foreach \i in {1,...,4} \draw(\i,-1)to[out=-60,in=-120](\i,0);
\draw[draw=blue,ultra thick] (2,-1)to[out=60,in=120](2,0);
\draw[draw=blue,ultra thick] (2,-0)--(2,1)--(3,2)--(3,3)--(2,3)--(1,4)--(2,5)--(2,6)--(3,6)--((4,7)--(4,8)--(5,8);
\draw[draw=red,ultra thick] (4,-1)to[out=-60,in=-120](4,0); 
\draw[draw=red,ultra thick] (4,0)--(4,1)--(5,1);
\draw[draw=cyan,ultra thick] (5,2)--(4,3)--(3,4)--(3,5)--(5,5);
\foreach \i in {0,...,5} \foreach \j in {0,...,8} \draw(\i,\j)node{$\ssc\bullet$};
\foreach \i in {1,...,4} \draw(\i,-1)node{$\ssc\bullet$};
\draw(2,0)node{$\blue\bullet$};\draw(2,1)node{$\blue\bullet$};\draw(3,2)node{$\blue\bullet$};\draw(3,3)node{$\blue\bullet$};
\draw(1,4)node{$\blue\bullet$};\draw(2,5)node{$\blue\bullet$};\draw(2,6)node{$\blue\bullet$};\draw(4,7)node{$\blue\bullet$};
\draw(4,8)node{$\blue\bullet$};
\draw(4,0)node{$\red\bullet$};\draw(4,1)node{$\red\bullet$};
\draw(4,3)node{$\cyan\bullet$}; \draw(3,4)node{$\cyan\bullet$}; \draw(3,5)node{$\cyan\bullet$};
\draw(1.5,-0.5)node{$\blue{x_{20}}$};\draw(1.5+0.3,0.5)node{$\blue{x_{21}}$};\draw(2.5-0.2,1.5+0.2)node{$\blue{y_{32}}$};\draw(2.5+0.3,2.5)node{$\blue{y_{33}}$};
\draw(1.5-0.2,3.5-0.2)node{$\blue{x_{14}}$};\draw(1.5-0.2,4.5+0.2)node{$\blue{y_{25}}$};\draw(1.5+0.3,5.5)node{$\blue{x_{26}}$};\draw(3.5-0.2,6.5+0.2)node{$\blue{y_{47}}$};\draw(3.5+0.3,7.5)node{$\blue{x_{48}}$};
\draw(4.5-0.2,2.5-0.2)node{$\cyan{x_{43}}$};\draw(3.5-0.2,3.5-0.2)node{$\cyan{x_{34}}$};\draw(2.5+0.3,4.5)node{$\cyan{x_{35}}$};
\draw(4.5,-0.5)node{$\red{y_{40}}$};\draw(3.5+0.3,0.5)node{$\red{x_{41}}$};
\draw(2,-1.4)node{$\blue{P_1}$};
\draw(4,-1.4)node{$\red{P_2}$};
\draw(5,7.6)node{$\blue{Q_1}$};\draw(5,8)node{$\blue\bullet$};
\draw(5,0.6)node{$\red{Q_2}$};\draw(5,1)node{$\red\bullet$};
\draw(5,1.6)node{$\cyan{P_3}$};\draw(5,2)node{$\cyan\bullet$}; 
\draw(5,4.6)node{$\cyan{Q_3}$};\draw(5,5)node{$\cyan\bullet$};
\end{tikzpicture}
}}
\end{equation}
For given $k$ or $k'$ and $c$ the directed edges take one or other of the following forms, each carrying the indicated weight $x_{kc}$ or $y_{kc}$: 
\begin{equation}
\begin{array}{|c|c|ccc|}
\hline
\mbox{Entry}&\mbox{Weight}&c=0&c\in\{1,2,5,6,7,8\}&c\in\{3,4\}\cr
\hline
&&&&\cr
k&x_{kc}
&\vcenter{\hbox{
\begin{tikzpicture}[x={(0in,-0.4in)},y={(0.4in,0in)}] 
\draw[-latex,draw=cyan,ultra thick](7,-1)to[out=60,in=180](7-0.2,-0.3); \draw[draw=cyan,ultra thick](7-0.2,-0.4)to[out=0,in=120](7,0); 
\end{tikzpicture}
}}
&\vcenter{\hbox{
\begin{tikzpicture}[x={(0in,-0.4in)},y={(0.4in,0in)}] 
\draw[-latex,draw=blue,ultra thick] (7,0.5)--(7,0.5+0.8);\draw[draw=blue,ultra thick] (7,0.5+0.6)--(7,1.5);
\end{tikzpicture}
}}
&\vcenter{\hbox{
\begin{tikzpicture}[x={(0in,-0.4in)},y={(0.4in,0in)}] 
\draw[-latex,draw=red,ultra thick] (7.5,2.5)--(7.5-0.7,2.5+0.7);\draw[draw=red,ultra thick] (7.5-0.5,2.5+0.5)--(6.5,3.5); 
\end{tikzpicture}
}}
\cr
k'&y_{kc}
&\vcenter{\hbox{
\begin{tikzpicture}[x={(0in,-0.4in)},y={(0.4in,0in)}] 
\draw[-latex,draw=cyan,ultra thick](8,-1)to[out=-60,in=-180](8+0.2,-0.3); \draw[draw=cyan,ultra thick](8+0.2,-0.4)to[out=0,in=-120](8,0); 
\end{tikzpicture}
}}
&\vcenter{\hbox{
\begin{tikzpicture}[x={(0in,-0.4in)},y={(0.4in,0in)}] 
\draw[-latex,draw=blue,ultra thick] (7.5,0.5)--(7.5+0.7,0.5+0.7); \draw[draw=blue,ultra thick] (7.5+0.5,0.5+0.5)--(8.5,1.5); 
\end{tikzpicture}
}}
&\vcenter{\hbox{
\begin{tikzpicture}[x={(0in,-0.4in)},y={(0.4in,0in)}] 
\draw[-latex,draw=red,ultra thick] (8,2.5)--(8,2.5+0.8); \draw[draw=red,ultra thick] (8,2.5+0.6)--(8,3.5);
\end{tikzpicture}
}}
\cr
&&&&\cr
\hline
\end{array}
\end{equation}

Although this lattice path environment differs from that of Hamel~\cite{Ham} it can be easily obtained by shrinking to zero length 
all vertical edges between rows of Hamel's lattice labelled by $k$ and $k'$ for all $k$, converting horizontal and diagonal
edges terminating in column $c=0$ to our concave upwards and concave downwards edges, respectively, and turning the whole 
lattice upside down. That our standard primed shifted skew tableaux $T\in{\cal T}^{\lambda/\mu}$ are in bijective correspondence
with $s$-tuples of non-intersecting paths $L\in{\cal L}^\phi$, the rectangular lattice determined by the cutting strip $\phi$,
follows from the arguments analogous to those used in~\cite{Ham}, as does her conclusion that under the same map non-standard 
tableaux lead to intersecting paths. The proof of Theorem~\ref{The-Ham} then follows from Stembridge's Theorem~6.2~\cite{Ste90}
obtained through his generalisation of the Lindstr\"om-Gessel-Viennot non-intersecting path model~\cite{Lin,GV}, 
together with the fact that as in the skew Schur function case, for any $\X=(x_{kc})$ and $\Y=(y_{kc})$,  it is necessary to apply 
the shift operator $\tau^{m(k,p)}$ to the content parameter $c$ of all terms contributing to the $Q$-functions of single strips 
$(\theta_p\#\theta_k)$ that appear on the right of (\ref{eqn-HamXY}).
\qed

As an application of Theorem~\ref{The-Ham}, consider the case $\lambda=(9,6,4,2)$ and $\mu=(4,3)$ with the cutting strip $\phi$
and outside decomposition as in (\ref{eqn-csQ}), so that $r=d=2$ and $s=3$. The data from (\ref{eqn-ovp}) and (\ref{eqn-ovpq}), 
together with the observation that the shift parameters $m(3,p)=3$ for $p=1,3$, implies the validity of the following Pfaffian 
identity where use has also been made of both (\ref{eqn-Qpq}) and (\ref{eqn-Qkp}) to deal with non-standard shapes:
\begin{eqnarray}
     \lefteqn{ Q_{9642/43}(\X,\!\Y) =} \\ \nonumber
     & & \left[
             \begin{array}{cccc}
						      0&Q_{9432/43}(\X,\!\Y)&Q_{9643/43}(\X,\!\Y)&Q_{61/1}(\tau^3(\X,\!\Y))\cr
								-Q_{9432/43}(\X,\!\Y)	&0&-Q_{6432/43}(\X,\!\Y)&0\cr
								-Q_{9643/43}(\X,\!\Y)	&Q_{6432}(\X,\!\Y)&0&Q_{31/1}(\tau^3(\X,\!\Y))\cr
								 -Q_{61/1}(\tau^3(\X,\!\Y))&0&-Q_{31/1}(\tau^3(\X,\!\Y))&0\cr
						 \end{array}
		\right]\,.
\end{eqnarray}

\section{Corollaries}\label{Sec-cor}

To establish contact with the historical results referred to in Section~\ref{Sec-previous-work} and to generalise
each of them to the case of parameters $\X=(x_{kc})$ and $\Y=(y_{kc})$ it is only necessary to identify 
the appropriate outside decomposition and use it in Theorems~\ref{The-HGX} and~\ref{The-HamXY}. 
These outside decompositions are all characterised by the shape of the corresponding cutting strip $\phi$.
The classical identities are associated with the shape of $\phi$ being a row, a column or a hook, and, as we shall see,
more recent identities arise from a ``shifted hook'' or a consideration of the cases where $\phi$ takes the shape of the
inner or outer rim of $F^{\lambda/\mu}$ or $SF^{\lambda/\mu}$ as appropriate. 

\subsection{Determinantal corollaries}\label{Subsec-det-cor}

Before embarking on the study of particular cases it is worth remarking that the components $\theta_p$
of an outside decomposition $\Theta$ generated by some cutting strip $\phi$ are significantly
constrained by the very fact that $\Theta$ is an outside decomposition of $F^{\lambda/\mu}$. 
Whatever its shape, the cutting strip $\phi$ extends from a box of content $-l+1$ to one of content
$l'-1$, where $l$ and $l'$ are the number of boxes in the first column and first row of $F^\lambda$, respectively.
Each strip $\theta_r\simeq\phi_{ab}$ has starting and ending boxes of content $a$ and $b$, 
respectively, on either the inner or outer rim of the skew diagram $F^{\lambda/\mu}$. 
This means that their content must take one or other of the forms $\lambda_i-i$, $-\lambda'_j+j$, $\mu_i-i+1$, or $-\mu'_j+j-1$ 
for $i=1,2,\ldots,l$ and $j=1,2,\ldots,l'$,  with $\mu$ and $\mu'$ extended if necessary to include trailing zeros. 
In Frobenius notation, if $\lambda=(\alpha|\beta)$ and $\mu=(\gamma|\delta)$ have Frobenius rank $p$ and $q$ respectively
these forms include $\alpha_i$, $-\beta_i$, 
$\gamma_j+1$ or $-\delta_j-1$ for $i=1,2,\ldots,p$ and $j=1,2,\ldots,q$. 

The outside decompositions of $F^{\lambda/\mu}$ in the case $\lambda=(5,4,4,2)=(4,2,1\, |\, 3,2,0)$ and $\mu=(3,2)=(2,0\, |\, 1,0)$
are illustrated below in the case of what we might call canonical shapes of the cutting strip $\phi$. Note that the inner rim for skew shape, although canonical, has, to the best of our knowledge, only recently appeared in the literature in a somewhat hidden form in \cite{Oka1}.  Of course for standard shape it coincides with the hook.

\begin{equation}
\begin{array}{|c|c|c|c|c|}
\hline
\mbox{Row}&\mbox{Column}&\mbox{Hook}&\mbox{Inner rim}&\mbox{Outer rim}\cr
\hline
&&&&\cr
\vcenter{\hbox{
\begin{tikzpicture}[x={(0in,-0.15in)},y={(0.15in,0in)}] 
\foreach \j in {1,...,5} \draw (1,\j) rectangle +(-1,-1);
\foreach \j in {1,...,4} \draw (2,\j) rectangle +(-1,-1);
\foreach \j in {1,...,4} \draw (3,\j) rectangle +(-1,-1);
\foreach \j in {1,...,2} \draw (4,\j) rectangle +(-1,-1);
\fill[blue!40!white](0.3,-0.5)--(0.3,5.5)--(0.7,5.5)--(0.7,-0.5)--cycle;
\fill[red!40!white](1.3,-0.5)--(1.3,5.5)--(1.7,5.5)--(1.7,-0.5)--cycle;
\fill[cyan!40!white](2.3,-0.5)--(2.3,5.5)--(2.7,5.5)--(2.7,-0.5)--cycle;
\fill[magenta!40!white](3.3,-0.5)--(3.3,5.5)--(3.7,5.5)--(3.7,-0.5)--cycle;
\foreach \j in {1,...,3} \draw(0.5,\j-0.5)node{$\ast$};
\foreach \j in {1,...,2} \draw(1.5,\j-0.5)node{$\ast$};
\end{tikzpicture}
}}
&
\vcenter{\hbox{
\begin{tikzpicture}[x={(0in,-0.15in)},y={(0.15in,0in)}] 
\foreach \j in {1,...,5} \draw (1,\j) rectangle +(-1,-1);
\foreach \j in {1,...,4} \draw (2,\j) rectangle +(-1,-1);
\foreach \j in {1,...,4} \draw (3,\j) rectangle +(-1,-1);
\foreach \j in {1,...,2} \draw (4,\j) rectangle +(-1,-1);
\fill[blue!40!white](-0.5,0.3)--(4.5,0.3)--(4.5,0.7)--(-0.5,0.7)--cycle;
\fill[red!40!white](-0.5,1.3)--(4.5,1.3)--(4.5,1.7)--(-0.5,1.7)--cycle;
\fill[cyan!40!white](-0.5,2.3)--(4.5,2.3)--(4.5,2.7)--(-0.5,2.7)--cycle;
\fill[magenta!40!white](-0.5,3.3)--(4.5,3.3)--(4.5,3.7)--(-0.5,3.7)--cycle;
\fill[orange!40!white](-0.5,4.3)--(4.5,4.3)--(4.5,4.7)--(-0.5,4.7)--cycle;
\foreach \j in {1,...,3} \draw(0.5,\j-0.5)node{$\ast$};
\foreach \j in {1,...,2} \draw(1.5,\j-0.5)node{$\ast$};
\end{tikzpicture}
}}
&
\vcenter{\hbox{
\begin{tikzpicture}[x={(0in,-0.15in)},y={(0.15in,0in)}] 
\foreach \j in {1,...,5} \draw (1,\j) rectangle +(-1,-1);
\foreach \j in {1,...,4} \draw (2,\j) rectangle +(-1,-1);
\foreach \j in {1,...,4} \draw (3,\j) rectangle +(-1,-1);
\foreach \j in {1,...,2} \draw (4,\j) rectangle +(-1,-1);
\fill[blue!40!white](0.3,0.3)--(0.3,5.5)--(0.7,5.5)--(0.7,0.7)--(4.5,0.7)--(4.5,0.3)--cycle;
\fill[red!40!white](1.3,1.3)--(1.3,5.5)--(1.7,5.5)--(1.7,1.7)--(4.5,1.7)--(4.5,1.3)--cycle;
\fill[cyan!40!white](2.3,2.3)--(2.3,5.5)--(2.7,5.5)--(2.7,2.7)--(4.5,2.7)--(4.5,2.3)--cycle;
\foreach \j in {1,...,3} \draw(0.5,\j-0.5)node{$\ast$};
\foreach \j in {1,...,2} \draw(1.5,\j-0.5)node{$\ast$};
\end{tikzpicture}
}}
&
\vcenter{\hbox{
\begin{tikzpicture}[x={(0in,-0.15in)},y={(0.15in,0in)}] 
\foreach \j in {1,...,5} \draw (1,\j) rectangle +(-1,-1);
\foreach \j in {1,...,4} \draw (2,\j) rectangle +(-1,-1);
\foreach \j in {1,...,4} \draw (3,\j) rectangle +(-1,-1);
\foreach \j in {1,...,2} \draw (4,\j) rectangle +(-1,-1);
\fill[blue!40!white](4.5,0.3)--(2.3,0.3)--(2.3,2.3)--(1.3,2.3)--(1.3,3.3)--(0.3,3.3)--(0.3,5.5)--(0.7,5.5)--(0.7,3.7)--(1.7,3.7)--(1.7,2.7)--(2.7,2.7)--(2.7,0.7)--(4.5,0.7)--cycle;
\fill[red!40!white](5.5,1.3)--(3.3,1.3)--(3.3,3.3)--(2.3,3.3)--(2.3,4.3)--(1.3,4.3)--(1.3,6.5)--(1.7,6.5)--(1.7,4.7)--(2.7,4.7)--(2.7,3.7)--(3.7,3.7)--(3.7,1.7)--(5.5,1.7)--cycle;
\foreach \j in {1,...,3} \draw(0.5,\j-0.5)node{$\ast$};
\foreach \j in {1,...,2} \draw(1.5,\j-0.5)node{$\ast$};
\end{tikzpicture}
}}
&
\vcenter{\hbox{
\begin{tikzpicture}[x={(0in,-0.15in)},y={(0.15in,0in)}] 
\foreach \j in {1,...,5} \draw (1,\j) rectangle +(-1,-1);
\foreach \j in {1,...,4} \draw (2,\j) rectangle +(-1,-1);
\foreach \j in {1,...,4} \draw (3,\j) rectangle +(-1,-1);
\foreach \j in {1,...,2} \draw (4,\j) rectangle +(-1,-1);
\fill[blue!40!white](4.5,0.3)--(3.3,0.3)--(3.3,1.3)--(2.3,1.3)--(2.3,3.3)--(0.3,3.3)--(0.3,5.5)--(0.7,5.5)--(0.7,3.7)--(2.7,3.7)--(2.7,1.7)--(3.7,1.7)--(3.7,0.7)--(4.5,0.7)--cycle;
\fill[red!40!white](3.5,-0.7)--(2.3,-0.7)--(2.3,0.3)--(1.3,0.3)--(1.3,2.3)--(-0.7,2.3)--(-0.7,4.5)--(-0.3,4.5)--(-0.3,2.7)--(1.7,2.7)--(1.7,0.7)--(2.7,0.7)--(2.7,-0.3)--(3.5,-0.3)--cycle;
\foreach \j in {1,...,3} \draw(0.5,\j-0.5)node{$\ast$};
\foreach \j in {1,...,2} \draw(1.5,\j-0.5)node{$\ast$};
\end{tikzpicture}
}}
\cr
&&&&\cr
\hline
&&&&\cr
\vcenter{\hbox{
\begin{tikzpicture}[x={(0in,-0.15in)},y={(0.15in,0in)}] 
\foreach \j in {1,...,5} \draw (1,\j) rectangle +(-1,-1);
\foreach \j in {1,...,4} \draw (2,\j) rectangle +(-1,-1);
\foreach \j in {1,...,4} \draw (3,\j) rectangle +(-1,-1);
\foreach \j in {1,...,2} \draw (4,\j) rectangle +(-1,-1);
\draw[draw=blue,ultra thick] (0.5,2.8)--(0.5,5.2);
\draw[draw=red,ultra thick] (1.5,1.8)--(1.5,4.2);
\draw[draw=cyan,ultra thick] (2.5,-0.2)--(2.5,4.2);
\draw[draw=magenta,ultra thick] (3.5,-0.2)--(3.5,2.2);
\foreach \j in {1,...,3} \draw(0.5,\j-0.5)node{$\ast$};
\foreach \j in {1,...,2} \draw(1.5,\j-0.5)node{$\ast$};
\end{tikzpicture}
}}
&
\vcenter{\hbox{
\begin{tikzpicture}[x={(0in,-0.15in)},y={(0.15in,0in)}] 
\foreach \j in {1,...,5} \draw (1,\j) rectangle +(-1,-1);
\foreach \j in {1,...,4} \draw (2,\j) rectangle +(-1,-1);
\foreach \j in {1,...,4} \draw (3,\j) rectangle +(-1,-1);
\foreach \j in {1,...,2} \draw (4,\j) rectangle +(-1,-1);
\draw[draw=blue,ultra thick](4.2,0.5)--(1.8,0.5);
\draw[draw=red,ultra thick] (4.2,1.5)--(1.8,1.5);
\draw[draw=cyan,ultra thick] (3.2,2.5)--(0.8,2.5);
\draw[draw=magenta,ultra thick](3.2,3.5)--(-0.2,3.5);
\draw[draw=orange,ultra thick] (1.2,4.5)--(-0.2,4.5);
\foreach \j in {1,...,3} \draw(0.5,\j-0.5)node{$\ast$};
\foreach \j in {1,...,2} \draw(1.5,\j-0.5)node{$\ast$};
\end{tikzpicture}
}}
&
\vcenter{\hbox{
\begin{tikzpicture}[x={(0in,-0.15in)},y={(0.15in,0in)}] 
\foreach \j in {1,...,5} \draw (1,\j) rectangle +(-1,-1);
\foreach \j in {1,...,4} \draw (2,\j) rectangle +(-1,-1);
\foreach \j in {1,...,4} \draw (3,\j) rectangle +(-1,-1);
\foreach \j in {1,...,2} \draw (4,\j) rectangle +(-1,-1);
\draw[draw=blue,ultra thick](4.2,0.5)--(1.8,0.5);
\draw[draw=red,ultra thick] (4.2,1.5)--(1.8,1.5);
\draw[draw=cyan,ultra thick] (3.2,2.5)--(2.5,2.5)--(2.5,4.2);
\draw[draw=magenta,ultra thick](1.5,1.8)--(1.5,4.2);
\draw[draw=orange,ultra thick] (0.5,2.8)--(0.5,5.2);
\foreach \j in {1,...,3} \draw(0.5,\j-0.5)node{$\ast$};
\foreach \j in {1,...,2} \draw(1.5,\j-0.5)node{$\ast$};
\end{tikzpicture}
}}
&
\vcenter{\hbox{
\begin{tikzpicture}[x={(0in,-0.15in)},y={(0.15in,0in)}] 
\foreach \j in {1,...,5} \draw (1,\j) rectangle +(-1,-1);
\foreach \j in {1,...,4} \draw (2,\j) rectangle +(-1,-1);
\foreach \j in {1,...,4} \draw (3,\j) rectangle +(-1,-1);
\foreach \j in {1,...,2} \draw (4,\j) rectangle +(-1,-1);
\draw[draw=blue,ultra thick](4.2,0.5)--(2.5,0.5)--(2.5,2.5)--(1.5,2.5)--(1.5,3.5)--(0.5,3.5)--(0.5,5.2);
\draw[draw=red,ultra thick] (4.2,1.5)--(3.5,1.5)--(3.5,2.2);
\draw[draw=magenta,ultra thick] (3.2,3.5)--(2.5,3.5)--(2.5,4.2);
\foreach \j in {1,...,3} \draw(0.5,\j-0.5)node{$\ast$};
\foreach \j in {1,...,2} \draw(1.5,\j-0.5)node{$\ast$};
\end{tikzpicture}
}}
&
\vcenter{\hbox{
\begin{tikzpicture}[x={(0in,-0.15in)},y={(0.15in,0in)}] 
\foreach \j in {1,...,5} \draw (1,\j) rectangle +(-1,-1);
\foreach \j in {1,...,4} \draw (2,\j) rectangle +(-1,-1);
\foreach \j in {1,...,4} \draw (3,\j) rectangle +(-1,-1);
\foreach \j in {1,...,2} \draw (4,\j) rectangle +(-1,-1);
\draw[draw=blue,ultra thick](4.2,0.5)--(3.5,0.5)--(3.5,1.5)--(2.5,1.5)--(2.5,3.5)--(0.5,3.5)--(0.5,5.2);
\draw[draw=red,ultra thick] (2.5,-0.2)--(2.5,0.5)--(1.8,0.5);
\draw[draw=magenta,ultra thick] (1.5,1.8)--(1.5,2.5)--(0.8,2.5);
\foreach \j in {1,...,3} \draw(0.5,\j-0.5)node{$\ast$};
\foreach \j in {1,...,2} \draw(1.5,\j-0.5)node{$\ast$};
\end{tikzpicture}
}}
\cr
&&&&\cr
\hline
\end{array}
\end{equation}

Guided by these we have the following corollaries of Theorem~\ref{The-HGX}:

\begin{Corollary}[Jacobi-Trudi type identity \cite{Mac92}]\label{Cor-JTX}
Let $\lambda$ and $\mu$ be partitions such that $\mu\subseteq\lambda$. 
Let $\ell(\lambda)=l$ and let $\X=(x_{kc})$ with $k\in\{1,2,\ldots,n\}$ and $c\in\Z$. Then
\begin{equation}\label{eqn-JTX}
    s_{\lambda/\mu}(\X)
    =\left|\,  h_{\lambda_i-\mu_j -i +j}(\tau^{\mu_j-j+1}\X) \,\right|_{1\leq i,j\leq l} \,,
\end{equation}
where $h_{r}(\tau^m\X)=s_{(r)}(\tau^m\X)$ for all $r,m\in\Z$, with $h_{0}(\tau^m\X)=1$ and $h_{r}(\tau^m\X)=0$ if $r<0$.
\end{Corollary}

\noindent{\bf Proof}: Setting $\phi$ to be a single row of boxes extending rightwards from a box of 
content $-l+1$ to one of content $l'-1$, the corresponding outside decomposition of $F^{\lambda/\mu)}$ 
is of the form $\Theta=(\theta_1,\theta_2,\ldots,\theta_s)$ with $s=l$ and, if we order the strips from top to bottom,
$\theta_i=F^{(\lambda_i)/(\mu_i)}\simeq\phi_{\mu_i-i+1,\lambda_i-i}$ for $i=1,2,\ldots,l$.
Hence $(\theta_i\#\theta_j)\simeq\phi_{\mu_i-i+1,\lambda_j-j}$. This is a null strip if $\mu_i-i=\lambda_j-j$ and empty
if $\mu_i-i>\lambda_j-j$. In the other cases it is a row containing at least one box and the content of its leftmost box, 
$\mu_i-i+1$, gives the value of the shift parameter $m(i,j)$.
It follows that, in (\ref{eqn-HGX}) we have $s_{(\theta_i\#\theta_j)}(\X)=s_{(\lambda_j-\mu_i-j+i)}(\tau^{\mu_i-i+1}\X)$
with the null strip and empty strip cases yielding $1$ and $0$, respectively.
Taking the transpose of the resulting determinant gives (\ref{eqn-JTX}).
\qed

\begin{Corollary} [Dual Jacobi-Trudi type identity \cite{Mac92}]\label{Cor-DJTX}
Let $\lambda$ and $\mu$ be partitions such that $\mu\subseteq\lambda$ with conjugates $\lambda'$ and $\mu'$, respectively.
Let $\ell(\lambda')=l'$ and let $\X=(x_{kc})$ with $k\in\{1,2,\ldots,n\}$ and $c\in\Z$. Then
\begin{equation}\label{eqn-DJTX}
    s_{\lambda/\mu}(\X)
    =\left|\,  e_{\lambda'_i-\mu'_j -i +j}(\tau^{-\mu'_j+j-1}\X) \,\right|_{1\leq i,j\leq l'} \,,
\end{equation}
where $e_{r}(\tau^m\X)=s_{(1^r)}(\tau^m\X)$ for all $r,m\in\Z$, with $e_{0}(\tau^m\X)=1$ and $e_{r}(\tau^m\X)=0$ if $r<0$.
\end{Corollary}

\noindent{\bf Proof}: Setting $\phi$ to be a column of boxes extending upwards from a box of content $-l+1$ to one of content $l'-1$, 
the corresponding outside decomposition of $F^{\lambda/\mu)}$ is of the form $\Theta=(\theta_1,\theta_2,\ldots,\theta_s)$
with $s=l'$ and if the column strips are read from left to right, then $\theta_i=F^{(1^{\lambda'_i})/(1^{\mu'_i})}\simeq\phi_{-\lambda'_i+i,-\mu'_i+i-1}$ 
for $i=1,2,\ldots,l'$. 
It follows that $(\theta_i\#\theta_j)\simeq\phi_{-\lambda'_i+i,-\mu'_j+j-1}$. 
This a null strip if $\mu'_j-j=\lambda'_i-i$ and empty if $\mu'_j-j>\lambda'_i-i$. 
In the other cases it is a column containing at least one box and the content of its topmost box, $-\mu'_j+j-1$, is the 
value of its shift parameter $m(i,j)$.
Therefore, in (\ref{eqn-HGX}) we have
$s_{(\theta_i\#\theta_j)}(\X)=s_{(1^{\lambda'_i-\mu'_j-i+j})}(\tau^{-\mu'_j+j-1}\X)$ with the null strip and empty strip cases
yielding $1$ and $0$, respectively. This confirms the validity of (\ref{eqn-DJTX}), as required.
\qed

\begin{Corollary} [Giambelli type identity \cite{Mac92}]\label{Cor-GX}
Let $\lambda=(\alpha|\beta)$ and $\mu=(\gamma|\delta)$ be partitions of Frobenius rank $p$ and $q$, 
respectively, with $0\leq q\leq p$ and $\mu\subseteq\lambda$. Then
for all $\X=(x_{kc})$ with $k\in\{1,2,\ldots,n\}$ and $c\in\Z$
\begin{equation}\label{eqn-GX}
    s_{\lambda/\mu}(\X)
		= (-1)^q \left|\begin{array}{cc} s_{(\alpha_i|\beta_j)}(\X)& h_{\alpha_i-\gamma_k}(\tau^{\gamma_k+1}\X))\cr\cr
														  e_{\beta_j-\delta_\ell}(\tau^{-\delta_\ell-1}\X)&0
						\end{array}								
		\right|_{\begin{array}{c}\sc 1\leq i,j\leq p\cr \sc 1\leq k,l\leq q\end{array}}	\,,
	\end{equation} 
with $h_0(\tau^m\X)=e_0(\tau^m\X)=1$ and $h_r(\tau^m\X)=e_r(\tau^m\X)=0$ for all $r<0$ and any $m$.
\end{Corollary}

\noindent{\bf Proof}: This time we set $\phi$ to be a hook of boxes extending upwards from a box of 
content $-l+1$ to one of content $0$ and then rightwards to one of content $l'-1$. 
The corresponding outside decomposition of $F^{\lambda/\mu)}$ is of the form $\Theta=(\theta_1,\theta_2,\ldots,\theta_s)$
with $s=p+q$. The strips $\theta_i$ may be ordered in such a way that:
\begin{equation}
  \theta_i=\begin{cases}
	               F^{(0|\beta_i)/(0|\delta_i)}\simeq\phi_{-\beta_i,-\delta_i-1}&\mbox{for $i=1,2,\ldots,q$}\,;\cr
								 F^{(\alpha_i|\beta_i)}\simeq\phi_{-\beta_i,\alpha_i}&\mbox{for $i=q+1,q+2,\ldots,p$}\,;\cr
								 F^{(\alpha_j|0)/(\gamma_j|0)}\simeq\phi_{\gamma_j+1,\alpha_j}&\left\{\begin{array}{l}\mbox{for $i\=p\+1,p\+2,\ldots,p\+q$}\cr
								                                                                                \mbox{with $j\=p\+q\+1\-i =1, 2, \ldots, q$}\,.\cr
																																																\end{array}\right.\cr
	          \end{cases}
\end{equation}
These expressions enable one to write down $(\theta_i\#\theta_j)$ in the form of a $(p+q)\times(p+q)$ matrix using (\ref{eqn-hash}). 
After applying a sequence of $q$ transpositions $(i,p+q+1-i)$ for $i=1,2,\ldots,q$ to the columns of this matrix
one obtains the matrix
\begin{equation}
\left[
\begin{array}{cc}
   \phi_{-\beta_i,\alpha_j}&\phi_{-\beta_i,-\delta_\ell-1}\cr
	 \phi_{\gamma_k+1,\alpha_j}&\emptyset\cr
\end{array}
\right]
\end{equation}
with $i,j=1,2,\ldots,p$ and $k,l=1,2,\ldots,q$. 
Using this in (\ref{eqn-HGX}), taking account of the shift factors $\tau^{\gamma_k+1}$, $\tau^{-\delta_\ell-1}$ and $\tau^0$,
appropriate to the row, column and hook shapes, respectively, and not forgetting the factor of $(-1)^q$ arising from the above 
permutation of columns, yields (\ref{eqn-GX}).
\qed

We take this opportunity to point out that another Giambelli-related identity, the so-called $M$-shifted skew Schur function Giambelli identities of Matsuno and Moriyama~\cite{MM}, is also a special case of Theorem \ref{The-HG}.  In this case the cutting strip is also a hook, but instead of the hook shaped cutting strip appropriate to the Giambelli identity of Corollary~\ref{Cor-GX} we use a ``shifted'' hook of the same shape with corner at position $(M,0)$ rather than $(0,0)$.  This result has also been given a ninth variation extension by Furukawa and Moriyama~\cite{FM}, and it too is derivable in the same way using an outside decomposition.

By way of preparation for the next two Corollaries it is convenient, 
for any $\kappa_i$ and $\nu=(\nu_1,\nu_2,\ldots,\nu_n)$, 
to let $(\kappa_i,\nu)=(\kappa_i,\nu_1,\nu_2,\ldots,\nu_n)$, $(\nu,\kappa_i)=(\nu_1,\nu_2,\ldots,\nu_n,\kappa_i)$, 
$(\kappa_i,\nu\backslash\nu_j)=(\kappa_i,\nu_1,\ldots,\hat{\nu}_j,\ldots,\nu_n)$ and 
$(\nu\backslash\nu_j,\kappa_i)=(\nu_1,\ldots,\hat{\nu}_j,\ldots,\nu_n,\kappa_i,)$
where $\hat{\nu}_j$ signifies the omission of the part $\nu_j$. Sequences of this type will
appear as Frobenius type labels for either of the two parts of skew Schur functions.
Any non-standard labels of this type are to be regularised using the 
identities $s_{(\eta|\zeta)/(\gamma|\delta)}(\X)=\sgn(\sigma)\sgn(\tau)s_{(\sigma(\eta)|\tau(\zeta))/(\gamma|\delta)}(\X)$ 
or $s_{(\alpha|\beta)/(\eta|\zeta)}(\X)=\sgn(\sigma)\sgn(\tau)s_{(\alpha|\beta)/(\sigma(\eta)|\tau(\zeta))}(\X)$ 
for any permutations $\sigma$ and $\tau$ of the parts of $\eta$ and $\zeta$. 
This implies that if either $\eta$ or $\zeta$ has any two equal parts then both $s_{(\eta|\zeta)/(\gamma|\delta)}(\X)=0$
and $s_{(\alpha|\beta)/(\eta|\zeta)}(\X)=0$. 
In all other cases one may choose $\sigma$ and $\tau$
so that the parts of $\sigma(\eta)$ and $\tau(\zeta)$ are strictly decreasing, that is to say, standard labels.

\begin{Corollary} [Okada type inner rim identity]\label{Cor-OkaX}
Let $\lambda=(\alpha|\beta)$ and $\mu=(\gamma|\delta)$ be partitions of Frobenius rank $p$ and $q$, respectively, 
with $0\leq q\leq p$ and $\mu\subseteq\lambda$. Then
for all $\X=(x_{kc})$ with $k\in\{1,2,\ldots,n\}$ and $c\in\Z$
\begin{equation}\label{eqn-OkaX}
     s_{\lambda/\mu}(\X)                
							  		=(-1)^q \left|  \begin{array}{cc} s_{((\alpha_i,\gamma)|(\beta_j,\delta))/(\gamma|\delta)}(\X)
		                         & s_{((\alpha_i,\gamma\backslash\gamma_k)|\delta)/(\gamma|\delta)}(\X)\cr\cr
														  s_{(\gamma|(\beta_j,\delta\backslash\delta_l))/(\gamma|\delta)}(\X)
		                         & 0
						                 \end{array}								
		                  \right|_{\begin{array}{c}\sc 1\leq i,j\leq p\cr \sc 1\leq k,l\leq q\end{array}}			\,.								
\end{equation}
\end{Corollary}

\noindent{\bf Proof}: 
The appropriate choice of $\phi$ is now the inner rim of the skew diagram $F^{\lambda/\mu}$.
In Frobenius notation this means that $\phi=F^{(\alpha_1,\gamma|\beta_1,\delta)/(\gamma|\delta)}$. 
As a consequence, each term of the form $(\theta_p\#\theta_q)\simeq\phi_{ad}$
in (\ref{eqn-HGX}) of Theorem~\ref{The-HGX}
can be expressed as a skew shape of the form $F^{(\eta|\zeta)/(\gamma|\delta)}$
with shift parameter $m(p,q)=0$.
To be more precise, one finds just three types of expression, namely: 
\begin{equation}\label{eqn-three-phi}
\begin{array}{rcl}
    \phi_{-\beta_j,\alpha_i}&=&F^{(\sigma(\alpha_i,\gamma))|\tau(\beta_j,\delta))/(\gamma|\delta)}\,;\cr 
    \phi_{-\beta_j,-\delta_\ell-1}&=&F^{(\gamma|\tau(\beta_j,\delta\backslash\delta_\ell))/(\gamma|\delta)}\,;\cr
    \phi_{\gamma_k+1,\alpha_i}&=&F^{(\sigma(\alpha_i,\gamma\backslash\gamma_k))|\delta)/(\gamma|\delta)}\,,\cr
\end{array}
\end{equation}
where in each case $\sigma$ and $\tau$ are those permutations whose action on arm and leg lengths leads to 
$F^{(\eta|\zeta)/(\gamma|\delta)}$ with $\eta$ and $\zeta$ either strict partitions or partitions in which at least one
part is repeated. In the latter case the corresponding generalised skew Schur function $s_{(\eta|\zeta)/(\gamma|\delta)}(\X)=0$.

The transformations in the first case can best be viewed as the diagonal, content preserving translation 
of the boxes of the hook $(\alpha_i|\beta_j)$ through $F^\mu$ until they form a sub-strip of $\phi$, the inner 
rim of $F^\lambda$, spanning its main diagonal. 
In the second and third cases the boxes of $(0|\beta_j)$ and $(\alpha_i|0)$ are diagonally translated in the same way, filling  
the $\delta_\ell$ and $\gamma_k$ boxes that are deleted from $F^\mu$ with the remaining boxes forming substrips of $\phi$ 
lying below and above, respectively, the main diagonal.
If $(\eta|\zeta)=(\gamma|\delta)$, which will be the case if and only if either $\gamma_k=\alpha_i$ or $\delta_l=\beta_j$
for some pair $(i,k)$ or $(j,\ell)$, then $s_{(\eta|\zeta)/(\gamma|\delta)}(\X)=\pm1$ with the sign given by $\sgn(\sigma)$ 
or $\sgn(\tau)$, respectively.
Using Theorem~\ref{The-HGX},
ordering rows and permuting columns of the resulting determinant in a suitable way, followed by taking the transpose, 
yields (\ref{eqn-OkaX}) up to a sign factor determined by the column permutation and the various factors $\sgn(\sigma)$ 
and $\sgn(\tau)$ that are common to particular rows or particular columns. 
That the sign factor in (\ref{eqn-OkaX}) is just $(-1)^q$ follows from the following 
more detailed considerations.

In order to compute all the sign factors for given $\lambda=(\alpha|\beta)$ and $\mu=(\gamma|\delta)$ we 
let $(\alpha\merge\gamma)$ and $(\beta\merge\delta)$ be the the two weakly decreasing sequences obtained by 
interleaving the parts of $\alpha$ and $\gamma$ and those of $\beta$ and $\delta$ while keeping track of the origin of each part. For example,
for $\lambda=(9,8,4,3,2,0\, | \,9,7,6,5,2,1)$ and $\mu=(5,4,2\, |\, 7,4,3)$ we have $(\alpha\merge\gamma)=(9,8,5,4,4,3,2,2,0)$
and $(\beta\merge\delta)=(9,7,7,6,5,4,3,2,1)$ in which the consecutive parts are to be identified with 
$(\alpha_1,\alpha_2,\gamma_1,\alpha_3,\gamma_2,\alpha_4,\alpha_5,\gamma_3,\alpha_6)$
and $(\beta_1,\beta_2,\delta_1,\beta_3,\beta_4,\delta_2,\delta_3,\beta_5,\beta_6)$.
It is to be noted that when we have equal parts $\alpha_3=\gamma_2=4$, $\alpha_5=\gamma_3=2$ and $\beta_2=\delta_1=7$
the interleaving operation $\merge$ leaves $\alpha_3$ to the left of $\gamma_2$, $\alpha_5$ to the left of $\gamma_3$, 
and $\beta_2$ to the left of $\delta_1$.

From the two sequences $(\alpha\merge\gamma)$ and $(\beta\merge\delta)$
we can read off the sign factors arising on writing
$s_{((\alpha_i,\gamma)|(\beta_j,\delta))/(\gamma|\delta)}(\X)$ in the form
$\sgn(\sigma)\sgn(\tau)s_{(\sigma(\alpha_i\merge\gamma)|\tau(\beta_j\merge\delta))/(\gamma|\delta)}(\X)$.
One finds 
$\sgn(\sigma)=(-1)^{\epsilon(\alpha_i\merge\gamma)}$
with $\epsilon(\alpha_i\merge\gamma)=\#\{\gamma_k >\alpha_i |k=1,2,\ldots,q\}$
and
$\sgn(\tau)=(-1)^{\epsilon(\beta_j\merge\delta)}$ with $\epsilon(\beta_j\merge\delta)=\#\{\delta_\ell>\beta_j | \ell=1,2,\ldots,q\}$.
In writing $s_{((\alpha_i,\gamma\backslash\gamma_k)|\delta)/(\gamma|\delta)}(\X)$ 
in the form $\sgn(\sigma)s_{((\alpha_i\merge\gamma\backslash\gamma_k)|\delta)/(\gamma|\delta)}(\X)$ and, similarly,
$s_{(\gamma|(\beta_j,\delta\backslash\delta_l))/(\gamma|\delta)}(\X)$ in the form
$\sgn(\tau)s_{(\gamma|(\beta_j\merge\delta\backslash\delta_l))/(\gamma|\delta)}(\X)$, it is important to recall that these 
skew Schur functions are identically zero unless $\alpha_i\geq\gamma_k$ and $\beta_j\geq\delta_\ell$, 
respectively. This ensures that after the above regularisation of every term in (\ref{eqn-OkaX}), including setting to $0$ those that are zero,
one can extract common factors $\sgn(\sigma)$ and $\sgn(\tau)$ from each row and column of the determinant leaving an overall
sign factor of
\begin{equation}\label{eqn-signs}
      (-1)^q\ (-1)^{\epsilon(\alpha\merge\gamma)}\ (-1)^{\epsilon(\beta\merge\delta)},
\end{equation}
where $\epsilon(\alpha\merge\gamma)=\sum_{i=1}^p \#\{\gamma_k>\alpha_i|k=1,2,\ldots,q\}$
and $\epsilon(\beta\merge\delta)=\sum_{j=1}^p \#\{\delta_\ell>\beta_j|\ell=1,2,\ldots,q\}$.

In our example $\lambda=(9,8,4,3,2,0\, |\, 9,7,6,5,2,1)$ and $\mu=(5,4,2\, |\, 7,4,3)$,
for which we have $q=3$ and $(\alpha\merge\gamma)=(9,8,5,4,4,3,2,2,0)$
and $(\beta\merge\delta)=(9,7,7,6,5,4,3,2,1)$, we find $\epsilon(\alpha\merge\gamma)=8$
and $\epsilon(\beta\merge\delta)=8$ leading to an overall sign factor of $-1$.
 
The expression (\ref{eqn-signs}) has to be compared with the sign factor associated with
re-ordering the rows and columns of the determinant in (\ref{eqn-HGX}) to match those of the transpose of
the determinant in (\ref{eqn-OkaX}). To this end we must first determine the specification of row and column 
labels $a$ and $b$ as determined by the inner rim outside decomposition strips $\theta_r\simeq\phi_{ab}$ 
for our given $\lambda=(\alpha|\beta)$ and $\mu=(\gamma|\delta)$.

The sequence $(\ov{\beta}\merge\ov{\delta})$, which is obtained from $(\beta\merge\delta)$ by changing the sign of every term,
may be followed by the reverse of the sequence $(\alpha\merge\gamma)$ to give a weakly increasing sequence.
This may be augmented by left brackets, ``$($'', to the left of each $\ov{\beta_j}$ and $\gamma_k$, and equinumerous 
right brackets, ``$)$'', to the right of each $\ov{\delta_\ell}$ and $\alpha_i$. In our example this yields
\begin{equation}\label{eqn-ex-merge}
  (\ov{\beta_1},(\ov{\beta_2},\ov{\delta_1})(\ov{\beta_3},(\ov{\beta_4},\ov{\delta_2}),\ov{\delta_3}),(\ov{\beta_5},(\ov{\beta_6},
	\alpha_6),(\gamma_3,\alpha_5),\alpha_4),(\gamma_2,\alpha_3),(\gamma_1,\alpha_2),\alpha_1)\,.
\end{equation}

In the case $\lambda=(9,8,4,3,2,0\, |\, 9,7,6,5,2,1)$ and $\mu=(5,4,2\, |\, 7,4,3)$ the above bracketing can be seen 
reflected below in the nesting of the various strips $\theta_r$ constituting the corresponding inner rim outside
decompostion of $F^{\lambda/\mu}$.
\begin{equation}\label{eqn-inner-theta}
\vcenter{\hbox{
\begin{tikzpicture}[x={(0in,-0.18in)},y={(0.18in,0in)}] 
\foreach \j in {1,...,10} \draw[thick] (1,\j) rectangle +(-1,-1);
\foreach \j in {1,...,10} \draw[thick] (2,\j) rectangle +(-1,-1);
\foreach \j in {1,...,7} \draw[thick] (3,\j) rectangle +(-1,-1);
\foreach \j in {1,...,7} \draw[thick] (4,\j) rectangle +(-1,-1);
\foreach \j in {1,...,7} \draw[thick] (5,\j) rectangle +(-1,-1);
\foreach \j in {1,...,6} \draw[thick] (6,\j) rectangle +(-1,-1);
\foreach \j in {1,...,6} \draw[thick] (7,\j) rectangle +(-1,-1);
\foreach \j in {1,...,4} \draw[thick] (8,\j) rectangle +(-1,-1);
\foreach \j in {1,...,4} \draw[thick] (9,\j) rectangle +(-1,-1);
\foreach \j in {1,...,1} \draw[thick] (10,\j) rectangle +(-1,-1);
\draw[draw=blue,ultra thick](10.2,0.5)--(8.5,0.5)--(8.5,1.5)--(6.5,1.5)--(6.5,3.5)--(3.5,3.5)--(3.5,5.5)--(2.5,5.5)--(2.5,6.5)--(0.5,6.5)--(0.5,10.2);
\draw[draw=red,ultra thick] (7.2,4.5)--(4.5,4.5)--(4.5,6.5)--(3.5,6.5)--(3.5,7.2);
\draw[draw=cyan,ultra thick] (7.2,5.5)--(5.5,5.5)--(5.5,6.2);
\draw[draw=red,ultra thick] (2.2,7.5)--(1.5,7.5)--(1.5,10.2);
\draw[draw=black,ultra thick](2.5,6.8)--(2.5,7.2);
\draw[draw=black,ultra thick](4.5,6.8)--(4.5,7.2);
\draw[draw=red,ultra thick] (9.2,2.5)--(7.5,2.5)--(7.5,4.2);
\draw[draw=cyan,ultra thick] (9.2,3.5)--(8.5,3.5)--(8.5,4.2);
\draw[draw=black,ultra thick](9.2,1.5)--(8.8,1.5);

\foreach \j in {1,...,6} \draw(0.5,\j-0.5)node{$\ast$};
\foreach \j in {1,...,6} \draw(1.5,\j-0.5)node{$\ast$};
\foreach \j in {1,...,5} \draw(2.5,\j-0.5)node{$\ast$};
\foreach \j in {1,...,3} \draw(3.5,\j-0.5)node{$\ast$};
\foreach \j in {1,...,3} \draw(4.5,\j-0.5)node{$\ast$};
\foreach \j in {1,...,3} \draw(5.5,\j-0.5)node{$\ast$};
\foreach \j in {1,...,1} \draw(6.5,\j-0.5)node{$\ast$};
\foreach \j in {1,...,1} \draw(7.5,\j-0.5)node{$\ast$};
\draw (0.7,10.7)node{$\blue\theta_1$};
\draw (1.7,10.7)node{$\red\theta_6$};
\draw (2.7,7.7)node{$\theta_5$};
\draw (3.7,7.7)node{$\red\theta_2$};
\draw (4.7,7.7)node{$\theta_4$};
\draw (5.7,6.7)node{$\cyan\theta_3$};
\draw (9.7,1.7)node{$\theta_9$};
\draw (8.7,4.7)node{$\cyan\theta_8$};
\draw (7.7,4.7)node{$\red\theta_7$};
\end{tikzpicture}
}}
\qquad
\begin{array}{rcll}
(\ov{\beta_1},\alpha_1)&&\theta_1\simeq\phi_{-9,9}\cr
(\ov{\beta_5},\alpha_4)&&\theta_2\simeq\phi_{-2,3}\cr
(\ov{\beta_6},\alpha_6)&&\theta_3\simeq\phi_{-1,0}\cr
(\gamma_3,\alpha_5)&&\theta_4\simeq\phi_{3,2}&\mbox{Null strip}\cr
(\gamma_2,\alpha_3)&&\theta_5\simeq\phi_{5,4}&\mbox{Null strip}\cr
(\gamma_1,\alpha_2)&&\theta_6\simeq\phi_{6,8}\cr
(\ov{\beta_3},\ov{\delta_3})&&\theta_7\simeq\phi_{-6,-4}\cr
(\ov{\beta_4},\ov{\delta_2})&&\theta_8\simeq\phi_{-5,-5}\cr
(\ov{\beta_2},\ov{\delta_1})&&\theta_9\simeq\phi_{-7,-8}&\mbox{Null strip}\cr
\end{array}
\end{equation}

By means of the sequential extraction of bracketed pairs $(~\cdots~,~\cdots~)$, one can identify pairs 
$(\ov{\beta_{j_m}},\ov{\delta_m})$ and $(\gamma_m,\alpha_{i_m})$ for $m=1,2,\ldots,q$ and 
$(\ov{\beta_{j_m}},\alpha_{i_m})$ for $m=q+1,q+2,\ldots,p$. 
As exemplified above on the right,
these define in each case, as in (\ref{eqn-three-phi}), the start and end points 
of a constituent strip $\theta_r$ of the inner rim outside decomposition of $F^{(\alpha|\beta)/(\gamma/\delta)}$. 
It is convenient to adopt the canonical ordering of these strips $\theta_r$ specified in  
Proposition~3.13 of~\cite{Ham-thesis}, whereby 
\begin{equation}
   \theta_r\simeq\begin{cases}
	         \phi_{-\beta_{j_{q+r}},\alpha_{i_{q+r}}}&\mbox{for $r=1,2,\ldots,p-q$;}\cr
					 \phi_{\gamma_{p+1-r}+1,\alpha_{i_{p+1-r}}}&\mbox{for $r=p-q+1,\ldots,p-1,p$;}\cr
					 \phi_{-\beta_{j_{p+q+1-r}},-\delta_{p+q+1-r}-1}&\mbox{for $r=p+1,p+2,\ldots,p+q$.}\cr
					\end{cases}
\end{equation}
These three cases correspond to strips crossing, lying above, and lying below the main diagonal.
Those crossing the main diagonal at $(d,d)$ are placed in order of increasing $d$, those lying above 
in order of increasing content, $\gamma_k+1$, of their initial box, and those lying below in order of decreasing 
content, $-\delta_\ell-1$, of their final box.

In our example the order of the corresponding bracketed pairs is
\begin{equation}\label{eqn-pairs}
  (\ov{\beta_1},\alpha_1),(\ov{\beta_5},\alpha_4),(\ov{\beta_6},\alpha_6),
	(\gamma_3,\alpha_5),(\gamma_2,\alpha_3),(\gamma_1,\alpha_2),
	(\ov{\beta_3},\ov{\delta_3}),(\ov{\beta_4},\ov{\delta_2}),(\ov{\beta_2},\ov{\delta_1})	
\end{equation}
as indicated on the right in (\ref{eqn-inner-theta}).

This choice of canonical order results in the rows and columns of the inner rim outside decomposition 
determinant in (\ref{eqn-HGX}) being specified by the sequences of the left and right hand elements of each pair.
In our example these take the form:
\begin{equation}\label{eqn-alpha-beta}
   (\ov{\beta_1},\ov{\beta_5},\ov{\beta_6},\gamma_3,\gamma_2,\gamma_1,\ov{\beta_3},\ov{\beta_4},\ov{\beta_2})
	\quad\mbox{and}\quad
	 (\alpha_1,\alpha_4,\alpha_6,\alpha_5,\alpha_3,\alpha_2,\ov{\delta_3},\ov{\delta_2},\ov{\delta_1})	\,.
\end{equation}
In order to conform with the row and column labelling of the determinant in (\ref{eqn-OkaX}) one first
permutes all $\gamma_k$ labels to the right of all $\ov{\beta_j}$ labels without changing their order, thereby giving rise 
in general to a sign factor of $(-1)^q$ . Since the $\gamma$'s and $\delta$'s are both in reverse order compared with (\ref{eqn-OkaX}), their reordering gives rise to two sign factors $(-1)^{q(q-1)/2}$ which together give just $+1$ for all $q$. To permute the $\alpha_i$ labels so that they are 
arranged from left to right in the order $(\alpha_1,\alpha_2,\ldots,\alpha_p)$ it is 
necessary to permute each $\alpha_i$ through a succession of transpositions with each $\alpha_m<\alpha_i$ to its right.
However, quite generally, thanks to the bracketing of pairs prior to their being canonically re-ordered, 
each $\alpha_i$ has as many $\alpha_m$ with $m<i$ to its right as there are $\gamma_k>\alpha_i$, as can be seen, 
for example, on comparing (\ref{eqn-pairs}) and (\ref{eqn-alpha-beta}). The required number of transpositions
is therefore $\#\{\gamma_k>\alpha_i|k=1,2,\ldots,q\}$. 
In exactly the same way, to permute the $\ov{\beta_j}$ labels so that they are 
arranged from left to right in the order $(\ov{\beta_1},\ov{\beta_2},\ldots,\ov{\beta_p})$, it is 
necessary to permute each $\ov{\beta_j}$ through a succession of transpositions with each 
$\ov{\beta_m}$ to its right with $\beta_m<\beta_j$. The bracketing rules then imply 
that the required number of transpositions is $\#\{\delta_\ell>\beta_j|k=1,2,\ldots,q\}$.
Permuting rows and columns of the inner rim outside decomposition determinant in just the same way as we have permuted their labels 
then leads to an overall sign of
\begin{equation}
         (-1)^q\ (-1)^{\epsilon(\alpha\merge\gamma)}\ (-1)^{\epsilon(\beta\merge\delta)}
\end{equation}
precisely as in (\ref{eqn-signs}).

This completes the proof of Corollary~\ref{Cor-OkaX} in the explicit form first 
established by Okada~\cite{Oka1} for skew Schur functions of argument $\x$, but here seen as a consequence of our choice 
of an inner rim cutting strip $\phi$ and applying more generally for all $\X=(x_{kc})$.
\qed

\begin{Corollary} [Lascoux-Pragacz type outer rim identity \cite{LP88}] \label{Cor-LPX}
Let $\lambda=(\alpha|\beta)$ and $\mu=(\gamma|\delta)$ be partitions of Frobenius rank $p$ and $q$, respectively, 
with $0\leq q\leq p$ and $\mu\subseteq\lambda$. Then
for all $\X=(x_{kc})$ with $k\in\{1,2,\ldots,n\}$ and $c\in\Z$
\begin{equation}\label{eqn-outer-sfnX}
     s_{\lambda/\mu}(\X)                
							  		=(-1)^q \left|\begin{array}{cc} s_{(\alpha|\beta)/((\alpha\backslash\alpha_i)|(\beta\backslash\beta_j))}(\X)
		                         & s_{(\alpha|\beta)/(((\alpha\backslash\alpha_i),\gamma_k)|\beta)}(\X)\cr\cr
														  s_{(\alpha|\beta)/(\alpha|((\beta\backslash\beta_j),\delta_\ell))}(\X)
		                         & 0
						                 \end{array}								
		                  \right|_{\begin{array}{c}\sc 1\leq i,j\leq p\cr \sc 1\leq k,l\leq q\end{array}}			\,.								
\end{equation} 
\end{Corollary}

\noindent{\bf Proof}: In this case the appropriate choice of $\phi$ is the outer rim of $F^{\lambda/\mu}$, 
that is to say in Frobenius notation $\phi=F^{(\alpha|\beta)/(\alpha\backslash\alpha_1|\beta\backslash\beta_1)}$. 
The outside decomposition corresponding to this outer rim $\phi$
is a simple diagonal translation of that for the inner rim $\phi$ previously encountered. Ordering the strips in the same 
canonical way as before, the content 
of the start and end boxes of each strip $\theta_r$ are identical for the inner and outer rim cutting strips, 
but of course their shapes are different. The same is true of their combinations under the $\#$ operation. 
We find this time that $(\theta_p\#\theta_q)\simeq\phi_{ad}$ takes one or other of the 
following three forms: 
\begin{equation}
\begin{array}{rcl}
    \phi_{-\beta_i,\alpha_j}&=&F^{(\alpha|\beta)/(\alpha\backslash\alpha_j)|(\beta\backslash\beta_i)}\,;\cr 
    \phi_{-\beta_j,-\delta_\ell-1}&=&F^{(\alpha|\beta)/(\alpha|\tau(\beta\backslash\beta_j,\delta_\ell))}\,,\cr
    \phi_{\gamma_k+1,\alpha_j}&=&F^{(\alpha|\beta)/(\sigma(\alpha\backslash\alpha_j,\gamma_k)|\beta)}\, ,\cr
\end{array}
\end{equation}
where the permutations $\sigma$ and $\tau$ are those that in Frobenius notation result in either strict partitions or those 
with neighbouring equal parts, with the latter giving rise to corresponding generalised skew Schur functions that are zero. 
Since each $(\theta_p\#\theta_q)\simeq\phi_{ad}$ is expressed as above in the form $F^{(\alpha|\beta)/(\eta|\zeta)}$
it follows that all shift parameters $m(p,q)$ are zero.
Then using Theorem~\ref{The-HGX} and permuting rows and columns 
of the resulting determinant in a suitable way, followed by taking the transpose,
yields (\ref{eqn-outer-sfnX}) up to a sign factor.

The sign factors arising from writing $s_{(\alpha|\beta)/(((\alpha\backslash\alpha_i),\gamma_k)|\beta)}(\X)$
and $s_{(\alpha|\beta)/(\alpha|((\beta\backslash\beta_j),\delta_\ell))}(\X)$ in the forms
$\sgn(\sigma) s_{(\alpha|\beta)/(\sigma(\alpha\backslash\alpha_j,\gamma_k)|\beta)}(\X)$ and 
$\sgn(\tau) s_{(\alpha|\tau(\beta\backslash\delta_\ell))}(\X)$,
respectively, are obtained, as in the proof of 
Corollary~\ref{Cor-OkaX}, from the merged sequences $(\alpha\merge\gamma)$ and $(\beta\merge\delta)$. 
The difference is that one is now moving $\gamma_k$ and $\delta_\ell$ to the left through the components of 
$\alpha$ and $\beta$, respectively. As a result $\sgn(\sigma)=(-1)^{\epsilon(\alpha\merge\gamma_k)}$
with $\epsilon(\alpha\merge\gamma_k)=\#\{\alpha_i>\gamma_k|i=1,2,\ldots,p\}$
and $\sgn(\tau)=(-1)^{\epsilon(\beta\merge\delta_\ell)}$ with $\epsilon(\beta\merge\delta_\ell)=\#\{\beta_j>\delta_\ell|i=1,2,\ldots,p\}$.
Extracting these sign factors common to columns and rows indexed by $k$ and $\ell$, respectively, in the determinant of 
(\ref{eqn-outer-sfnX}) and combining them all with the given factor of $(-1)^q$ gives the overall sign factor 
appearing in (\ref{eqn-signs}). 
On the other hand, in fixing the numbering of the strips $\theta_r$ in this
outer rim outside decomposition to be canonical, as in the case of the inner rim outside decomposition, the overall sign factor arising
from permuting the rows and columns of the corresponding determinant in Theorem~\ref{The-HGX} so as 
to coincide with the ordering of columns and rows in (\ref{Cor-LPX}) is the same
as that arising in Corollary~\ref{Cor-OkaX}. This is again that of (\ref{eqn-signs}).
This confirms the validity of the sign factor $(-1)^q$ appearing in (\ref{eqn-outer-sfnX}), 
thereby completing the proof of Corollary~\ref{Cor-LPX} which represents
an explicit form of a result established first by Lascoux and Pragacz~\cite{LP88} for skew Schur functions of argument $\x$, 
but here applying more generally for all arguments $\X=(x_{kc})$.
\qed


Returning to our example with $\lambda=(4,2,1|3,2,0)$ and $\mu=(2,0|1,0)$, the connection between 
the inner and outer rim outside decompositions can be seen in the following illustration in which
the fact that $\gamma_1=\alpha_2=2$ and $\delta_2=\beta_3=0$ is the origin of the null strips 
consisting of a single edge. Clearly, the sets of outside decomposition strips are just the 
diagonal translation of one another.
\begin{equation}\label{eqn-inner-outer-phi-theta}
\begin{array}{cc}
\mbox{Inner rim outside decomposition}&\mbox{Outer rim outside decomposition}\cr\cr
\end{array}
\end{equation}
In both cases the canonical ordering of the strips $\theta_r$ leads to the identification:
\begin{equation}\label{eqn-phi-ab}
\theta_1\simeq\phi_{\ov34},~~\theta_2\simeq\phi_{11},~~\theta_3\simeq\phi_{32},~~\theta_4\simeq\phi_{0\ov1},~~\theta_5\simeq\phi_{\ov2\ov2}. 
\end{equation}
It follows that for $\lambda=(\blue{4,2,1}|\blue{3,2,0})$ and $\mu=(\red{2,0}|\red{1,0})$ we have in both inner and outer rim cases
\begin{equation}\label{eqn-inner-phi-ad}
(\theta_p\#\theta_q)=\left[\, \begin{array}{ccccc} 
                        \phi_{\ov34}&\phi_{\ov31}&\phi_{\ov32}&\phi_{\ov3\ov1}&\phi_{\ov3\ov2}\cr
												\phi_{14}&\phi_{11}&\phi_{12}&\phi_{1\ov1}&\phi_{1\ov2}\cr
												\phi_{34}&\phi_{31}&\phi_{32}&\phi_{3\ov1}&\phi_{3\ov2}\cr
												\phi_{04}&\phi_{01}&\phi_{02}&\phi_{0\ov1}&\phi_{0\ov2}\cr
												\phi_{\ov24}&\phi_{\ov21}&\phi_{\ov22}&\phi_{\ov2\ov1}&\phi_{\ov2\ov2}\cr
																								\end{array}
                     \,\right]
\end{equation}
However, the cutting strips are different in the two cases. Viewed diagrammatically, in the inner rim case we find
\begin{equation}\label{eqn-inner-diagrams}
\left[\, 
\begin{array}{ccccc} 
\vcenter{\hbox{
\begin{tikzpicture}[x={(0in,-0.15in)},y={(0.15in,0in)}] 
\foreach \j in {1,...,5} \draw (1,\j) rectangle +(-1,-1);
\foreach \j in {1,...,4} \draw (2,\j) rectangle +(-1,-1);
\foreach \j in {1,...,3} \draw (3,\j) rectangle +(-1,-1);
\foreach \j in {1} \draw (4,\j) rectangle +(-1,-1);
\foreach \j in {1,2,3} \draw(1-0.5,\j-0.5)node{$\ast$};
\foreach \j in {1,2} \draw(2-0.5,\j-0.5)node{$\ast$};
\draw (1-0.5,4-0.5)node{$\sc3$};\draw (1-0.5,5-0.5)node{$\sc4$};
\draw (2-0.5,3-0.5)node{$\sc1$};\draw (2-0.5,4-0.5)node{$\sc2$};
\draw (3-0.5,1-0.5)node{$\sc\ov2$};\draw (3-0.5,2-0.5)node{$\sc\ov1$};\draw (3-0.5,3-0.5)node{$\sc0$};
\draw (4-0.5,1-0.5)node{$\sc\ov3$};
\end{tikzpicture}
}}
&
\vcenter{\hbox{
\begin{tikzpicture}[x={(0in,-0.15in)},y={(0.15in,0in)}] 
\foreach \j in {1,...,3} \draw (1,\j) rectangle +(-1,-1);
\foreach \j in {1,...,3} \draw (2,\j) rectangle +(-1,-1);
\foreach \j in {1,...,3} \draw (3,\j) rectangle +(-1,-1);
\foreach \j in {1} \draw (4,\j) rectangle +(-1,-1);
\foreach \j in {1,2,3} \draw(1-0.5,\j-0.5)node{$\ast$};
\foreach \j in {1,2} \draw(2-0.5,\j-0.5)node{$\ast$};
\draw (2-0.5,3-0.5)node{$\sc1$};
\draw (3-0.5,1-0.5)node{$\sc\ov2$};\draw (3-0.5,2-0.5)node{$\sc\ov1$};\draw (3-0.5,3-0.5)node{$\sc0$};
\draw (4-0.5,1-0.5)node{$\sc\ov3$};
\end{tikzpicture}
}}
&
\vcenter{\hbox{
\begin{tikzpicture}[x={(0in,-0.15in)},y={(0.15in,0in)}] 
\foreach \j in {1,...,3} \draw (1,\j) rectangle +(-1,-1);
\foreach \j in {1,...,4} \draw (2,\j) rectangle +(-1,-1);
\foreach \j in {1,...,3} \draw (3,\j) rectangle +(-1,-1);
\foreach \j in {1} \draw (4,\j) rectangle +(-1,-1);
\foreach \j in {1,2,3} \draw(1-0.5,\j-0.5)node{$\ast$};
\foreach \j in {1,2} \draw(2-0.5,\j-0.5)node{$\ast$};
\draw (2-0.5,3-0.5)node{$\sc1$};\draw (2-0.5,4-0.5)node{$\sc2$};
\draw (3-0.5,1-0.5)node{$\sc\ov2$};\draw (3-0.5,2-0.5)node{$\sc\ov1$};\draw (3-0.5,3-0.5)node{$\sc0$};
\draw (4-0.5,1-0.5)node{$\sc\ov3$};
\end{tikzpicture}
}}
&
\vcenter{\hbox{
\begin{tikzpicture}[x={(0in,-0.15in)},y={(0.15in,0in)}] 
\foreach \j in {1,...,3} \draw (1,\j) rectangle +(-1,-1);
\foreach \j in {1,...,2} \draw (2,\j) rectangle +(-1,-1);
\foreach \j in {1,...,2} \draw (3,\j) rectangle +(-1,-1);
\foreach \j in {1} \draw (4,\j) rectangle +(-1,-1);
\foreach \j in {1,2,3} \draw(1-0.5,\j-0.5)node{$\ast$};
\foreach \j in {1,2} \draw(2-0.5,\j-0.5)node{$\ast$};
\draw (3-0.5,1-0.5)node{$\sc\ov2$};\draw (3-0.5,2-0.5)node{$\sc\ov1$};
\draw (4-0.5,1-0.5)node{$\sc\ov3$};
\end{tikzpicture}
}}
&
\vcenter{\hbox{
\begin{tikzpicture}[x={(0in,-0.15in)},y={(0.15in,0in)}] 
\foreach \j in {1,...,3} \draw (1,\j) rectangle +(-1,-1);
\foreach \j in {1,...,2} \draw (2,\j) rectangle +(-1,-1);
\foreach \j in {1,...,1} \draw (3,\j) rectangle +(-1,-1);
\foreach \j in {1} \draw (4,\j) rectangle +(-1,-1);
\foreach \j in {1,2,3} \draw(1-0.5,\j-0.5)node{$\ast$};
\foreach \j in {1,2} \draw(2-0.5,\j-0.5)node{$\ast$};
\draw (3-0.5,1-0.5)node{$\sc\ov2$};
\draw (4-0.5,1-0.5)node{$\sc\ov3$};
\end{tikzpicture}
}}
\cr
\end{array}
	\,\right]
\end{equation}
and in the outer rim case
\begin{equation}\label{eqn-outer-diagrams}
\left[\, 
\begin{array}{ccccc} 
\vcenter{\hbox{
\begin{tikzpicture}[x={(0in,-0.15in)},y={(0.15in,0in)}] 
\foreach \j in {1,...,5} \draw (1,\j) rectangle +(-1,-1);
\foreach \j in {1,...,4} \draw (2,\j) rectangle +(-1,-1);
\foreach \j in {1,...,4} \draw (3,\j) rectangle +(-1,-1);
\foreach \j in {1,2} \draw (4,\j) rectangle +(-1,-1);
\foreach \j in {1,2,3} \draw(1-0.5,\j-0.5)node{$\ast$};
\foreach \j in {1,2,3} \draw(2-0.5,\j-0.5)node{$\ast$};
\foreach \j in {1} \draw(3-0.5,\j-0.5)node{$\ast$};
\draw (4-0.5,1-0.5)node{$\sc\ov3$};\draw (4-0.5,2-0.5)node{$\sc\ov2$};
\draw (3-0.5,2-0.5)node{$\sc\ov1$};\draw (3-0.5,3-0.5)node{$\sc0$};\draw (3-0.5,4-0.5)node{$\sc1$};
\draw (2-0.5,4-0.5)node{$\sc2$};
\draw (1-0.5,4-0.5)node{$\sc3$};\draw (1-0.5,5-0.5)node{$\sc4$};
\end{tikzpicture}
}}
&
\vcenter{\hbox{
\begin{tikzpicture}[x={(0in,-0.15in)},y={(0.15in,0in)}] 
\foreach \j in {1,...,5} \draw (1,\j) rectangle +(-1,-1);
\foreach \j in {1,...,4} \draw (2,\j) rectangle +(-1,-1);
\foreach \j in {1,...,4} \draw (3,\j) rectangle +(-1,-1);
\foreach \j in {1,2} \draw (4,\j) rectangle +(-1,-1);
\foreach \j in {1,...,5} \draw(1-0.5,\j-0.5)node{$\ast$};
\foreach \j in {1,...,4} \draw(2-0.5,\j-0.5)node{$\ast$};
\foreach \j in {1} \draw(3-0.5,\j-0.5)node{$\ast$};
\draw (4-0.5,1-0.5)node{$\sc\ov3$};\draw (4-0.5,2-0.5)node{$\sc\ov2$};
\draw (3-0.5,2-0.5)node{$\sc\ov1$};\draw (3-0.5,3-0.5)node{$\sc0$};\draw (3-0.5,4-0.5)node{$\sc1$};
\end{tikzpicture}
}}
&
\vcenter{\hbox{
\begin{tikzpicture}[x={(0in,-0.15in)},y={(0.15in,0in)}] 
\foreach \j in {1,...,5} \draw (1,\j) rectangle +(-1,-1);
\foreach \j in {1,...,4} \draw (2,\j) rectangle +(-1,-1);
\foreach \j in {1,...,4} \draw (3,\j) rectangle +(-1,-1);
\foreach \j in {1,2} \draw (4,\j) rectangle +(-1,-1);
\foreach \j in {1,...,5} \draw(1-0.5,\j-0.5)node{$\ast$};
\foreach \j in {1,2,3} \draw(2-0.5,\j-0.5)node{$\ast$};
\foreach \j in {1} \draw(3-0.5,\j-0.5)node{$\ast$};
\draw (4-0.5,1-0.5)node{$\sc\ov3$};\draw (4-0.5,2-0.5)node{$\sc\ov2$};
\draw (3-0.5,2-0.5)node{$\sc\ov1$};\draw (3-0.5,3-0.5)node{$\sc0$};\draw (3-0.5,4-0.5)node{$\sc1$};
\draw (2-0.5,4-0.5)node{$\sc2$};
\end{tikzpicture}
}}
&
\vcenter{\hbox{
\begin{tikzpicture}[x={(0in,-0.15in)},y={(0.15in,0in)}] 
\foreach \j in {1,...,5} \draw (1,\j) rectangle +(-1,-1);
\foreach \j in {1,...,4} \draw (2,\j) rectangle +(-1,-1);
\foreach \j in {1,...,4} \draw (3,\j) rectangle +(-1,-1);
\foreach \j in {1,2} \draw (4,\j) rectangle +(-1,-1);
\foreach \j in {1,...,5} \draw(1-0.5,\j-0.5)node{$\ast$};
\foreach \j in {1,...,4} \draw(2-0.5,\j-0.5)node{$\ast$};
\foreach \j in {1,3,4} \draw(3-0.5,\j-0.5)node{$\ast$};
\draw (4-0.5,1-0.5)node{$\sc\ov3$};\draw (4-0.5,2-0.5)node{$\sc\ov2$};
\draw (3-0.5,2-0.5)node{$\sc\ov1$};\
\end{tikzpicture}
}}
&
\vcenter{\hbox{
\begin{tikzpicture}[x={(0in,-0.15in)},y={(0.15in,0in)}] 
\foreach \j in {1,...,5} \draw (1,\j) rectangle +(-1,-1);
\foreach \j in {1,...,4} \draw (2,\j) rectangle +(-1,-1);
\foreach \j in {1,...,4} \draw (3,\j) rectangle +(-1,-1);
\foreach \j in {1,2} \draw (4,\j) rectangle +(-1,-1);
\foreach \j in {1,...,5} \draw(1-0.5,\j-0.5)node{$\ast$};
\foreach \j in {1,...,4} \draw(2-0.5,\j-0.5)node{$\ast$};
\foreach \j in {1,...,4} \draw(3-0.5,\j-0.5)node{$\ast$};
\draw (4-0.5,1-0.5)node{$\sc\ov3$};\draw (4-0.5,2-0.5)node{$\sc\ov2$};
\end{tikzpicture}
}}
\cr
\end{array}
	\,\right]
\end{equation}
These immediately give the identities
\begin{equation}\label{eqn-inner-ex}
%
\,\right|\cr
\end{array}
\end{equation}
where advantage has been taken of the fact that $s_{\kappa/\mu}(\X)=s_{\lambda/\kappa}(\X)=0$ if 
$\kappa=(\ldots,k,k+1,\ldots,)$ or $\kappa'=(\ldots,l,l+1,\ldots)$, while $s_{\kappa/\mu}(\X)=s_{\lambda/\lambda}(\X)=1$. 
Each null strip gives rise to a row or column of the determinant containing
a single entry $1$ with all other entries $0$, so that the null strips play no role in the
evaluation of $s_{(\theta_i\#\theta_j)}(\X)$ and might have been omitted from the outset. 
However they did play a crucial role in the derivation of the Corollaries~\ref{Cor-OkaX} and \ref{Cor-LPX}.

\subsection{Pfaffian corollaries}\label{Subsec-pfaff-cor}

Throughout this section, $\lambda$ and $\mu$ are strict partitions of lengths $p$ and $q$, respectively with $\mu\subseteq\lambda$. 
The image of of the skew shifted diagram $SF^{\lambda/\mu}$ under reflection in its SW to NE anti-diagonal through the box at 
position $(1,\lambda_1)$ is the skew shifted diagram $SF^{\kappa/\nu}$, with $\kappa$ and $\nu$ defined in this way strict 
partition of lengths $r$ and $s$, respectively. For example, in the case $\lambda=(6,5,3)$ and $\mu=(3,2)$ the reflection
in the anti-diagonal of $SF^{\lambda/\mu}$ yields $SF^{\kappa/\nu}$ with $\kappa=(6,5,4,1)$ and $\nu=(4,2,1)$ as illustrated below
\begin{equation}
SF^{\lambda/\mu}=
\vcenter{\hbox{
\begin{tikzpicture}[x={(0in,-0.15in)},y={(0.15in,0in)}] 
\foreach \j in {1,...,6} \draw[thick] (1,\j) rectangle +(-1,-1);
\foreach \j in {2,...,6} \draw[thick] (2,\j) rectangle +(-1,-1);
\foreach \j in {3,...,5} \draw[thick] (3,\j) rectangle +(-1,-1);\foreach \j in {6} \draw (3,\j) rectangle +(-1,-1);
\foreach \j in {4,...,6} \draw (4,\j) rectangle +(-1,-1);
\foreach \j in {5,...,6} \draw (5,\j) rectangle +(-1,-1);
\foreach \j in {6} \draw (6,\j) rectangle +(-1,-1);
\foreach \j in {1,...,3} \draw(0.5,\j-0.5)node{$\ast$};
\foreach \j in {2,...,3} \draw(1.5,\j-0.5)node{$\ast$};
\draw[ultra thick](2,2)--(2,3)--(0,3)--(0,6)--(2,6)--(2,5)--(3,5)--(3,2)--(2,2);
\draw (4,2)--(-1,7);
\foreach \i in {4,...,4} \draw(\i-0.5,3.5)node{$\ast$};
\foreach \i in {4,...,5} \draw(\i-0.5,4.5)node{$\ast$};
\foreach \i in {3,...,6} \draw(\i-0.5,5.5)node{$\ast$};
\end{tikzpicture}
}}
\qquad
SF^{\kappa/\nu}=
\vcenter{\hbox{
\begin{tikzpicture}[x={(0in,-0.15in)},y={(0.15in,0in)}] 
\foreach \j in {1,...,6} \draw[thick] (1,\j) rectangle +(-1,-1);
\foreach \j in {2,...,6} \draw[thick] (2,\j) rectangle +(-1,-1);
\foreach \j in {3,...,6} \draw[thick] (3,\j) rectangle +(-1,-1);
\foreach \j in {4} \draw[thick](4,\j) rectangle +(-1,-1);\foreach \j in {5,6} \draw(4,\j) rectangle +(-1,-1);
\foreach \j in {5,...,6} \draw(5,\j) rectangle +(-1,-1);
\foreach \j in {6} \draw(6,\j) rectangle +(-1,-1);
\foreach \j in {1,...,4} \draw(0.5,\j-0.5)node{$\ast$};
\foreach \j in {2,...,3} \draw(1.5,\j-0.5)node{$\ast$};
\foreach \j in {3} \draw(2.5,\j-0.5)node{$\ast$};
\draw[ultra thick](3,3)--(1,3)--(1,4)--(0,4)--(0,6)--(3,6)--(3,4)--(4,4)--(4,3)--(3,3);
\draw(4,2)--(-0.5,6.5);
\foreach \i in {4,...,5} \draw(\i-0.5,4.5)node{$\ast$};
\foreach \i in {4,...,6} \draw(\i-0.5,5.5)node{$\ast$};
\end{tikzpicture}
}}
\end{equation}


For each outside decomposition $\Theta$ of $SF^{\lambda/\mu}$ the corresponding cutting strip $\phi$ extends
from a box of content $0$ to one of content $\ell'-1$ where $\ell'=\lambda_1$ is the number of boxes in the first row of $SF^\lambda$.
The constituent strips of the outside decomposition are necessarily of the form $\theta_k\simeq\phi_{ab}$ with $a$ of the form $\mu_i$
or $\nu_i$, and $b$ of the form $\lambda_i-1$ or $\kappa_j-1$, where it is sometimes convenient to extend the various partitions 
$\lambda$, $\mu$, $\kappa$ and $\nu$ through the inclusion of a single trailing $0$.

In the case $\lambda=(6,5,3)$ and $\mu=(3,2)$ the various outside decompositions are illustrated below 
that arise from the superposition on $SF^{\lambda/\mu}$ of translations of cutting strips $\phi$ in the form 
of a row, a column, an inner rim and an outer rim.

\begin{equation}
\begin{array}{|c|c|c|c|}
\hline
\mbox{Row}&\mbox{Column}&\mbox{Inner rim}&\mbox{Outer rim}\cr
\hline
&&&\cr
\vcenter{\hbox{
\begin{tikzpicture}[x={(0in,-0.15in)},y={(0.15in,0in)}] 
\foreach \j in {1,...,6} \draw (1,\j) rectangle +(-1,-1);
\foreach \j in {2,...,6} \draw (2,\j) rectangle +(-1,-1);
\foreach \j in {3,...,5} \draw (3,\j) rectangle +(-1,-1);
\fill[blue!40!white](0.3,-0.5)--(0.3,6.5)--(0.7,6.5)--(0.7,-0.5)--cycle;
\fill[red!40!white](1.3,1.5)--(1.3,6.5)--(1.7,6.5)--(1.7,1.5)--cycle;
\fill[cyan!40!white](2.3,2.5)--(2.3,6.5)--(2.7,6.5)--(2.7,2.5)--cycle;
\foreach \j in {1,...,3} \draw(0.5,\j-0.5)node{$\ast$};
\foreach \j in {2,...,3} \draw(1.5,\j-0.5)node{$\ast$};
\end{tikzpicture}
}}
&
\vcenter{\hbox{
\begin{tikzpicture}[x={(0in,-0.15in)},y={(0.15in,0in)}] 
\foreach \j in {1,...,6} \draw (1,\j) rectangle +(-1,-1);
\foreach \j in {2,...,6} \draw (2,\j) rectangle +(-1,-1);
\foreach \j in {3,...,5} \draw (3,\j) rectangle +(-1,-1);
\fill[blue!40!white](-0.5,0.3)--(1.5,0.3)--(1.5,0.7)--(-0.5,0.7)--cycle;
\fill[red!40!white](-0.5,1.3)--(2.5,1.3)--(2.5,1.7)--(-0.5,1.7)--cycle;
\fill[cyan!40!white](-0.5,2.3)--(3.5,2.3)--(3.5,2.7)--(-0.5,2.7)--cycle;
\fill[magenta!40!white](-0.5,3.3)--(3.5,3.3)--(3.5,3.7)--(-0.5,3.7)--cycle;
\fill[orange!40!white](-0.5,4.3)--(3.5,4.3)--(3.5,4.7)--(-0.5,4.7)--cycle;
\fill[green!40!white](-0.5,5.3)--(3.5,5.3)--(3.5,5.7)--(-0.5,5.7)--cycle;
\foreach \j in {1,...,3} \draw(0.5,\j-0.5)node{$\ast$};
\foreach \j in {2,...,3} \draw(1.5,\j-0.5)node{$\ast$};
\end{tikzpicture}
}}
&
\vcenter{\hbox{
\begin{tikzpicture}[x={(0in,-0.15in)},y={(0.15in,0in)}] 
\foreach \j in {1,...,6} \draw (1,\j) rectangle +(-1,-1);
\foreach \j in {2,...,6} \draw (2,\j) rectangle +(-1,-1);
\foreach \j in {3,...,5} \draw (3,\j) rectangle +(-1,-1);
\fill[blue!40!white](3.5,2.3)--(2.3,2.3)--(2.3,3.3)--(0.3,3.3)--(0.3,6.5)--(0.7,6.5)--(0.7,3.7)--(2.7,3.7)--(2.7,2.7)--(3.5,2.7)--cycle;
\fill[red!40!white](4.5,3.3)--(3.3,3.3)--(3.3,4.3)--(1.3,4.3)--(1.3,6.5)--(1.7,6.5)--(1.7,4.7)--(3.7,4.7)--(3.7,3.7)--(4.5,3.7)--cycle;
\foreach \j in {1,...,3} \draw(0.5,\j-0.5)node{$\ast$};
\foreach \j in {2,...,3} \draw(1.5,\j-0.5)node{$\ast$};
\end{tikzpicture}
}}
&
\vcenter{\hbox{
\begin{tikzpicture}[x={(0in,-0.15in)},y={(0.15in,0in)}] 
\foreach \j in {1,...,6} \draw (1,\j) rectangle +(-1,-1);
\foreach \j in {2,...,6} \draw (2,\j) rectangle +(-1,-1);
\foreach \j in {3,...,5} \draw (3,\j) rectangle +(-1,-1);
\fill[blue!40!white](3.5,2.3)--(2.3,2.3)--(2.3,4.3)--(1.3,4.3)--(1.3,5.3)--(0.3,5.3)--(0.3,6.5)--(0.7,6.5)--(0.7,5.7)--(1.7,5.7)--(1.7,4.7)--(2.7,4.7)--(2.7,2.7)--(3.5,2.7)--cycle;
\fill[red!40!white](2.5,1.3)--(1.3,1.3)--(1.3,3.3)--(0.3,3.3)--(0.3,4.3)--(-0.7,4.3)--(-0.7,5.5)--(-0.3,5.5)--(-0.3,4.7)--(0.7,4.7)--(0.7,3.7)--(1.7,3.7)--(1.7,1.7)--(2.5,1.7)--cycle;
\foreach \j in {1,...,3} \draw(0.5,\j-0.5)node{$\ast$};
\foreach \j in {2,...,3} \draw(1.5,\j-0.5)node{$\ast$};
\end{tikzpicture}
}}
\cr
&&&\cr
\hline
&&&\cr
\vcenter{\hbox{
\begin{tikzpicture}[x={(0in,-0.15in)},y={(0.15in,0in)}] 
\foreach \j in {1,...,6} \draw (1,\j) rectangle +(-1,-1);
\foreach \j in {2,...,6} \draw (2,\j) rectangle +(-1,-1);
\foreach \j in {3,...,5} \draw (3,\j) rectangle +(-1,-1);
\draw[draw=blue,ultra thick](0.5,2.8)--(0.5,6.2);
\draw[draw=red,ultra thick] (1.5,2.8)--(1.5,6.2);
\draw[draw=cyan,ultra thick] (2.5,1.8)--(2.5,5.2);
\foreach \j in {1,...,3} \draw(0.5,\j-0.5)node{$\ast$};
\foreach \j in {2,...,3} \draw(1.5,\j-0.5)node{$\ast$};
\end{tikzpicture}
}}
&
\vcenter{\hbox{
\begin{tikzpicture}[x={(0in,-0.15in)},y={(0.15in,0in)}] 
\foreach \j in {1,...,6} \draw (1,\j) rectangle +(-1,-1);
\foreach \j in {2,...,6} \draw (2,\j) rectangle +(-1,-1);
\foreach \j in {3,...,5} \draw (3,\j) rectangle +(-1,-1);
\draw[draw=cyan,ultra thick] (3.2,2.5)--(1.8,2.5);
\draw[draw=magenta,ultra thick](3.2,3.5)--(-0.2,3.5);
\draw[draw=orange,ultra thick] (3.2,4.5)--(-0.2,4.5);
\draw[draw=green,ultra thick] (2.2,5.5)--(-0.2,5.5);
\foreach \j in {1,...,3} \draw(0.5,\j-0.5)node{$\ast$};
\foreach \j in {2,...,3} \draw(1.5,\j-0.5)node{$\ast$};
\end{tikzpicture}
}}
&
\vcenter{\hbox{
\begin{tikzpicture}[x={(0in,-0.15in)},y={(0.15in,0in)}] 
\foreach \j in {1,...,6} \draw (1,\j) rectangle +(-1,-1);
\foreach \j in {2,...,6} \draw (2,\j) rectangle +(-1,-1);
\foreach \j in {3,...,5} \draw (3,\j) rectangle +(-1,-1);
\draw[draw=blue,ultra thick](3.2,2.5)--(2.5,2.5)--(2.5,3.5)--(0.5,3.5)--(0.5,6.2);
\draw[draw=red,ultra thick] (3.2,4.5)--(1.5,4.5)--(1.5,6.2);
\foreach \j in {1,...,3} \draw(0.5,\j-0.5)node{$\ast$};
\foreach \j in {2,...,3} \draw(1.5,\j-0.5)node{$\ast$};
\end{tikzpicture}
}}
&
\vcenter{\hbox{
\begin{tikzpicture}[x={(0in,-0.15in)},y={(0.15in,0in)}] 
\foreach \j in {1,...,6} \draw (1,\j) rectangle +(-1,-1);
\foreach \j in {2,...,6} \draw (2,\j) rectangle +(-1,-1);
\foreach \j in {3,...,5} \draw (3,\j) rectangle +(-1,-1);
\draw[draw=blue,ultra thick](3.2,2.5)--(2.5,2.5)--(2.5,4.5)--(1.5,4.5)--(1.5,5.5)--(0.5,5.5)--(0.5,6.2);
\draw[draw=red,ultra thick] (1.5,2.8)--(1.5,2.8)--(1.5,3.5)--(0.5,3.5)--(0.5,4.5)--(-0.2,4.5);
\foreach \j in {1,...,3} \draw(0.5,\j-0.5)node{$\ast$};
\foreach \j in {2,...,3} \draw(1.5,\j-0.5)node{$\ast$};
\end{tikzpicture}
}}
\cr
&&&\cr
\hline
\end{array}
\end{equation}

Generalising each of these one arrives at the following corollaries of Theorem~\ref{The-HamXY}.

\begin{Corollary}[J\'osefiak-Pragacz-Nimmo type identity \cite{JP}, \cite{Nim}]\label{Cor-Qrow}
Let $\lambda$ and $\mu$ be strict partitions of lengths $p$ and $q$ with $\mu\subseteq\lambda$
and let $\X=(x_{kc})$ and $\Y=(y_{kc})$ with $k=1,2,\ldots,n$ and $c\in\Z$. Then
\begin{equation}\label{eqn-Qrow}
    Q_{\lambda/\mu}(\X,\!\Y) = \pf\left( \begin{array}{cc} Q_{(\lambda_i,\lambda_j)}(\X,\!\Y)&Q_{(\lambda_i-\mu_{m-k+1})}(\tau^{\mu_{m-k+1}}(\X,\!\Y))\cr
		                                      -\,^t(Q_{(\lambda_i-\mu_{m-k+1})}(\tau^{\mu_{m-k+1}}(\X,\!\Y))&0\cr \end{array} \right)  
\end{equation} 
with $1\leq i,j\leq \ell$ and $1\leq k\leq m$, where $\ell=p$ or $p+1$ according as $p+q$ is even or odd respectively, and $m=q$. 
Here $Q_{(\lambda_i,\lambda_j)}(\X,\!\Y)=-Q_{(\lambda_j,\lambda_i)}(\X,\!\Y)$
so that $Q_{(\lambda_i,\lambda_i)}(\X,\!\Y)=0$, and $Q_{(\lambda_i-\mu_{m-k+1})}(\tau^{\mu_{m-k+1}}(\X,\!\Y))=1$ or $0$ if 
$\lambda_i=\mu_{m-k+1}$ or $\lambda_i<\mu_{m-k+1}$, respectively. 
\end{Corollary}

\noindent{\bf Proof}: Here one takes $\phi$ to be a horizontal strip, that is a row extending from a box of content $0$ to one of content $\lambda_1-1$.
Numbering the horizontal strips $\theta_i$ from top to bottom one has $\theta_i\simeq\phi_{\mu_i,\lambda_i-1}$ for $1=1,2,\ldots,q$.
and $\theta_i\simeq\phi_{0,\lambda_i-1}$ for $i=q+1,q+2,\ldots,p$.
Hence $\ov{\theta}_i\simeq\phi_{0,\lambda_i-1}$ for $i=1,2,\ldots,p$, so that $(\ov{\theta}_i,\ov{\theta}_j)=(\phi_{0,\lambda_i-1},\phi_{0,\lambda_j-1})$
is the two-rowed shifted diagram $SF^{\lambda_i,\lambda_j}$ for $1\leq i,j\leq p$, which is regular if and only if $i<j$. 
Similarly, if $1\leq j\leq q$ and $1\leq i\leq p$ we have $(\theta_j\#\ov{\theta}_i)=\phi_{\mu_j,\lambda_i-1}$ which is the row diagram $SF^{(\lambda_i-\mu_j)}$
in which the content of the leftmost box is $\mu_j$. The use of this data in Theorem~\ref{The-HamXY} then yields (\ref{eqn-Qrow}), a result
first established by other means in the case $x_{kc}=y_{kc}=x_k$ for all $k$ and $c$ by J\'osefiak and Pragacz~\cite{JP}, 
and precisely as above in the same special case by Hamel~\cite{Ham}. 
\qed

\begin{Corollary}\label{Cor-Qcol}
Let $\lambda$ and $\mu$ be strict partitions with $\mu\subseteq\lambda$ such that the image of $SF^{\lambda/\mu}$
under anti-diagonal reflection is $SF^{\kappa/\nu}$.
Let $\X=(x_{kc})$, $\Y=(y_{kc})$, $\widetilde{\X}=(x_{n-k+1},c)$ and $\widetilde{\Y}=(y_{n-k+1},c)$ with $k=1,2,\ldots,n$ and $c\in\Z$. 
Then, recalling that we denote the lengths of $\lambda$, $\mu$, $\kappa$ and $\nu$ by $p$, $q$, $r$ and $s$,
\begin{equation}\label{eqn-Qcol}
    Q_{\lambda/\mu}(\X,\!\Y) = \pf\left( \begin{array}{cc} Q_{(\kappa_i,\kappa_j)}(\widetilde{\Y},\!\widetilde{\X})                                                               &Q_{(\kappa_i-\nu_{m-k+1})}(\tau^{\nu_{m-k+1}}(\widetilde{\Y},\!\widetilde{\X})\cr
		                                      -\,^t(Q_{(\kappa_i-\nu_{m-k+1})}(\tau_{\nu_{m-k+1}}(\widetilde{\Y},\!\widetilde{\X})&0\cr \end{array} \right)  
\end{equation} 
with $1\leq i,j\leq \ell$ and $1\leq k\leq m$, where $\ell=r$ or $r+1$ according as $r+s$ is even or odd, respectively,
and $m=s$.
Here $Q_{(\kappa_i,\kappa_j)}(\widetilde{\Y},\!\widetilde{\X})=-Q_{(\kappa_j,\kappa_i)}(\widetilde{\Y},\!\widetilde{\X})$
so that $Q_{(\kappa_i,\kappa_i)}(\widetilde{\Y},\!\widetilde{\X})=0$, and $Q_{(\kappa_i-\nu_{m-k+1})}(\tau^{\nu_{m-k+1}}(\widetilde{\Y},\!\widetilde{\X})=1$ or $0$ if 
$\kappa_i=\nu_{m-k+1}$ or $\kappa_i<\nu_{m-k+1}$, respectively. 
\end{Corollary}

\noindent{\bf Proof}: The appropriate cutting strip $\phi$ is a vertical strip extending from a box of content $0$ to one of content $\ell'=\lambda_1-1$.
Numbering the vertical strips $\theta_i$ from right to left one has $\theta_i\simeq\phi_{\nu_i,\kappa_i-1}$ for $1=1,2,\ldots,s$.
and $\theta_i\simeq\phi_{0,\kappa_i-1}$ for $i=s+1,q+2,\ldots,r$.
Hence $\ov{\theta}_i\simeq\phi_{0,\kappa_i-1}$ for $i=1,2,\ldots,r$, so that $(\ov{\theta}_i,\ov{\theta}_j)=(\phi_{0,\kappa_i-1},\phi_{0,\kappa_j-1})$
is the two-column shifted diagram obtained by anti-diagonal reflection of $SF^{\kappa_i,\kappa_j}$ for $1\leq i,j\leq r$, which is regular if and only if $i>j$. 
Similarly, if $1\leq k\leq m$ and $1\leq i\leq \ell$ we have $(\theta_k\#\ov{\theta}_i)=\phi_{\nu_k,\kappa_i-1}$ which is image of the row diagram 
$SF^{(\kappa_i-\nu_k)}$ under anti-diagonal reflection, with the content of the leftmost box given by $\nu_k$. The required result (\ref{eqn-Qcol})
then follows from Theorem~\ref{The-HamXY}. 
\qed

Alternatively, Corollary~\ref{Cor-Qcol} can be seen to follow directly from Corollary~\ref{Cor-Qrow} through use of the following:
\begin{Lemma}\label{Lem-Qrefl}
Let $\lambda$ and $\mu$ be strict partitions with $\mu\subseteq\lambda$ such that the image of $SF^{\lambda/\mu}$
under anti-diagonal reflection is $SF^{\kappa/\nu}$, and let $\X=(x_{kc})$, $\Y=(y_{kc})$, $\widetilde{\X}=(x_{n-k+1},c)$ 
and $\widetilde{\Y}=(y_{n-k+1},c)$ with $k=1,2,\ldots,n$ and $c\in\Z$. Then
\begin{equation}
       Q_{\lambda/\mu}(\X,\!\Y) =  Q_{\kappa/\nu}(\widetilde{\Y},\!\widetilde{\X})
\end{equation}
\end{Lemma}.

\noindent{\bf Proof}: Under the anti-diagonal reflection map from $SF^{\lambda/\mu}$ to $SF^{\kappa/\nu}$ the content of boxes is 
preserved but the image of each skew shifted tableau $T\in{\cal ST}^{\lambda/\mu}$ under the same reflection will be non-standard.
However, if the reflection is augmented by mapping each entry $k$ to $(n-k+1)'$ and $k'$ to $(n-k+1)$ the result will be 
a standard skew shifted tableau $\widetilde{T}\in{\cal ST}^{\kappa/\nu}$. This augmented map is bijective,
and in order that it be weight preserving it suffices to map weights $x_{kc}$ and $y_{kc}$ to $y_{(n-k+1),c}$ and $x_{(n-k+1),c}$,
respectively. \qed

This is illustrated below:
\begin{equation}\nonumber
\begin{array}{cc}
F^{\lambda/\mu}=
\vcenter{\hbox{
\begin{tikzpicture}[x={(0in,-0.2in)},y={(0.2in,0in)}] 
\foreach \j in {1,...,6} \draw[thick] (1,\j) rectangle +(-1,-1);
\foreach \j in {2,...,6} \draw[thick] (2,\j) rectangle +(-1,-1);
\foreach \j in {3,...,5} \draw[thick] (3,\j) rectangle +(-1,-1);
\foreach \j in {1,...,3} \draw(0.5,\j-0.5)node{$\ast$};
\foreach \j in {2,...,3} \draw(1.5,\j-0.5)node{$\ast$};
\draw (1-0.5,4-0.5)node{$\sc3$};\draw (1-0.5,5-0.5)node{$\sc4$};\draw (1-0.5,6-0.5)node{$\sc5$};
\draw (2-0.5,4-0.5)node{$\sc2$};\draw (2-0.5,5-0.5)node{$\sc3$};\draw (2-0.5,6-0.5)node{$\sc4$};
\draw (3-0.5,3-0.5)node{$\sc0$};\draw (3-0.5,4-0.5)node{$\sc1$};\draw (3-0.5,5-0.5)node{$\sc2$};
\end{tikzpicture}
}}
&
F^{\kappa/\nu}=
\vcenter{\hbox{
\begin{tikzpicture}[x={(0in,-0.2in)},y={(0.2in,0in)}] 
\foreach \j in {1,...,6} \draw[thick] (1,\j) rectangle +(-1,-1);
\foreach \j in {2,...,6} \draw[thick] (2,\j) rectangle +(-1,-1);
\foreach \j in {3,...,6} \draw[thick] (3,\j) rectangle +(-1,-1);
\foreach \j in {4} \draw[thick](4,\j) rectangle +(-1,-1);
\foreach \j in {1,...,4} \draw(0.5,\j-0.5)node{$\ast$};
\foreach \j in {2,...,3} \draw(1.5,\j-0.5)node{$\ast$};
\foreach \j in {3} \draw(2.5,\j-0.5)node{$\ast$};
\draw (1-0.5,5-0.5)node{$\sc4$};\draw (1-0.5,6-0.5)node{$\sc5$};
\draw (2-0.5,4-0.5)node{$\sc2$};\draw (2-0.5,5-0.5)node{$\sc3$};\draw (2-0.5,6-0.5)node{$\sc4$};
\draw (3-0.5,4-0.5)node{$\sc1$};\draw (3-0.5,5-0.5)node{$\sc2$};\draw (3-0.5,6-0.5)node{$\sc3$};
\draw (4-0.5,4-0.5)node{$\sc0$};
\end{tikzpicture}
}}\cr\cr
T\in{\cal ST}^{\lambda/\mu}=
\vcenter{\hbox{
\begin{tikzpicture}[x={(0in,-0.2in)},y={(0.2in,0in)}] 
\foreach \j in {1,...,6} \draw[thick] (1,\j) rectangle +(-1,-1);
\foreach \j in {2,...,6} \draw[thick] (2,\j) rectangle +(-1,-1);
\foreach \j in {3,...,5} \draw[thick] (3,\j) rectangle +(-1,-1);\foreach \j in {6} 
\foreach \j in {1,...,3} \draw(0.5,\j-0.5)node{$\ast$};
\foreach \j in {2,...,3} \draw(1.5,\j-0.5)node{$\ast$};
\draw (1-0.5,4-0.5)node{$\sc1$};\draw (1-0.5,5-0.5)node{$\sc2$};\draw (1-0.5,6-0.5)node{$\sc2$};
\draw (2-0.5,4-0.5)node{$\sc2'$};\draw (2-0.5,5-0.5)node{$\sc3'$};\draw (2-0.5,6-0.5)node{$\sc4$};
\draw (3-0.5,3-0.5)node{$\sc1'$};\draw (3-0.5,4-0.5)node{$\sc3$};\draw (3-0.5,5-0.5)node{$\sc3$};
\end{tikzpicture}
}}
&
\widetilde{T}\in{\cal ST}^{\kappa/\nu}=
\vcenter{\hbox{
\begin{tikzpicture}[x={(0in,-0.2in)},y={(0.2in,0in)}] 
\foreach \j in {1,...,6} \draw[thick] (1,\j) rectangle +(-1,-1);
\foreach \j in {2,...,6} \draw[thick] (2,\j) rectangle +(-1,-1);
\foreach \j in {3,...,6} \draw[thick] (3,\j) rectangle +(-1,-1);
\foreach \j in {4} \draw[thick](4,\j) rectangle +(-1,-1); 
\foreach \j in {1,...,4} \draw(0.5,\j-0.5)node{$\ast$};
\foreach \j in {2,...,3} \draw(1.5,\j-0.5)node{$\ast$};
\foreach \j in {3} \draw(2.5,\j-0.5)node{$\ast$};
\draw (1-0.5,5-0.5)node{$\sc1'$};\draw (1-0.5,6-0.5)node{$\sc3'$};
\draw (2-0.5,4-0.5)node{$\sc2'$};\draw (2-0.5,5-0.5)node{$\sc2$};\draw (2-0.5,6-0.5)node{$\sc3'$};
\draw (3-0.5,4-0.5)node{$\sc2'$};\draw (3-0.5,5-0.5)node{$\sc3$};\draw (3-0.5,6-0.5)node{$\sc4'$};
\draw (4-0.5,4-0.5)node{$\sc4$};
\end{tikzpicture}
}}
\cr\cr
\wgt(T)=
\vcenter{\hbox{
\begin{tikzpicture}[x={(0in,-0.2in)},y={(0.25in,0in)}] 
\draw (1-0.5,4-0.5)node{$  x_{13}$};\draw (1-0.5,5-0.5)node{$  x_{24}$};\draw (1-0.5,6-0.5)node{$  x_{25}$};
\draw (2-0.5,4-0.5)node{$  y_{22}$};\draw (2-0.5,5-0.5)node{$  y_{33}$};\draw (2-0.5,6-0.5)node{$  x_{44}$};
\draw (3-0.5,3-0.5)node{$  y_{10}$};\draw (3-0.5,4-0.5)node{$  x_{32}$};\draw (3-0.5,5-0.5)node{$  x_{33}$};
\end{tikzpicture}
}}
&
\wgt(\widetilde{T})=
\vcenter{\hbox{
\begin{tikzpicture}[x={(0in,-0.2in)},y={(0.25in,0in)}] 
\draw (1-0.5,5-0.5)node{$ x_{44}$};\draw (1-0.5,6-0.5)node{$ x_{25}$};
\draw (2-0.5,4-0.5)node{$ x_{32}$};\draw (2-0.5,5-0.5)node{$ y_{33}$};\draw (2-0.5,6-0.5)node{$ x_{24}$};
\draw (3-0.5,4-0.5)node{$ x_{31}$};\draw (3-0.5,5-0.5)node{$ y_{22}$};\draw (3-0.5,6-0.5)node{$ x_{13}$};
\draw (4-0.5,4-0.5)node{$ y_{10}$};
\end{tikzpicture}
}}
\end{array}
\end{equation}

In preparation for the next two Corollaries we adopt a notation whereby
$(\lambda_i,\lambda_j,\mu)=(\lambda_i,\lambda_j,\mu_1,\mu_2,\ldots,\mu_m)$ and 
$(\lambda_i,\mu\backslash\mu_{k})=(\lambda_i,\mu_1,\ldots,\hat{\mu}_{k},\ldots,\mu_m)$ 
where $\hat{\mu}_{k}$ signifies the omission of the part $\mu_{k}$.
Similarly $(\lambda\backslash(\lambda_i,\lambda_j))=(\lambda_1,\ldots,\hat{\lambda}_i,\ldots,\hat{\lambda}_j,\ldots,\lambda_\ell)$
with $\hat{\lambda}_i$ and $\hat{\lambda_j}$ signifying the deletion of the parts $\lambda_i$ and $\lambda_j$,
while $((\lambda\backslash\lambda_i),\mu_{k})=(\lambda_1,\ldots,\hat{\lambda}_i,\ldots,\lambda_\ell,\mu_{k})$.
Such sequences appear as labels $\kappa$ and $\nu$ of various skew $Q$-functions $Q_{\kappa/\nu}(\X,\Y)$ in 
the statement of the two Corollaries.

\begin{Corollary}[Okada type identity]\label{Cor-Qinner}
Let $\lambda$ and $\mu$ be strict partitions of lengths $\ell(\lambda)$ and $\ell(\mu)$, respectively, such that $\mu\subseteq\lambda$. Then
for all $\X=(x_{kc})$ and $\Y=(y_{kc})$
\begin{equation}\label{eqn-Qinner}
    Q_{\lambda/\mu}(\X,\!\Y) = \pf \left( \begin{array}{cc} Q_{(\lambda_i,\lambda_j,\mu)/\mu}(\X,\!\Y)&Q_{(\lambda_i,(\mu\backslash\mu_{k}))/\mu}(\X,\!\Y)\cr
		                                      -\,^t(Q_{(\lambda_i,(\mu\backslash\mu_{k}))/\mu}(\X,\!\Y_))&0\cr \end{array} \right)  
\end{equation} 
with $1\leq i,j\leq\ell$ and $1\leq k\leq m$, where $\ell=\ell(\lambda)$ and $m=\ell(\mu)$ or $\ell(\mu)+1$ according 
as $\ell(\lambda)+\ell(\mu)$ is even or odd, respectively,
with $\mu$ including a trailing $0$ if $\ell(\lambda)+\ell(\mu)$ is odd.
\end{Corollary}

\noindent{\bf Proof}: The appropriate choice of the cutting strip $\phi$ is the inner rim of $SF^{\lambda/\mu}$, 
in which case $\phi=SF^{(\lambda_1,\mu)/\mu}$. Each strip $\theta_r$ of the outside decomposition
takes one or other of the forms $\phi_{0,\lambda_i-1}$ or $\phi_{\mu_k,\lambda_j-1}$, each of which is expressible 
in the form $SF^{\kappa/\mu}$ with shift parameter $0$. 
It follows that in (\ref{eqn-HamXY}) we have 
\begin{equation}
 (\ov{\theta}_{p},\ov{\theta}_{q})=(\phi_{0,\lambda_i-1},\phi_{0,\lambda_j-1})
\quad\mbox{and}\quad
(\theta_{\ell'-k+1}\#\theta_{p})=\phi_{\mu_{k},\lambda_i-1}
\end{equation}
for some $i,j\in\{1,2,\ldots,\ell'\}$ and for $k\in\{1,2,\ldots,m'\}$,
with $\ell'=\ell(\lambda)$ or $\ell(\lambda)+1$ according as $\ell(\lambda)+\ell(\mu)$ is even or odd respectively,
and $m'=\ell(\mu)$.

Adopting the order specified in~\cite{Ham} one finds:
\begin{equation}
   \theta_r\simeq \begin{cases} \phi_{0,\lambda_{\pi_r}-1}&\mbox{for $r=1,2,\ldots,\ell'-m'$};\cr
	                              \phi_{\mu_{\ell'+1-r},\lambda_{\pi_r}-1}&\mbox{for $r=\ell'-m'+1,\ldots,\ell'-1,\ell'$}.\cr
									\end{cases}
\end{equation}
This corresponds to those strips starting from a box at $(d,d)$ of content $0$ on the main diagonal being placed in order of increasing $d$ and 
followed in order of increasing content of their initial boxes by those lying wholly above the main diagonal, including null strips as necessary.
The particular null strip $\theta_{\ell'-m'}\simeq\phi_{0\ov1}$ arises if and only if $\ell(\lambda)+\ell(\mu)$ is odd. 
The determination of the sequence $\pi=(\pi_1,\pi_2,\ldots,\pi_{\ell'})$ is accomplished, very
much as in the unshifted diagram case, by mean of a mergeing procedure and bracketing operation. When applied to our 
running example in the shifted diagram case, $\lambda=(7,6,4,2)$ and $\mu=(4,3)$, 
their merger gives $(7,6,4,4,3,2)$.
This is of the form $(\lambda_1,\lambda_2,\lambda_3,\mu_1,\mu_2,\lambda_4)$, which when bracketed 
and closed with additional $,0)$'s in the form 
$(\lambda_1,(\lambda_2,(\lambda_3,\mu_1),\mu_2),(\lambda_4,0),0)$ yields the sequence of ordered pairs 
$((\lambda_1,0),(\lambda_4,0)),(\lambda_2,\mu_2),(\lambda_3,\mu_1))$ and correspondingly
\begin{equation}
  \theta_1\simeq\phi_{0,\lambda_1-1},~~\theta_2\simeq\phi_{0,\lambda_4-1},~~\theta_3\simeq\phi_{\mu_2,\lambda_2-1},~~\theta_4\simeq\phi_{\mu_1,\lambda_3-1},
\end{equation}

More generally, in order to rearrange the elements of the Pfaffian in (\ref{eqn-HamXY}) so as to
correspond with those of the Pfaffian in (\ref{eqn-Qinner}) one applies the permutation $\pi$
mapping $(1,2,\ldots,\ell')$ to $(\pi_1,\pi_2,\cdots,\pi_{\ell'})$
to the first $\ell'$ rows and columns giving a sign change of $\sgn(\pi)$. 
In the case $\ell'=\ell(\lambda)$, as a result of the merging process that defines $\pi$, we have 
$\sgn(\pi)=(-1)^{\epsilon(\lambda\merge\mu)}$. In the case $\ell'=\ell(\lambda)+1$ the overall sign change remains 
the same. This comes about because the addition of a trailing $0$ to $\lambda$ gives
$\sgn(\pi)=(-1)^{\epsilon((\lambda,0)\merge\mu)}=(-1)^{\epsilon(\lambda\merge\mu)}(-1)^{\ell(\mu)}$.
However, to complete the rearrangement of the elements of the Pfaffian of (\ref{eqn-HamXY}) to coincide 
with those of the Pfaffian of (\ref{eqn-Qinner}) one must apply a further cyclic permutation moving 
the $\ell'$th row and column to the $(\ell'+m)$th position with $m=\ell(\mu)$, thereby introducing an additional sign factor 
$(-1)^{\ell(\mu)}$, leaving once again an overall sign change of $(-1)^{\epsilon(\lambda\merge\mu)}$.

It remains to look at the precise relationship between corresponding elements 
of the Pfaffians of (\ref{eqn-HamXY}) and (\ref{eqn-Qinner}). It should be noted that
for our inner rim cutting strip $\phi$ we have $(\phi_{0,\lambda_i-1},\phi_{0,\lambda_j-1})$ 
of regular skew shifted shape if only if $\lambda_i>\lambda_j$, i.e $i<j$. 
Regularising the elements of the Pfaffian of (\ref{eqn-Qinner}) by way of the merge procedure gives 
$Q_{(\lambda_i,\lambda_j,\mu)/\mu}(\X,\!\Y)=\sgn(\sigma)Q_{(\phi_{0,\lambda_i-1},\phi_{0,\lambda_j-1})}(\X,\!\Y)$.
with $\sgn(\sigma)=(-1)^{\epsilon(\lambda_i\merge\mu)+\epsilon(\lambda_j\merge\mu)}$ for $i<j$ 
and $(-1)^{\epsilon(\lambda_i\merge\mu)+\epsilon(\lambda_j\merge,\mu)+1}$ for $i>j$, where
the extra factor of $-1$ along with the fact that $Q_{(\phi_{0,\lambda_i-1},\phi_{0,\lambda_j-1})}(\X,\!\Y)=0$ if $i=j$,
ensures the required antisymmetry of the Pfaffian matrix.
In addition $Q_{(\lambda_i,(\mu\backslash\mu_{k}))/\mu}(\X,\!\Y)=\sgn(\tau)Q_{(\phi_{\mu_{k},\lambda_i-1})}(\X,\!\Y)$
with $\sgn(\tau)=(-1)^{\epsilon(\lambda_i\merge\mu)}$ if $\lambda_i\geq\mu_k$, with  
$Q_{(\phi_{\mu_{k},\lambda_i-1})}(\X,\!\Y)=0$ if $\lambda_i<\mu_k$.
Thus the regularisation of the terms in the Pfaffian of (\ref{eqn-Qinner}) gives rise to 
sign factors $(-1)^{\epsilon(\lambda_i\merge\mu)}$ and $(-1)^{\epsilon(\lambda_j\merge,\mu)}$
common to each non-zero element of row $i$ and column $j$, respectively. The extraction
of these common factors, identical for corresponding rows and columns, leaves the Pfaffian 
unchanged except for an overall sign factor of $(-1)^{\epsilon(\lambda\merge\mu)}$. 
It is important to note that this remains true whether or not $\mu$ contains a trailing $0$
because $(-1)^{\epsilon(\lambda\merge(\mu,0))}=(-1)^{\epsilon(\lambda\merge\mu)}$. 

This factor obtained by regularisation therefore precisely cancels that obtained by row and column 
permutation, thereby completing the proof of  the Corollary~\ref{Cor-Qinner}.
\qed
 
\begin{Corollary}\label{Cor-Qouter}
Let $\lambda$ and $\mu$ be strict partitions of lengths $\ell(\lambda)$ and $\ell(\mu)$, respectively, such that $\mu\subseteq\lambda$. 
Let 
\begin{equation}
Q_{\lambda/(\lambda\backslash([\lambda_i,\lambda_j]))}(\X,\!\Y)=
 \begin{cases}-Q_{\lambda/(\lambda\backslash(\lambda_i,\lambda_j))}(\X,\!\Y)&\mbox{if $i<j$};\cr
               ~~~0&\mbox{if $i=j$};\cr
							 +Q_{\lambda/(\lambda\backslash(\lambda_i,\lambda_j))}(\X,\!\Y)&\mbox{if $i>j$}. \cr
							\end{cases}
\end{equation}
Then for all $\X=(x_{kc})$ and $\Y=(y_{kc})$
\begin{equation}\label{eqn-Qouter}
    Q_{\lambda/\mu}(\X,\!\Y) =  \pf\left( \begin{array}{cc} 
		                       Q_{\lambda/(\lambda\backslash([\lambda_i,\lambda_j]))}(\X,\!\Y)&Q_{\lambda/((\lambda\backslash\lambda_i),\mu_{k})}(\X,\!\Y)\cr
		                 -\,^t Q_{\lambda/((\lambda\backslash\lambda_i),\mu_{k})}(\X,\!\Y)&0\cr \end{array} \right)  
\end{equation} 
with $1\leq i,j\leq l$ and $1\leq k\leq m$, where $l=\ell(\lambda)$ and $m=\ell(\mu)$ or $\ell(\mu)+1$ according 
as $\ell(\lambda)+\ell(\mu)$ is even or odd, respectively, with $\mu$ including a trailing $0$ if $\ell(\lambda)+\ell(\mu)$ is odd. 

\end{Corollary}

\noindent{\bf Proof}: 
The appropriate choice of the cutting strip $\phi$ is the outer rim of $SF^{\lambda/\mu}$, 
with $\phi=SF^{\lambda/(\lambda\backslash\lambda_1)}$. As in the inner rim case
each strip $\theta_r$ of the outside decomposition takes one or other of the forms $\phi_{0,\lambda_i-1}$ or 
$\phi_{\mu_k,\lambda_j-1}$, each of which is expressible in the form $SF^{\lambda/\nu}$ 
for some partition $\nu$ with shift parameter $0$. 
The argument is then essentially the same as before except that this time
the required antisymmetry of the Pfaffian matrix in (\ref{eqn-Qouter}) has to be imposed through the definition of 
$Q_{\lambda/(\lambda\backslash([\lambda_i,\lambda_j]))}(\X,\!\Y)$. The fact that the accompanying
$+$ sign appears in the case $i>j$ comes about because in the case of an outer rim cutting strip $\phi$ 
it is only $(\phi_{0,\lambda_i-1},\phi_{0,\lambda_j-1})$ with $\lambda_i<\lambda_j$ that is of regular 
skew shifted shape $SF^{\lambda/(\lambda\backslash(\lambda_i,\lambda_j)}$.
It is notable also that in this outer rim case 
$Q_{\lambda/(\lambda\backslash(\lambda_i,\lambda_j))}(\X,\!\Y)=Q_{(\phi_{0,\lambda_i-1},\phi_{0,\lambda_j-1})}(\X,\!\Y)$
if $i>j$ without any necessity of regularisation, and thus no sign factors $\sgn(\sigma)$ appear.
However $Q_{\lambda/((\lambda\backslash\lambda_i),\mu_{k})}(\X,\!\Y)$ does require regularisation
leading to a sign factor $(-1)^{\epsilon(\lambda\merge\mu_k)}$. Extracting these factors this time from the 
last $m$ rows and columns of the Pfaffian in (\ref{eqn-Qouter}) yields an overall sign factor of $(-1)^{\epsilon(\lambda\merge\mu)}$
regardless of whether or not $\mu$ involves a trailing $0$. Once again this overall sign is just what is required to counter the
sign factor of $(-1)^{\epsilon(\lambda\merge\mu)}$ obtained by simultaneously permuting 
the rows and columns of the Pfaffian of (\ref{eqn-HamXY}) where the relevant permutation is exactly 
the same as the one appearing, for the same $\lambda$ and $\mu$, in the inner rim outside
decomposition case.
\qed


To see all this played out in detail, we offer an inner rim example as follows in the case $\lambda=(7,6,4,2)$ and $\mu=(4,3)$.
The content of an inner rim cutting strip $\phi$ and the corresponding outside decomposition take the form:
\begin{equation}\label{eqn-Qphi}
\phi = 
\vcenter{\hbox{

}}
\end{equation}
This allows us to identify
\begin{equation}
 \theta_1\simeq \phi_{06},~~\theta_2\simeq\phi_{01}~~\theta_3\simeq\phi_{35},~~\theta_4\simeq\phi_{43},
\end{equation}
and conclude that
\begin{equation}
(\ov{\theta}_{\red{p}},\ov{\theta}_{\red{q}})= 
\left[
\begin{array}{cccc}
    (\phi_{06},\phi_{06})&(\phi_{06},\phi_{01})&(\phi_{06},\phi_{05})&(\phi_{06},\phi_{03})\cr
		(\phi_{01},\phi_{06})&(\phi_{01},\phi_{01})&(\phi_{01},\phi_{05})&(\phi_{01},\phi_{03})\cr
		(\phi_{05},\phi_{06})&(\phi_{05},\phi_{01})&(\phi_{05},\phi_{05})&(\phi_{05},\phi_{03})\cr
		(\phi_{03},\phi_{06})&(\phi_{03},\phi_{01})&(\phi_{03},\phi_{05})&(\phi_{03},\phi_{03})\cr
		
\end{array}
\right]
\quad
(\theta_k\#\theta_{\red{p}})= 
\left[
\begin{array}{cc}
    \phi_{46}&\phi_{36}\cr
		\phi_{41}&\phi_{31}\cr
		\phi_{45}&\phi_{35}\cr
		\phi_{43}&\phi_{33}\cr
\end{array}
\right]
\end{equation}
Given the shape and content of $\phi$ in (\ref{eqn-Qphi}) these give diagrammatically:
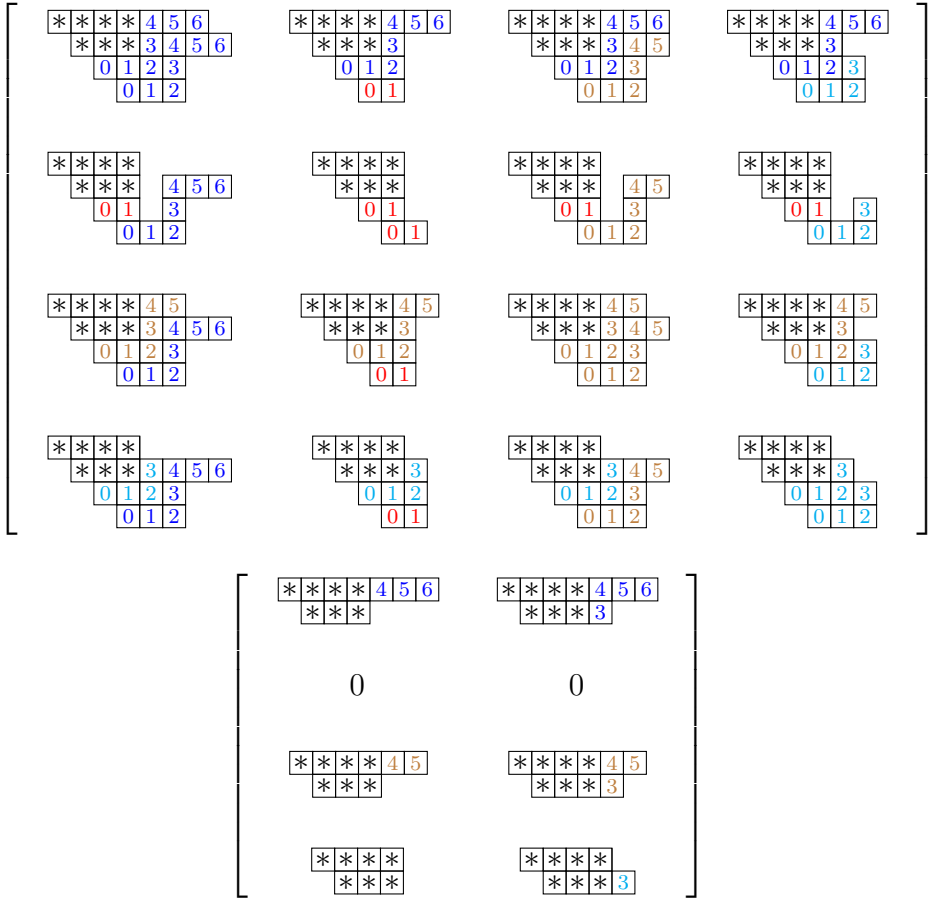
\begin{equation}\nonumber
\left[
\begin{array}{cccc}
\vcenter{\hbox{
\begin{tikzpicture}[x={(0in,-0.12in)},y={(0.12in,0in)}] 
\foreach \j in {1,...,4} \draw (1,\j) rectangle +(-1,-1);
\foreach \j in {2,...,4} \draw (2,\j) rectangle +(-1,-1);
\foreach \j in {1,...,4} \draw(0.5,\j-0.5)node{$\ast$};
\foreach \j in {2,...,4} \draw(1.5,\j-0.5)node{$\ast$};
\draw (1,7) rectangle +(-1,-1);\draw(1-0.5,7-0.5)node{$\sc\blue6$};
\draw (1,6) rectangle +(-1,-1);\draw(1-0.5,6-0.5)node{$\sc\blue5$};
\draw (1,5) rectangle +(-1,-1);\draw(1-0.5,5-0.5)node{$\sc\blue4$};
\draw (2,5) rectangle +(-1,-1);\draw(2-0.5,5-0.5)node{$\sc\blue3$};
\draw (3,5) rectangle +(-1,-1);\draw(3-0.5,5-0.5)node{$\sc\blue2$};
\draw (3,4) rectangle +(-1,-1);\draw(3-0.5,4-0.5)node{$\sc\blue1$};
\draw (3,3) rectangle +(-1,-1);\draw(3-0.5,3-0.5)node{$\sc\blue0$};
%
\draw (1+1,7+1) rectangle +(-1,-1);\draw(1+0.5,7+0.5)node{$\sc\blue6$};
\draw (1+1,6+1) rectangle +(-1,-1);\draw(1+0.5,6+0.5)node{$\sc\blue5$};
\draw (1+1,5+1) rectangle +(-1,-1);\draw(1+0.5,5+0.5)node{$\sc\blue4$};
\draw (2+1,5+1) rectangle +(-1,-1);\draw(2+0.5,5+0.5)node{$\sc\blue3$};
\draw (3+1,5+1) rectangle +(-1,-1);\draw(3+0.5,5+0.5)node{$\sc\blue2$};
\draw (3+1,4+1) rectangle +(-1,-1);\draw(3+0.5,4+0.5)node{$\sc\blue1$};
\draw (3+1,3+1) rectangle +(-1,-1);\draw(3+0.5,3+0.5)node{$\sc\blue0$};
\end{tikzpicture}
}}
&
\vcenter{\hbox{
\begin{tikzpicture}[x={(0in,-0.12in)},y={(0.12in,0in)}] 
\foreach \j in {1,...,4} \draw (1,\j) rectangle +(-1,-1);
\foreach \j in {2,...,4} \draw (2,\j) rectangle +(-1,-1);
\foreach \j in {1,...,4} \draw(0.5,\j-0.5)node{$\ast$};
\foreach \j in {2,...,4} \draw(1.5,\j-0.5)node{$\ast$};
\draw (1,7) rectangle +(-1,-1);\draw(1-0.5,7-0.5)node{$\sc\blue6$};
\draw (1,6) rectangle +(-1,-1);\draw(1-0.5,6-0.5)node{$\sc\blue5$};
\draw (1,5) rectangle +(-1,-1);\draw(1-0.5,5-0.5)node{$\sc\blue4$};
\draw (2,5) rectangle +(-1,-1);\draw(2-0.5,5-0.5)node{$\sc\blue3$};
\draw (3,5) rectangle +(-1,-1);\draw(3-0.5,5-0.5)node{$\sc\blue2$};
\draw (3,4) rectangle +(-1,-1);\draw(3-0.5,4-0.5)node{$\sc\blue1$};
\draw (3,3) rectangle +(-1,-1);\draw(3-0.5,3-0.5)node{$\sc\blue0$};
\draw (3+1,4+1) rectangle +(-1,-1);\draw(3+0.5,4+0.5)node{$\sc\red1$};
\draw (3+1,3+1) rectangle +(-1,-1);\draw(3+0.5,3+0.5)node{$\sc\red0$};
\end{tikzpicture}
}}
&
\vcenter{\hbox{
\begin{tikzpicture}[x={(0in,-0.12in)},y={(0.12in,0in)}] 
\foreach \j in {1,...,4} \draw (1,\j) rectangle +(-1,-1);
\foreach \j in {2,...,4} \draw (2,\j) rectangle +(-1,-1);
\foreach \j in {1,...,4} \draw(0.5,\j-0.5)node{$\ast$};
\foreach \j in {2,...,4} \draw(1.5,\j-0.5)node{$\ast$};
\draw (1,7) rectangle +(-1,-1);\draw(1-0.5,7-0.5)node{$\sc\blue6$};
\draw (1,6) rectangle +(-1,-1);\draw(1-0.5,6-0.5)node{$\sc\blue5$};
\draw (1,5) rectangle +(-1,-1);\draw(1-0.5,5-0.5)node{$\sc\blue4$};
\draw (2,5) rectangle +(-1,-1);\draw(2-0.5,5-0.5)node{$\sc\blue3$};
\draw (3,5) rectangle +(-1,-1);\draw(3-0.5,5-0.5)node{$\sc\blue2$};
\draw (3,4) rectangle +(-1,-1);\draw(3-0.5,4-0.5)node{$\sc\blue1$};
\draw (3,3) rectangle +(-1,-1);\draw(3-0.5,3-0.5)node{$\sc\blue0$};
\draw (1+1,6+1) rectangle +(-1,-1);\draw(1+0.5,6+0.5)node{$\sc\brown5$};
\draw (1+1,5+1) rectangle +(-1,-1);\draw(1+0.5,5+0.5)node{$\sc\brown4$};
\draw (2+1,5+1) rectangle +(-1,-1);\draw(2+0.5,5+0.5)node{$\sc\brown3$};
\draw (3+1,5+1) rectangle +(-1,-1);\draw(3+0.5,5+0.5)node{$\sc\brown2$};
\draw (3+1,4+1) rectangle +(-1,-1);\draw(3+0.5,4+0.5)node{$\sc\brown1$};
\draw (3+1,3+1) rectangle +(-1,-1);\draw(3+0.5,3+0.5)node{$\sc\brown0$};
\end{tikzpicture}
}}
&
\vcenter{\hbox{
\begin{tikzpicture}[x={(0in,-0.12in)},y={(0.12in,0in)}] 
\foreach \j in {1,...,4} \draw (1,\j) rectangle +(-1,-1);
\foreach \j in {2,...,4} \draw (2,\j) rectangle +(-1,-1);
\foreach \j in {1,...,4} \draw(0.5,\j-0.5)node{$\ast$};
\foreach \j in {2,...,4} \draw(1.5,\j-0.5)node{$\ast$};
\draw (1,7) rectangle +(-1,-1);\draw(1-0.5,7-0.5)node{$\sc\blue6$};
\draw (1,6) rectangle +(-1,-1);\draw(1-0.5,6-0.5)node{$\sc\blue5$};
\draw (1,5) rectangle +(-1,-1);\draw(1-0.5,5-0.5)node{$\sc\blue4$};
\draw (2,5) rectangle +(-1,-1);\draw(2-0.5,5-0.5)node{$\sc\blue3$};
\draw (3,5) rectangle +(-1,-1);\draw(3-0.5,5-0.5)node{$\sc\blue2$};
\draw (3,4) rectangle +(-1,-1);\draw(3-0.5,4-0.5)node{$\sc\blue1$};
\draw (3,3) rectangle +(-1,-1);\draw(3-0.5,3-0.5)node{$\sc\blue0$};
\draw (0+1,4+1)--(1+1,4+1);
\draw (2+1,5+1) rectangle +(-1,-1);\draw(2+0.5,5+0.5)node{$\sc\cyan3$};
\draw (3+1,5+1) rectangle +(-1,-1);\draw(3+0.5,5+0.5)node{$\sc\cyan2$};
\draw (3+1,4+1) rectangle +(-1,-1);\draw(3+0.5,4+0.5)node{$\sc\cyan1$};
\draw (3+1,3+1) rectangle +(-1,-1);\draw(3+0.5,3+0.5)node{$\sc\cyan0$};
\end{tikzpicture}
}}
\cr
\end{array}
\right]
\end{equation}

\begin{equation}\nonumber
\left[
\begin{array}{cc}
\vcenter{\hbox{
\begin{tikzpicture}[x={(0in,-0.12in)},y={(0.12in,0in)}] 
\foreach \j in {1,...,4} \draw (1,\j) rectangle +(-1,-1);
\foreach \j in {2,...,4} \draw (2,\j) rectangle +(-1,-1);
\foreach \j in {1,...,4} \draw(0.5,\j-0.5)node{$\ast$};
\foreach \j in {2,...,4} \draw(1.5,\j-0.5)node{$\ast$};
%
\draw (1,7) rectangle +(-1,-1);\draw(1-0.5,7-0.5)node{$\sc\blue6$};
\draw (1,6) rectangle +(-1,-1);\draw(1-0.5,6-0.5)node{$\sc\blue5$};
\draw (1,5) rectangle +(-1,-1);\draw(1-0.5,5-0.5)node{$\sc\blue4$};
\end{tikzpicture}
}}
&
\vcenter{\hbox{
\begin{tikzpicture}[x={(0in,-0.12in)},y={(0.12in,0in)}] 
\foreach \j in {1,...,4} \draw (1,\j) rectangle +(-1,-1);
\foreach \j in {2,...,4} \draw (2,\j) rectangle +(-1,-1);
\foreach \j in {1,...,4} \draw(0.5,\j-0.5)node{$\ast$};
\foreach \j in {2,...,4} \draw(1.5,\j-0.5)node{$\ast$};
%
\draw (1,7) rectangle +(-1,-1);\draw(1-0.5,7-0.5)node{$\sc\blue6$};
\draw (1,6) rectangle +(-1,-1);\draw(1-0.5,6-0.5)node{$\sc\blue5$};
\draw (1,5) rectangle +(-1,-1);\draw(1-0.5,5-0.5)node{$\sc\blue4$};
\draw (2,5) rectangle +(-1,-1);\draw(2-0.5,5-0.5)node{$\sc\blue3$};
\end{tikzpicture}
}}
\cr\cr
0
&
0
\cr\cr
\end{array}
\right]
\end{equation}
Setting $(\X,\!\Y)=(Z)$ for typographical simplicity, these diagrams imply that
\begin{equation}
%
		\right]\cr
\end{array}
\end{equation}
where every term has first been written as a regular $Q$-function, $1$ or $0$ by means of (\ref{eqn-Qpq})
which implies that $Q_{(\ov{\theta}_i,\ov{\theta}_j)}(Z)=-Q_{(\ov{\theta}_j,\ov{\theta}_i)}(Z)$
with $Q_{(\ov{\theta}_i,\ov{\theta}_i)}(Z)=0$, together with the fact that $Q_{\kappa/\mu}(Z)=0$ if 
$\kappa_i=\kappa_j$ for any $i\neq j$, while $Q_{\kappa/\mu}(Z)=1$ if $\kappa=\mu$ and $=0$ if $\mu\not\subseteq\kappa$.

In this case the required simultaneous permutation of rows and columns is the one that maps $(7,2,6,4)$ to $(7,6,4,2)$. 
This has sign $+1$, in agreement with the fact that $\lambda\merge\mu=(\blue{7,6,4,2})\merge(\red{4,3})=(\blue{7,6,4},\red{4,3},\blue{2})$
with $(-1)^{\epsilon(\lambda\merge\mu)}=(-1)^2$. Applying this permutation to the first four rows and columns of
the above display implies that 
\begin{equation}
\begin{array}{l}
Q_{\blue{7642}/\red{43}}(Z)\cr
=\pf\left[
             \begin{array}{cccccc}
   0& Q_{\blue{76}\red{43}/43}(Z)&0&Q_{\blue{7}\red{43}\blue{2}/43}(Z)    &Q_{\blue{7}\red{3}/43}(Z)&Q_{\blue{7}\red{4}/43}(Z)\cr	
	-Q_{\blue{76}\red{43}/43}(Z)&0&0&Q_{\blue{6}\red{43}\blue{2}/43}(Z)     &Q_{\blue{6}\red{3}/43}(Z)&Q_{\blue{6}\red{4}/43}(Z)\cr
	    0& 0          &0     &0          &1        &0\cr
	 -Q_{\blue{7}\red{43}\blue{2}/43}(Z) &-Q_{\blue{6}\red{43}\blue{2}/43}(Z)&0 &0             &0        &0\cr      
			   -Q_{\blue{7}\red{3}/43}(Z)  &-Q_{\blue{6}\red{3}/43}(Z)&-1 &0 &0        &0\cr
				 -Q_{\blue{7}\red{4}/43}(Z)  &-Q_{\blue{6}\red{4}/43}(Z)&0  &0 &0        &0\cr
						 \end{array}
		\right]\cr\cr	
=\pf\left[
             \begin{array}{ccccccc}
 Q_{\blue{77}\red{43}/43}(Z) & Q_{\blue{76}\red{43}/43}(Z)&Q_{\blue{74}\red{43}/43}(Z)&Q_{\blue{72}\red{43}/43}(Z)&Q_{\blue{7}\red{3}/43}(Z)&Q_{\blue{7}\red{4}/43}(Z)\cr			   Q_{\blue{67}\red{43}/43}(Z)&Q_{\blue{66}\red{43}/43}(Z)&Q_{\blue{64}\red{43}/43}(Z)&Q_{\blue{62}\red{43}/43}(Z)&Q_{\blue{6}\red{3}/43}(Z)&Q_{\blue{6}\red{4}/43}(Z)\cr
	Q_{\blue{47}\red{43}/43}(Z)	& Q_{\blue{46}\red{43}/43}(Z)  &Q_{\blue{44}\red{43}/43}(Z)&Q_{\blue{42}\red{43}/43}(Z) &Q_{\blue{4}\red{3}/43}(Z)&Q_{\blue{4}\red{4}/43}(Z)\cr      
Q_{\blue{27}\red{43}/43}(Z)&Q_{\blue{26}\red{43}/43}(Z)&Q_{\blue{24}\red{43}/43}(Z)&Q_{\blue{22}\red{43}/43}(Z)&Q_{\blue{2}\red{3}/43}(Z)&Q_{\blue{2}\red{4}/43}(Z)\cr 
			   -Q_{\blue{7}\red{3}/43}(Z) &-Q_{\blue{6}\red{3}/43}(Z)&-Q_{\blue{4}\red{3}/43}(Z) &-Q_{\blue{2}\red{3}/43}(Z)  &0      &0\cr
				 -Q_{\blue{7}\red{4}/43}(Z) &-Q_{\blue{6}\red{4}/43}(Z)&-Q_{\blue{4}\red{4}/43}(Z) &-Q_{\blue{2}\red{4}/43}(Z) &0 &0 \cr
				\end{array}
				\right]\cr	
\end{array}
\end{equation}
where in the final step the required form (\ref{eqn-Qinner}) has been recovered by rewriting each term, including the $0$'s and $\pm1$'s, appropriately. 
That each replacement is valid may be checked by using the fact that $Q_{\kappa/\mu}(Z)=\sgn(\sigma)Q_{\sigma(\kappa)/\mu}(Z)$ for 
any permutation $\sigma$, together with the previously mentioned special cases, $Q_{\kappa/\mu}(Z)=0$ if 
$\kappa_i=\kappa_j$ for any $i\neq j$, while $Q_{\kappa/\mu}(Z)=1$ if $\kappa=\mu$ and $=0$ if $\mu\not\subseteq\kappa$.

\section{Conclusion}

We have introduced a shifted tableaux definition for the ninth variation $Q$-functions of Okada~\cite{Oka2} and extended the work of Hamel~\cite{Ham} to produce a family of Pfaffians defined by outside decompositions of these shifted tableaux. As one of several corollaries we have obtained combinatorially an explicit Pfaffian identity derived algebraically in~\cite{Oka1}, showing that it is based on an inner rim outside decomposition, and a corresponding explicit formula in the case of an outer rim outside decomposition.  The same approach has been applied to the Schur function determinantal identities, and we take this opportunity to point out that some more recently proposed identities~\cite{Oka1,MM} are subsumed by the Hamel-Goulden outside decomposition identity. 
In particular the so-called $M$-shifted skew Schur function Giambelli identities of Matsuno and Moriyama~\cite{MM}, given a ninth variation extension by Furukawa and Moriyama~\cite{FM}, can be recovered by the simple expedient of replacing the hook shaped cutting strip appropriate to the Giambelli identity of Corollary~\ref{Cor-GX} by one of the same length with corner at position $(M,0)$ rather than $(0,0)$.  We have thus brought under the outside decomposition umbrella a large number of determinantal and Pfaffian identities for symmetric functions as general as ninth variation. We note that we have stopped short of generalizing this to determinants defined by outside decompositions with thick strips work as in \cite{Jin}, but this could be an additional direction for generalization.   We also note that we could use tableaux techniques to consider determinantal identities for ninth variation supersymmetric functions, as have been considered in Bachmann \cite{Bac}.  This is intended to be a 
topic of future work.

\bigskip


\noindent{\bf Acknowledgements}
\label{sec:ack}
This work was supported by the Canadian Tri-Council Research
Support Fund. The first author (AMF) acknowledges the support of a
Discovery Grant from the Natural Sciences and Engineering Research Council of
Canada (NSERC). The second author (RCK) is grateful for the hospitality extended to him 
by Professor Bill Chen at the Center for Applied Mathematics at Tianjin University 
and for the opportunity to pursue this project while visiting him there.

\end{document}